\documentclass[12 pt]{article}
\usepackage{listings}
\usepackage{amsmath}
\usepackage{amssymb}
\usepackage{graphicx}
\usepackage{longtable}
\usepackage{chapterbib}
\usepackage[timestamp,first]{draftcopy}
\usepackage[first]{draftcopy}
\draftcopySetScale{100}
\draftcopyPageX{200}
\draftcopyPageY{100}
\draftcopyBottomX{200}
\draftcopyBottomY{100}
\draftcopySetScaleFactor{0.8}

\newtheorem{theorem}{Theorem}[section]
\newtheorem{remark}{Remark}[section]
\newtheorem{lemma}{Lemma}[section]
\newtheorem{definition}{Definition}[section]

\usepackage[dvips]{epsfig}
\begin{document}
\title{\bf Exactification of Stirling's Approximation for the Logarithm of the Gamma Function}
\author{\bf Victor Kowalenko\\ School of Mathematics and Statistics\\ 
The University of Melbourne\\ Victoria 3010, Australia.}
\vspace{0.5 cm}

\maketitle




\begin{abstract}
Exactification is the process of obtaining exact values of a function from its complete asymptotic expansion. 
This work studies the complete form of Stirling's approximation for the logarithm of the gamma function, which consists 
of standard leading terms plus a remainder term involving an infinite asymptotic series. To obtain values of the function, 
the divergent remainder must be regularized. Two regularization techniques are introduced: Borel summation and Mellin-Barnes (MB) 
regularization. The Borel-summed remainder is found to be composed of an infinite convergent sum of exponential integrals and 
discontinuous logarithmic terms from crossing Stokes sectors and lines, while the MB-regularized remainders possess one MB integral, 
with similar logarithmic terms. Because MB integrals are valid over overlapping domains of convergence, two MB-regularized asymptotic 
forms can often be used to evaluate the logarithm of the gamma function. Although the Borel-summed remainder is truncated,
albeit at very large values of the sum, it is found that all the remainders when combined with (1) the truncated asymptotic series, 
(2) the leading terms of Stirling's approximation and (3) their logarithmic terms yield identical values that agree
with the very high precision results obtained from mathematical software packages. 
\end{abstract}
 \vspace{0.3cm} 

{\it Keywords}: Asymptotic series, Asymptotic form, Borel summation, Complete asymptotic expansion, Discontinuity, Divergent series, 
Domain of convergence, Exactification, Gamma function, Mellin-Barnes regularization, Regularization, Remainder, Stokes discontinuity, 
Stokes line, Stokes phenomenon, Stokes sector, Stirling's approximation

\vspace{0.3 cm}

{\bfseries 2010 Mathematics Subject Classification:} 30B10, 30B30, 30E15, 30E20, 34E05, 34E15, 40A05, 40G10, 40G99, 41A60

\vspace{0.3 cm}
{\it email}: vkowa@unimelb.edu.au
\vspace{5cm}
\lstset{language=C}

 
\newpage
\section{Introduction}
In asymptotics exactification is defined as the process of obtaining the exact values of a function/integral from its 
complete asymptotic expansion and has already been achieved in two notable cases. For those unfamiliar with the concept,
a complete asymptotic expansion is defined as a power series expansion for a function or integral that not only possesses all the 
terms in a dominant asymptotic series, but also all the terms in frequently neglected subdominant or transcendental 
asymptotic series, should they exist. The latter series are said to lie beyond all orders, while the methods and theory behind
them belong to the discipline or field now known as asymptotics beyond all orders or exponential asymptotics \cite{seg91}. One 
outcome of this relatively new field is that it seeks to obtain far more accurate values from the asymptotic expansions for 
functions/integrals than standard Poincare${\acute{\rm e}}$ asymptotics \cite{whi73}. These calculations, which
often yield values that are accurate to more than twenty decimal places, are referred to as hyperasymptotic evaluations or 
hyperasymptotics, for short. Hence exactification represents the extreme of hyperasymptotics.  

In the first successful case of exactification exact values of a particular case of the generalized Euler-Jacobi series, 
viz.\ $S_3(a) =\sum_{n=1}^{\infty} \exp(-a n^3)$, were evaluated from its the complete asymptotic expansion, which was given 
in powers of $a$. Although it had been found earlier in Ref.\ \cite{kow95} that there could be more than one subdominant
series in the complete asymptotic expansion for the generalized Euler-Jacobi series, the complete asymptotic expansion 
for $S_3(a)$ was found to be composed of an infinite dominant algebraic series and another infinite exponentially-decaying 
asymptotic series, whose coefficients resembled those appearing in the asymptotic series for the Airy function ${\rm Ai}(z)$. 
In carrying out the exactification of this complete asymptotic expansion, a range of values for $a$ was considered with the 
calculations performed to astonishing accuracy. This was necessary in order to observe the effect of the subdominant asymptotic 
series, which required in some instances that the analysis be conducted to 65 decimal places as described in Sec.\ 7 of Ref.\ 
\cite{kow95}. 

In the second case \cite{kow002} exact values of Bessel and Hankel functions were calculated from their well-known asymptotic 
expansions given in Ref.\ \cite{wat95}. In this instance there were no subdominant exponential series because the analysis was
restricted to positive real values of the variable. However, unlike Ref.\ \cite{kow95}, different values or levels of truncation 
were applied to the asymptotic series. Whilst the truncated asymptotic series yielded a different value for a fixed 
value of the variable, when it was added to the corresponding regularized value for its remainder, the actual 
value of the Bessel or Hankel function was calculated to within the machine precision of the computing system. The level of 
truncation was governed by an integer parameter $N$, which will also be introduced here. In fact, the truncation parameter as it is
called will play a much greater role here since it will be set equal to much larger values than those in Ref.\ \cite{kow002}. 

Because a complete asymptotic expansion is composed of divergent series, exactification involves being able to obtain 
meaningful values from such series. To evaluate these, one must introduce the concept of regularization, which is defined 
in this work as the removal of the infinity in the remainder of an asymptotic series in order to make the series summable. It was 
first demonstrated in Ref.\ \cite{kow001} that the infinity appearing in the remainder of an asymptotic series arises from
an impropriety in the method used to derive it. Consequently, regularization was seen as a necessary means of correcting an asymptotic 
method such as the method of steepest descent or iteration of a differential equation. Regularization was also shown to be 
analogous to taking the finite or Hadamard part of a divergent integral \cite{kow001}- \cite{kow11c}.

Two very different techniques will be used to regularize the divergent series appearing throughout this work. As described 
in Refs.\ \cite{kow002,kow09}, the most common method of regularizing a divergent series is Borel summation, but often, it 
produces results that are not amenable to fast and accurate computation. To overcome this drawback, the numerical 
technique of Mellin-Barnes regularization was developed for the first time in Ref.\ \cite{kow95}. In this regularization technique 
divergent series are expressed via Cauchy's residue theorem in terms of Mellin-Barnes integrals and divergent arc-contour integrals. 
In the process of regularization the latter integrals are discarded, while the Mellin-Barnes integrals yield finite 
values, again much like the Hadamard finite part of a divergent integral. Amazingly, the finite values obtained when this technique 
is applied to an asymptotic expansion of a function yield exact values of the original function, but with one major difference
compared with Borel summation. Instead of having to deal with Stokes sectors and lines, we now have to contend with the domains of 
convergence for the Mellin-Barnes integrals, which not only encompass the former, but also overlap each other. So, while 
Borel summation and Mellin-Barnes regularization represent techniques for regularizing asymptotic series and yield the same values  
for the original function from which the complete asymptotic expansion has been derived, they are nevertheless completely different. 
Moreover, they can be used as a check on one another, which will occur throughout this work.

In the two cases of exactification mentioned above only positive real values of the power variable in the asymptotic expansions
were considered, although it was stated that complex values would be studied in the future. As discussed in the preface to Ref.\ 
\cite{kow95}, such an undertaking represents a formidable challenge because as the variable in an asymptotic series moves about the 
complex plane or its argument changes, a complete asymptotic expansion experiences significant modification due to the Stokes 
phenomenon \cite{sto04}. This means that at particular rays or lines in the complex plane, an asymptotic expansion develops jump 
discontinuities, which can result in the emergence of an extra asymptotic series in the complete asymptotic expansion. Thus, a 
complete asymptotic expansion is only uniform over either a sector or a ray in the complex plane, which means in turn that
in order to exactify a complete asymptotic expansion over all arguments or phases of the variable in the power series, one requires a 
deep understanding of the Stokes phenomenon. This understanding entails: (1) being able to determine the locations of all jump 
discontinuities, and (2) solving the more intricate problem of their quantification when they do occur.

Because of the Stokes phenomenon, it becomes necessary not only to specify a complete asymptotic expansion, but also the 
range of the argument of the variable in the power series of the expansion. The combination of two such statements are referred 
to as asymptotic forms in this work. In particular, it should be noted that if the same complete asymptotic expansion is
valid for different sectors or rays of the complex plane, then in each instance the original function being studied will also 
be different. 

A major advance in enabling asymptotic forms to be evaluated over all values of the argument of the power series variable occurred 
with the publication of Ref.\ \cite{kow09}, which began by building upon Stokes' seminal discovery of the phenomenon now named after 
him \cite{sto04} and then proceeded to develop via a series of propositions a theory/approach that enables complete asymptotic 
expansions to be derived for higher order Stokes sectors than ever considered before. The reason why higher order Stokes sectors or
all values of the argument of the power series variable need to be considered is that when an asymptotic expansion is derived, it is 
often multivalued in nature. If one wishes to evaluate the original function over the entire principal branch for its variable via 
asymptotic forms, then the asymptotic forms pertaining to the higher/lower Stokes sectors may be required due to the fact that the asymptotic 
series are often composed of (inverse) powers of the variable to a power. E.g., the asymptotic series for the error function, ${\rm erf}(z)$, 
which is studied extensively in Ch.\ 1 of Ref.\ \cite{din73}, is composed of an infinite series in powers of $1/z^2$. If one wishes to 
determine values of ${\rm erf}(z)$ over the principal branch via its asymptotic forms, then one requires the asymptotic forms for higher 
Stokes sectors, not just the asymptotic form for the lowest Stokes sector given by $|{\rm arg}\,z | < \pi/2$. This issue will be discussed 
in detail later, particularly in Secs.\ 3 and 5. 

The developments in Ref.\ \cite{kow09} were primarily concerned with the regularization of the two types of generalized 
terminants for all Stokes sectors and lines. The term terminant was introduced by Dingle \cite{din73} after he found that 
the late terms of the asymptotic series for a host of functions in mathematical physics could be approximated by them. In the
present work the aim is to continue with the development of a general theory in asymptotics by using the results in Ref.\ 
\cite{kow09} as a base. It should be borne in mind that since complete asymptotic expansions are composed of divergent series, 
we are essentially talking about the development of a general theory for handling all divergent series, a quest that has stretched 
over centuries, but to this day remains elusive. One important reason why such a theory remains elusive is that the divergence 
can be different for diverse problems requiring different approaches or methods to regularize them. In fact, we 
shall see that the asymptotic forms derived in this work are composed of an infinite power series that can be regularized via 
generalized terminants and another infinite series, which is logarithmically divergent and thus, needs to be handled differently. 

Discovered in the 1730's \cite{wiki12} Stirling's approximation/formula is a famous result for obtaining values of the factorial 
function or its more general version, the gamma function, denoted by $\Gamma(z)$. Because the gamma function exhibits 
rapid exponential growth, those working in asymptotics frequently study the alternative version of the approximation, where it 
is expressed in terms of its logarithm, i.e.\ $\ln \Gamma(z)$. From an asymptotics point of view, $\ln \Gamma(z)$ is a far more 
formidable function than $\Gamma(z)$, because it possesses the multivaluedness of the logarithm function. To compensate for
taking the logarithm, it is often assumed that one has at their disposal an exponentiating routine that will render values of 
$\Gamma(z)$ from the complex numerical values calculated for $\ln \Gamma(z)$. Despite this, however, no one has ever been able to 
render exact values of $\Gamma(z)$ or $\ln \Gamma(z)$ when either function has been expressed in terms of the Stirling approximation 
because when all the terms in the approximation are considered, it becomes an asymptotic expansion. Consequently, it is not 
absolutely convergent. Here we aim to investigate how the results in Ref.\ \cite{kow09} can be used to determine exact values of 
the complete version of Stirling's approximation of $\ln \Gamma(z)$ for all values of ${\rm arg}\, z$. In short, this paper aims to exactify 
Stirling's approximation.   

Sec.\ 2 introduces the standard form of Stirling's approximation for $\ln \Gamma(z)$ as given in Ref.\ \cite{abr65}. Then $\ln \Gamma(z)$ 
is expressed in terms of the specific leading order terms associated with the approximation and a truncated asymptotic power series, 
whose coefficients are related to the Bernoulli numbers. Although it can be divergent, the latter series is often neglected in 
accordance with standard Poincar$\acute{\rm e}$ asymptotics \cite{whi73}. Moreover, it is not known whether inclusion of the 
asymptotic series represents a complete asymptotic expansion for $\ln \Gamma(z)$. Consequently, we turn to Binet's second expression
of $\ln \Gamma(z)$ to generate the complete asymptotic expansion for $\ln \Gamma(z)$, in which the asymptotic series is now expressed
as an infinite sum of generalized terminants. Moreover, a general theory for deriving regularized
values of the remainders for both types of these series via Borel summation is elucidated in Ref.\ \cite{kow09}. With the aid of 
this theory, new asymptotic forms for $\ln \Gamma(z)$ are presented and proved for all Stokes sectors and lines. Because of an
infinite number of singularities situated on each Stokes line, on each occasion when crossing from one Stokes sector to its 
neighbouring sector an extra infinite series appears in the asymptotic forms, which is referred to as the Stokes discontinuity
term. Since these series can become logarithmically divergent, they must be regularized according to Lemma\ 2.2, which describes
the regularization of the standard Taylor series expansion of $\ln(1+z)$. The section concludes with how the expression for
$\ln \Gamma(z)$ can be used to derive new results for the digamma function and Euler's constant.

For the interested reader, it ought to be stated here that there are numerous power series expansions for $\ln \Gamma(z)$ of which the
Stirling approximation is the oldest and arguably, the most famous. Blagouchine \cite{bla016} presents an extensive list
of these expansions together with references to their origin. Furthermore, he derives two new series expansions for this important
special function, which cannot be written explicitly in powers of $z$ up to a given order even though from his extensive numerical
analysis they appear to be converging, albeit very slowly. Perhaps, the most interesting property of the new expansions is they are composed 
of rational coefficients for $z =n/2 \pm \alpha /\pi$, when $n$ is an integer and $\alpha$ is a positive rational greater than $\pi/6$.
This cannot be said of the Stirling approximation since there is a term of $\ln 2 \pi$ in it. With regard to the issue of convergence, however, 
this paper aims to show how the divergence in the Stirling approximation can be tamed so that the final result will yield exact values
of $\ln \Gamma(z)$ in a time-expedient manner for any value of $z$ as indicated earlier. 

Because such symbols as $\sim$, $\simeq$, $\pm \cdot$, $\geq$, $\leq$ and the Landau gauge symbols of $O()$ and $o()$ abound in
the discipline, standard Poincar${\acute{\rm e}}$ asymptotics suffers from the drawbacks of vagueness and limited range of applicability
as discussed in the prologue to Ref.\ \cite{din73}. E.g., when an asymptotic expansion is described as being valid for large and 
small values of a variable, invariably one does not know the actual value when this applies or indeed, just how accurate it is compared 
with the original function. Moreover, this relative accuracy varies for different values of the variable. The afore-mentioned symbols 
will not be employed here at any stage, having become redundant due to the concept of regularization. Nevertheless, the reader may still feel 
uncomfortable or uncertain about the results given in Sec.\ 2. Therefore, Sec.\ 3 presents an extensive numerical investigation aimed at 
verifying, in particular, the results of Theorem\ 2.1. As will be observed, since the results from these studies are exact, they are far 
more accurate than the alternative hyperasymptotic approach of developing strategies for truncating asymptotic series beyond the optimal 
point of truncation as discussed in Refs.\ \cite{ber90}-\cite{ber91a}. In addition, numerical studies can reveal many interesting properties, 
provided they are carried out in an appropriate manner. That is, a non-appropriate study represents situations where claims are often made 
about a hyperasymptotic approach improving the accuracy of an asymptotic expansion when the asymptotic expansion is already very accurate. 
Specifically, in such studies a relatively large value of the variable is chosen when the expansion has already been derived for large values,
instead of choosing small values. Although the hyperasymptotic approach can improve the situation for such large values, the method will still 
break down for small values, which is often disregarded by the proponents of these approaches. In this work, however, the opposite will be done; 
exact values of a large variable expansion will be obtained for small and intermediate values of the variable, which represent the regions 
where an asymptotic expansion is no longer valid according to standard Poincar${\acute{\rm e}}$ asymptotics. In addition, it should be 
pointed out that numbers do not lie, whereas ``theorems" employing the afore-mentioned symbols are unable to indicate just how accurate an 
expansion/approximation is and where it actually does break down. This is because the various terms in a complete asymptotic expansion change 
their behaviour over the complex plane. Consequently, usually neglected subdominant terms can become the dominant contribution.

Sec.\ 3 proceeds with a numerical investigation into the Borel-summed asymptotic forms given by Eq.\ (\ref{fiftyeight}). In the first
example the remainder is expressed as an infinite sum involving the incomplete gamma function, which means that Mathematica's intrinsic
routine is used to evaluate the remainder. To make contact with standard Poincar${\acute{\rm e}}$ asymptotics, initially a large value of 
$|z|$ is chosen, viz. $z \!=\! 3$. The truncation parameter $N$ takes on various values between 1 and 50 with the optimal point of 
truncation occurring around $N_{OP} \!=\!10$. This means for $N \!<\! 15$, the leading terms in Stirling's approximation or $F(z)$ 
given by Eq.\ (\ref{fiftynine}) dominate, but for $N \!>\! 30$, the remainder and truncated sum dominate the calculations. Nevertheless, 
when all the quantities in Eq.\ (\ref{fiftyeight}) are added together, they always yield the exact value of $\ln \Gamma(3)$ regardless 
of the value of $N$, provided a sufficient number of decimal places has been specified in order to allow the cancellation of decimal 
places to occur when the truncated sum and the regularized value of its remainder are combined. The same analysis is repeated for 
$z \!=\! 1/10$, which could never be studied in standard Poincar${\acute {\rm e}}$ asymptotics. For such a value the leading terms in 
Stirling's approximation or $F(z)$ are not accurate. Hence it is vital that the truncated sum of the asymptotic series and the 
regularized value of its remainder must be included to obtain $\ln\Gamma(1/10)$. In fact, because there is no optimal point of truncation, 
the remainder and truncated sum diverge far more rapidly even for relatively small values of the truncation parameter.

At this stage the Borel-summed asymptotic forms have been shown to yield $\ln \Gamma(z)$ for real values of $z$. That is, exactification
of the asymptotic form for a function has been achieved for real values of $z$ as in the previous instances studied in Refs.\
\cite{kow95} and \cite{kow002}. Now, the analysis is extended to complex values of $z$. This is done by putting $z=|z| \exp(i \theta)$,
where $|z|$ equals either 3 or 1/10 and $\theta$ ranges over the principal branch of the complex plane, i.e.\ $-\pi< \theta \leq \pi$.
In addition, the regularized value of the remainder is evaluated by using the Borel-summed forms in Thm.\ 2.1. That is, instead of relying
on a routine that can calculate values of the incomplete gamma function, we evaluate the exponential integrals appearing in Eqs.\
(\ref{sixty}) and (\ref{sixtya}) by using the NIntegrate routine in Mathematica, which has the effect of increasing the time of computation
substantially. For $|z|=3$ and $N$ below the optimal point of truncation, it is found that the leading terms of the Stirling's approximation 
are accurate to the first couple of decimal places to $\ln \Gamma(z)$, but not the thirty figure accuracy sought throughout this work. For 
$\theta >\pi/2$ and $\theta<-\pi/2$, which represent the adjacent Stokes sectors to $|\theta|< \pi/2$, the Stokes discontinuity is non-zero, 
but is $O(10^{-7})$. Therefore, it only affects the accuracy beyond the seventh decimal place. Similarly, the remainder is very small. 

The $|z|=1/10$ situation, however, is completely different. Except for the values of the truncation parameter $N$ close to zero, the 
regularized value of the remainder and the truncated sum dominate, although when combined, they cancel many decimal places. 
For $|\theta|<\pi/2$, the Stokes discontinuity term vanishes as for the $|z|=3$ case, but outside the sector, it makes a substantial 
contribution to $\ln \Gamma(z)$, which cannot be neglected. Therefore, whilst the $|z|=1/10$ case exhibits completely different behaviour
to $|z|=3$, when all the contributions are combined, it is able to yield the values of $\ln \Gamma(z)$ to the 30 decimal places specified
by the Mathematica module/program as in the $|z|=3$ case.

For $z$ situated on the Stokes lines of $\theta=\pm \pi/2$, the asymptotic forms possess Cauchy principal value integrals for
their regularized remainder in addition to a Stokes discontinuity term that is calculated by summing semi-residue contributions
rather than full residue contributions in the adjacent Stokes sectors. Consequently, a different program is required to evaluate
the various terms belonging to $\ln \Gamma(z)$. For the Stokes line given by $\theta=\pi/2$, $|z|$ is set equal to 3, while for $\theta =-\pi/2$,
it is set equal to 1/10. In both cases the truncated series and regularized remainder are imaginary, while the Stokes discontinuity
term is real, which is consistent with the rules for the Stokes phenomenon presented in Ch.\ 1 of Ref.\ \cite{din73}. As found
previously, the truncated sum and regularized value of the remainder are very small when $N<20$ for the large value of $|z|$, 
but begin to diverge from thereon, while for $|z|=1/10$, they dominate for $N>3$. In addition, for $|z|=3$, the Stokes discontinuity is
very small ($O(10^{-9}$), while for $|z|=1/10$, it represents a sizeable contribution to the real part of $\ln \Gamma(-i/10)$, although
it is still not as significant as the real of the leading terms in the Stirling approximation. Nevertheless, all terms are necessary in
order to obtain $\ln \Gamma(z)$ to thirty decimal figures/places.  

The numerical studies of Sec.\ 3 demonstrate that it is indeed possible to obtain exact values of a function from its asymptotic
expansion based on the conventional view of the Stokes phenomenon. However, this contradicts a more radical view, first proposed by Berry 
\cite{ber89} and later made ``rigorous" by Olver \cite{olv90}. According to this view, which is now called Stokes smoothing in spite of
the fact Stokes never held such a view, the Stokes phenomenon is no longer believed to be discontinuous, but undergoes a smooth and very 
rapid transition at Stokes lines. That is, instead of a step-function multiplying the subdominant terms in an asymptotic expansion as 
in Eqs.\ (\ref{sixtyfour}) and (\ref{sixtysix}), the Stokes multiplier is now reckoned to behave for large values of $|z|$ 
as an error function that transitions rapidly from 0 to unity as opposed to toggling at a Stokes line. The issue requires investigation 
here because if it is indeed valid, then it implies that one can never obtain exact values for a function from an asymptotic 
expansion as the Stokes multiplier cannot be made exact according to Olver's treatment. Whilst an exact representation for the Stokes 
multiplier does not exist, the situation can be inverted by stating that one should not be able to obtain exact values for $\ln \Gamma(z)$ 
in the vicinity of any Stokes line because according to the conventional view originated by Stokes, the multiplier behaves as a 
step-function there. In fact, close to a Stokes line, smoothing implies that the Stokes multiplier is almost equal to 1/2, not zero just 
before it or unity just past it. Therefore, one should not be able to obtain exact values for $\ln \Gamma(z)$ according to the conventional
view if Stokes smoothing is valid. Since a Stokes line occurs at ${\rm arg}\,z \!=\! \pi/2$ in the Borel-summed asymptotic forms of Thm.\ 2.1, 
the numerical investigation in Sec.\ 3 is concluded by evaluating $\ln \Gamma(z)$ for $z \!=\! |z| \exp(i(1/2 +\delta)\pi)$, where $|z|$ 
is large and $\delta$ is both positive and negative with its magnitude ranging from 1/10 to $1/20\,000$. If the smoothing concept is 
correct, then the Stokes multiplier should be close to $1/2$ and the Stokes discontinuity term will be close to its value at $\theta
\!=\! \pi/2$. However, by using the asymptotic forms above and below the Stokes line, viz.\ the first and third forms in Eq.\ (\ref{fiftyeight}), 
we obtain exact values of $\ln \Gamma(z)$ for all values of $\delta$ in Table\ \ref{tab5}, thereby confirming the step-function behaviour 
of the Stokes multiplier close to a Stokes line. Consequently, the conventional view of the Stokes phenomenon is vindicated.

According to standard Poincar${\acute{\rm e}}$ asymptotics \cite{whi73}, it is generally not permissible to differentiate an asymptotic
expansion. Since it has been shown that the asymptotic forms in Sec.\ 2 yield exact values of $\ln \Gamma(z)$ over all values of the
argument of $z$ and $\ln \Gamma(z)$ is differentiable, we can differentiate the asymptotic forms,n Thm.\ 2.1 thereby obtaining asymptotic 
forms for the digamma function, $\psi(z)$. Sec.\ 3 concludes with a theorem, which gives the asymptotic forms for this special function. 

Whilst Secs.\ 2 and 3 are devoted to the derivation and verification of Borel-summed asymptotic forms for $\ln \Gamma(z)$, Sec.\ 4
presents the asymptotic forms for $\ln \Gamma(z)$ obtained via MB regularization. As mentioned earlier, MB regularization 
was introduced in Ref.\ \cite{kow95} because it was able to produce asymptotic forms that are more amenable and faster 
to compute than the corresponding forms obtained via Borel summation. However, this is not the only reason for studying the forms 
obtained by this technique. Suppose the solution $f(z)$ to a problem happens to equal $\ln \Gamma(z^n)$, where $n>2$. In this situation 
we cannot use the inherent routine in a software package such as Mathematica \cite{wol92} to calculate the values of the function 
when $|{\rm arg}\,z| > \pi/n$, since the routine for $\ln \Gamma(z)$ (LogGamma in Mathematica) is restricted to the argument lying
in the principal branch of complex plane or to $-\pi <{\rm arg}\,z \leq \pi$. When $z$ is replaced by $z^n$, the routine will only
provide values of $\ln \Gamma(z^n)$ for $-\pi/n< {\rm arg}\,z \leq \pi/n$. However, the aim is to determine $f(z)$ over the entire
the principal branch. We can use the Borel-summed asymptotic forms for the higher or lower Stokes sectors in this situation, but we 
also need to be able to verify the results. This can be accomplished by using the MB-regularized asymptotic forms, which represent
a different method of evaluating $f(z)$.

Thm.\ 4.1 presents the MB-regularized asymptotic forms for $\ln \Gamma(z)$, where the regularized value for the remainder of the
asymptotic series in the Stirling approximation is expressed in terms of MB integrals with domains of convergence given by 
$(\pm M-1)\pi \!<\! {\rm arg}\,z \!<\! (\pm M+1)\pi$ and $M$, a non-negative integer. Although there are no Stokes lines and sectors 
in MB regularization, there is an extra term to the MB-regularized remainder arising from the summation of the residues of the 
singularities on the Stokes lines. As a consequence, the extra term is expressed as an infinite series, which can become 
logarithmically divergent. Regularization of the series results in ambiguity whenever ${\rm arg}\,z \!=\! (\pm M \!-\! 1/2) \pi$ or at 
the Stokes lines, of which there are two in each domain of convergence. This ambiguity indicates that $\ln \Gamma(z)$ is discontinuous 
at the Stokes lines and is resolved by adopting the Zwaan-Dingle principle \cite{kow09}, which states that an initially real function 
cannot suddenly become imaginary. Hence the real value of the regularized value of the infinite series is chosen to yield the value 
of $\ln \Gamma(z)$ at Stokes lines.

As indicated previously, the LogGamma routine only evaluates $\ln \Gamma(z)$ for $z$ situated in the principal branch. Yet the Borel-summed 
and MB-regularized asymptotic forms presented in Secs.\ 2 and 4 cover all values of ${\rm arg}\,z$. This is particularly useful if we 
wish to evaluate $\ln \Gamma(z^3)$ for $z$ lying in the principal branch despite the fact that the LogGamma routine will only yield values
for $-\pi/3 \!<\! {\rm arg}\,z \!\leq\! \pi/3$. Sec.\ 5 begins by presenting the MB-regularized asymptotic forms for $\ln \Gamma(z^3)$. In this
instance there are seven domains of convergence spanning the principal branch. In addition, the specific forms corresponding to
the Stokes lines are also presented. Then the MB-regularized asymptotic forms are implemented in a Mathematica module and the code is
run for various values of the truncation parameter and $\theta$ ranging between $(-\pi,\pi)$. In those cases where there are two asymptotic
forms for the value of $\theta$, the regularized values of $\ln \Gamma(z^3)$ are found to agree with each other in accordance with the
concept of regularization. Moreover, at the Stokes lines other than $|\theta|=\pi/6$, it is also found that the MB-regularized asymptotic 
forms exhibit jump discontinuities. Basically, the mid-point between the discontinuous values from the left and right sides at each Stokes
line is selected as the regularized value for $\ln \Gamma(z^3)$, which is in accordance with the Zwaan-Dingle principle and yields
logarithmic terms.

Sec.\ 6 proceeds by presenting the Borel-summed asymptotic forms for $\ln \Gamma(z^3)$, where from Eq.\ (\ref{ninety})
there are seven Stokes sectors spanning the principal branch for $z$. This is the same number of domains of convergence for the MB-regularized
asymptotic forms, but since there is no overlapping, the Stokes sectors are much narrower with different endpoints. Nevertheless, in
accordance with the concept of regularization the Borel-summed forms should give identical values to the MB-regularized asymptotic forms.
Because the remainder for the Borel-summed asymptotic forms is an infinite sum of exponential integrals, as indicated previously, the Mathematica
code calculating both MB-regularized and Borel-summed asymptotic forms was implemented as a batch program to enable many calculations to
be performed simultaneously. In addition, the Stokes discontinuity and logarithmic terms require careful treatment. In order to avoid
a blow-out in the computation time and a substantial reduction in the accuracy, a relatively large value of $|z|$, viz.\ 5/2, was chosen. It is
found that despite varying both $\theta$ and the truncation parameter, the MB-regularized results agree with each other to the 30 significant
figures. The Borel-summed results are also in agreement with the MB-regularized asymptotic forms provided the remainder is very small by
choosing the truncation parameter not to be far away from the optimal point of truncation. For $N \!=\! 2$, however, it is found that the Borel-summed
and MB-regularized do agree with each other but to much less significant figures (22 as opposed to 30). This is attributed to the fact that
the Borel-summed remainder represents an approximation because the  upper limit in its sum was set to $10^5$. 

Because the Borel-summed and MB-regularized asymptotic forms for $\ln \Gamma(z^3)$ are distinctly different at the Stokes lines, Sec.\ 6
presents another batch program that evaluates these forms for $|z|=9/10$, which is considered to lie in the intermediate region between
large and small values of $|z|$. Such values could also never be considered in standard Poincar${\acute {\rm e}}$ asymptotics. The major differences
between the MB-regularized and Borel-summed asymptotic forms occur in the regularized values of the remainder of the asymptotic series and
in the multiplier of the main logarithmic term. The remainder in the Borel-summed asymptotic form is composed of an infinite number of 
Cauchy principal value integrals, while the multiplier of the main logarithmic is multiplied by 1/2 signifying that semi-residues have been
evaluated as opposed to full residues in the MB-regularized asymptotic forms. Nevertheless, it is found that the complex values of $\ln \Gamma(z^3)$
obtained from the Borel-summed forms agree to 27 decimal places/figures with the two corresponding MB-regularized forms for each Stokes line.
Moreover, for $\theta \!=\! \pi/6$, all three forms agree with the value obtained via the LogGamma routine in Mathematica \cite{wol92}. With 
this final numerical example Stirling's approximation has been exactified for all values of ${\rm arg}\, z$.

Finally, the various Mathematica programs used in the numerical studies are presented in a condensed format in the appendix. Although the main 
reason for presenting them is to enable the reader to verify the results displayed in the tables, the reader can also use them to conduct their
own numerical studies. In addition, they can be adapted to become more or less accurate by changing the various Options appearing in the main 
routines. In general, it was found that the MB-regularized asymptotic forms took substantially less time to compute than the Borel-summed  
asymptotic forms. 

\section{Stirling's Approximation}
Stirling's approximation or formula \cite{wiki12} is used to approximate large values of the factorial function, although what constitutes
a large number is often unclear. Nevertheless, it is often written as 
\begin{eqnarray}
\ln n! =\ln \Gamma(n+1) = n\ln n - n + \frac{1}{2}\,\ln (2\pi n) + \cdots \;.
\label{one}\end{eqnarray}
Whilst the above statement may be regarded as a good approximation for large integer values of $n$ in standard Poincar${\acute {\rm e}}$ 
asymptotics, e.g., for $n>5$, the difference between the actual value and that from the above result is less than one percent, it is unsuitable 
for hyperasymptotic evaluation since there is an infinite number of terms that have been neglected. In this work we wish to extend the
above approximation for the factorial function to the gamma function, $\Gamma(z)$, where the argument $z$ is complex and its magnitude $|z|$
is not necessarily large as in the above result. For the purposes of this work, the terms on the rhs of the above result will be referred to 
as the leading terms in Stirling's approximation. When we derive the complete asymptotic expansion of $\ln \Gamma(z)$ in Sec.\ 2, 
they will be denoted by $F(z)$, where $n$ is replaced by $z$ in the above result. Hence they will represent a specific contribution when 
calculating exact values of $\ln \Gamma(z)$, which can be checked with the actual value throughout this work to gauge the accuracy of the 
above approximation. 

It should also be noted that the above result also obscures the fact that the missing terms belong to an infinite power series in $1/n^2$ or 
$1/z^2$, which, as we shall see shortly, can become divergent. Even when the series is not divergent, it is conditionally convergent, but never 
absolutely convergent. Nevertheless, the missing terms will be critical for obtaining exact values of $\ln \Gamma(z)$.  

Occasionally, a problem arises where there is an interest in the missing terms in the above statement. Then Stirling's approximation 
is expressed differently. For example, according to No.\ 6.1.41 in Abramowitz and Stegun \cite{abr65}, $\ln \Gamma(z)$ can be expressed as
\begin{eqnarray}
\ln \, \Gamma(z) \sim \Bigl( z - \frac{1}{2} \Bigr) \ln\,z - z + \frac{1}{2} \ln(2\pi)  + \frac{1}{12z}- \frac{1}{360z^3} +
\frac{1}{1260z^5} -\frac{1}{1680z^7} + \dots \;,
\label{two}\end{eqnarray}
where $z \to \infty$ and $|{\rm arg}\,z | < \pi$. Here, we see that the leading terms are identical to the first statement or version when $z$ 
is replaced by $n$. In other textbooks the dots in the above result are replaced by the Landau gauge symbol, which would be $O(z^{-9})$ 
since it is next highest order term that has been omitted. Gradshteyn and Ryzhik \cite{gra94} express the power series after the $\ln(2 \pi)$ 
term as a truncated power series in which the coefficients depend on the Bernoulli numbers. As a consequence, the tilde is replaced by an 
equals sign and a remainder term, $R_N(z)$, is introduced. This is given by
\begin{eqnarray}
R_N(z)= \sum_{k=N}^{\infty} \frac{B_{2k}}{2k(2k-1) z^{2k-1}} \;. 
\label{three}\end{eqnarray}
Although the remainder is bounded according to No.\ 8.344 in Ref.\ \cite{gra94} in terms of $z$ and $N$, for $\Re\, z>0$, the series still 
diverges once the optimal point of truncation is exceeded. Furthermore, Gradshteyn and Ryzhik are even more vague in their presentation of 
Stirling's approximation than Abramowitz and Stegun because they stipulate that the expansion is valid for large values of $|z|$ without specifying 
what large means.  All the preceding forms/material only serve to re-inforce just how vague and confusing standard Poincar${\acute {\rm e}}$ 
asymptotics \cite{whi73} can be, which is somewhat of a contradiction in a supposedly precise subject like mathematics. As indicated in the
introduction, this work aims to investigate whether the vagueness and concomitant limited range of applicability in Stirling's approximation 
can be overcome by employing the concepts and techniques developed in asymptotics beyond all orders \cite{seg91} over recent years, a field that 
has been primarily concerned with obtaining hyperasymptotic values from the asymptotic expansions of functions/integrals. In particular, we aim to 
apply the developments in Refs. \cite{kow95},\cite{kow002}, \cite{kow001}-\cite{kow11c} and \cite{kow11}, where the exact values have been calculated 
from the complete asymptotic expansions of various functions, to the complete version of Stirling's approximation. Before this can be accomplished, 
however, we need to derive the complete form of Stirling's approximation, which means in turn that we require the following lemma.

\begin{lemma}
As a result of regularization the power/Taylor series expansion for ${\rm arctan}\,u$, which is given by $\sum_{k=0}^{\infty} u^{2k+1}/(2k+1)$, 
can be expressed as 
\begin{eqnarray}
\sum_{k=0}^{\infty} \frac{(-1)^k u^{2k+1}}{(2k+1)} \begin{cases} = {\rm arctan}\, u\quad, & \quad -1< \Re \, (iu) < 1 \quad, \cr
\equiv {\rm arctan}\, u \quad, & \quad \Re \, (iu) \leq -1, \quad {\rm and} \quad \Re \, (i u) \geq 1 \quad.  
\end{cases}
\label{four}\end{eqnarray}
\end{lemma}
{\bfseries Proof}. The definition for ${\rm arctan} \, u$ can be obtained from No.\ 2.141(2) in Ref.\ \cite{gra94}, which is
\begin{eqnarray}
{\rm arctan} \,u = \int_0^{u} \frac{dt}{1+ t^2} \;.
\label{five}\end{eqnarray}
Decomposing the integral into partial factions yields
\begin{eqnarray}
{\rm arctan} \, u = \frac{1}{2} \int_0^{u} dt \; \Bigl( \frac{1}{1-it} + \frac{1}{1+it} \Bigr) \;.
\label{six}\end{eqnarray}
We now replace the integrands on the rhs of the above integral by their respective form of the geometric series. According to Refs.\ 
\cite{kow002}, \cite{kow09}, \cite{kow11} and \cite{kow09a}, the first geometric series is convergent when $\Re\, (it) <1$ for all values 
of $t$ lying in $[0,u]$, while the second geometric series is convergent when $\Re\, (it) > -1$ for all values of $t$ lying in $[0,u]$.
Therefore, in the strip given by $\Re \, (iu)<1$ and $\Re\, (iu)>-1$, both geometric series will be convergent. Moreover, according to 
the same references, both geometric series will be absolutely convergent in the strip only when $u$ is situated inside the unit disk, 
i.e. $|u|<1$. For $u$ outside the unit disk, both series will be conditionally convergent. Hence within the strip, we find after 
interchanging the order of the integration and summation that
\begin{eqnarray}
{\rm arctan} \, u = \frac{1}{2} \int_0^{u} dt \Bigl( e^{i \pi k/2} t^k + e^{-i \pi k/2} t^k \Bigr)
 = \sum_{k=0}^{\infty} \cos(\pi k/2) \; \frac{u^{k+1}}{k+1} \;.
\label{seven}\end{eqnarray}
The series in Eq.\ (\ref{seven}) vanishes for odd values of $k$. Consequently, we can replace $k$ by $2k$ and sum over $k$ from zero 
to infinity. When this is done, we arrive at the first result in the lemma.  
The geometric series corresponding to the first term on the rhs of Eq.\ (\ref{six}) will be divergent or yields an infnity when
$\Re\, (it) >1$, while the second geometric series is divergent for $\Re\, (it)<-1$. Hence one of the series will always be divergent 
when $u$ not situated within the strip given by $-1 < \Re\, (iu) <1$. When the infinity is removed from the geometric series or it is
regularized, one is left with a finite value, which is identical to its value when it is convergent. Therefore, outside the strip the 
rhs of Eq.\ (\ref{six}) represents the regularized value of both series. Because the rhs is now not the actual value of ${\rm arctan}\, u$, 
the equals sign no longer applies. Instead, it is replaced by an equivalence symbol indicating that ${\rm arctan}\, u$ is only equivalent 
to the material on the rhs. Hence for either $\Re \,(iu) \geq 1$ or $\Re \, (iu) \leq -1$, we arrive at 
\begin{eqnarray}
{\rm arctan} \, u \equiv \frac{1}{2} \int_0^{u} dt \sum_{k=0}^{\infty} \Bigl( (it)^k + (-it)^k \Bigr) \;.
\label{eight}\end{eqnarray}
By interchanging the order of the integration and summation, we can evaluate the above integral, which yields the same power series in 
Eq.\ (\ref{seven}).  For $\Re\, (iu) = \pm 1$ or the border between the strip and the divergent regions, the series is undefined, but
because the regularized value of the divergent region is the same as that for the convergent strip, we can extend the regularized value to
include $\Re\, (iu) =\pm 1$. Since the resulting series vanishes for odd values of $k$, we obtain the second result in the lemma. Note, 
however, that the regularization process does not eliminate the logarithmic singularities at $(iu) =\pm 1$. This completes the proof of
the lemma.

It should also be mentioned that since the equivalence symbol is less stringent than the equals sign, we can replace the latter by
an equivalence symbol, which wold be valid for all values of $u$. That is, the result in the lemma can be expressed as
\begin{eqnarray}
\sum_{k=0}^{\infty} \frac{(-1)^k u^{2k+1}}{(2k+1)} \equiv {\rm arctan}\,u,  \quad \forall u.
\label{sevena}\end{eqnarray} 
Therefore, if we encounter the series in a problem, then we can replace it by the term on the rhs. As a result, we generate equivalence
statements, not equations, although it need not imply that the lhs was originally divergent. The reader will become more familiar with 
this point as we proceed. 

We are now in a position to derive the complete form of Stirling's approximation. Our starting point is Binet's second expression for 
$\ln \Gamma(z)$, which is derived in terms of an infinite integral in Sec.\ 12.32 of Ref.\ \cite{whi73}. There it is given as
\begin{eqnarray}
\ln \Gamma(z) = \Bigl( z- \frac{1}{2} \Bigr) \ln z -z + \frac{1}{2} \,\ln (2 \pi) + 2 \int_0^{\infty} dt \; 
\frac{{\rm arctan}(t/z)} {e^{2 \pi t}-1} \;\;. 
\label{nine}\end{eqnarray}
Next we make a change of variable, $y = 2\pi t$. Since $z$ is complex, we introduce Eq.\ (\ref{sevena}) thereby replacing the ${\rm arctan}$ 
function, thereby generating an equivalence statement or an equivalence for short. By substituting $k$ by $k+1$, we obtain
\begin{eqnarray}
\ln \Gamma(z) - \Bigl( z- \frac{1}{2} \Bigr) \ln z +z - \frac{1}{2} \,\ln (2 \pi) \equiv \frac{1}{\pi}  \sum_{k=1}^{\infty}
\frac{(-1)^{k+1}}{(2k-1)}  \, \Bigl( \frac{1}{2\pi z}\Bigr)^{2k-1} \int_0^{\infty} dy \; \frac{y^{2k-1}} {e^y -1} \;. 
\label{ten}\end{eqnarray}
The lhs of the above equivalence is finite (convergent), while the rhs can be either divergent or convergent. According to No.\ 3.411(1) 
in Ref.\ \cite{gra94}, the integral in the above equivalence is equal to $\Gamma(2k) \zeta(2k)$. Therefore, Equivalence\ (\ref{ten}) 
reduces to
\begin{eqnarray}
\ln \Gamma(z) - \Bigl( z- \frac{1}{2} \Bigr) \ln z +z - \frac{1}{2} \,\ln (2 \pi) \equiv 2 z \sum_{k=1}^{\infty}
\frac{(-1)^{k+1}} {(2k-1)} \; \frac{\Gamma(2k) \, \zeta(2k)}{(2\pi z)^{2k}}  \;. 
\label{eleven}\end{eqnarray}
On p.\ 282 of Ref.\ \cite{par01} and p.\ 220 of Ref.\ \cite{par11} Paris and Kaminski introduce the function $\Omega(z)$ to
signify all the terms on the lhs of the above equivalence. Moreover, with the aid of the reflection formula for the
gamma function, viz.\
\begin{eqnarray}
\Gamma(z) \, \Gamma(1-z) = \frac{\pi}{\sin(\pi z)} \;,
\label{elevena}\end{eqnarray}
one can derive the following continuation formula:
\begin{eqnarray}
\Omega(z) +  \Omega\left(ze^{\pm i \pi} \right) = -\ln \left( 1 - e^{\mp 2  i \pi z} \right) \;.
\label{elevenb}\end{eqnarray}
This formula allows one to obtain values of $\ln \Gamma(z)$ for $z$ situated in the left hand complex plane via the corresponding 
values in the right hand complex plane. As we shall see in Thm.\ 2.1, the term on the rhs of the above equation 
emerges when the Stokes phenomenon is considered in the development of the complete asymptotic expansion for $\ln \Gamma(z)$.  

The series in Equivalence\ (\ref{eleven}) can be expressed in terms of the monotonically decreasing positive fractions known 
as the cosecant numbers, $c_k$, which are studied extensively in Ref.\ \cite{kow11}. There they are found to be given by
\begin{align}
c_k & = (-1)^k \,L_{P,k} \Bigl[ (-1)^{N_k} N_k! \prod_{i=1}^k \Bigl( \frac{1}{(2i+1)!}\Bigr)^{\!n_i} \, 
\frac{1}{n_i!} \Bigr]
\nonumber\\
&= \;\; 2 \Bigl(1 -2^{1-2k} \Bigl) \frac{\zeta(2k)}{\pi ^{2k}} \;\;.
\label{twelve}\end{align} 
In the above equation $\zeta(k)$ represents the Riemann zeta function, while $N_k = \sum_{i=1}^k n_i$, where $n_i$ 
denotes the number of occurrences or multiplicity of each element $i$ in the integer partitions scanned by the partition 
operator $L_{P,k}[ \cdot]$ of Ref.\ \cite{kow12}, the latter being defined as
\begin{eqnarray}
L_{P,k}[\cdot] =  \sum_{\scriptstyle n_1, n_2,n_3, \dots,n_{k}=0 \atop{\sum_{i=1}^{k}i n_i =k}}^{k,[k/2], [k/3],\dots,1}  
\;\;.
\label{thirteen}\end{eqnarray}
According to Eqs.\ (338) in Ref.\ \cite{kow11}, we have
\begin{eqnarray}
c_k(1)= -2\zeta(2k)/\pi^{2k} \;,
\label{fourteen}\end{eqnarray}
where the cosecant polynomials, $c_k(x)$, are given by
\begin{eqnarray}
c_k(x) =\sum_{j=0}^{k} \frac{(-1)^j}{(2j)!} \, c_{k-j} \, x^j \;.
\label{fifteen}\end{eqnarray}
Hence we find that
\begin{eqnarray}
c_k(1)=  c_k / \left(2^{1-2k}-1 \right)\;,
\label{fifteena}\end{eqnarray}
while the rhs of Equivalence\ (\ref{eleven}) can be expressed as
\begin{eqnarray}
S(z) = z \sum_{k=1}^{\infty} \frac{(-1)^k}{(2z)^{2k}}\, \Gamma(2k-1)\, c_k(1) \;.
\label{sixteen}\end{eqnarray}
Moreover, the standard form of Stirling's approximation in terms of the Bernoulli numbers, viz.\ Eq.\ (\ref{three}), can be 
obtained with the aid of Eq.\ (366) in Ref.\ \cite{kow11}, which gives
\begin{eqnarray}
c_k(1)  = \frac{(-1)^{k}}{(2k)!} \; 2^{2k}\, B_{2k}\;.
\label{seventeen}\end{eqnarray}
Because the gamma function appears in the summand, the radius of absolute convergence for the infinite series, $S(z)$, is zero. 
As a result of this observation, we shall employ the following definition throughout the remainder of this work.
 
\begin{definition}
An asymptotic (power) series is defined here as an infinite series with zero radius of absolute convergence.
\end{definition}
This may seem a strange or even facile definition compared with the infinite series appearing in standard asymptotics 
based on the Poincar$\acute{\rm e}$ prescription, but it is, in fact, quite general. Infinite power series with a finite
of radius of absolute convergence simply do not display asymptotic behaviour. Asymptotic behaviour is where successive terms 
in a series continue to approach the limit or actual value of the function the series represents. However, they only reach 
the closest value at the optimal point of truncation before continually diverging from this closest value at higher truncation 
points. An infinite series with a finite radius of absolute convergence has its optimal point of truncation at infinity
and as a consequence, does not diverge inside its radius of absolute convergence. Outside the radius of absolute convergence, 
however, the series is divergent, but then it does not possess an optimal point of truncation. Hence an infinite series
with a finite radius of convergence does not display asymptotic behaviour. 

Another aspect of the above definition for an asymptotic series is that it does not explicitly mention divergence. This
is because an asymptotic series cannot always be divergent. If this were the case, then the function represented by 
the series would be infinite everywhere and, consequently of very little use or interest just as its inverse function 
would be zero everywhere. Thus, an asymptotic series is divergent for a specific range of values in the complex plane,
while it is conditionally convergent for the remaining values.  

It is well-known that an asymptotic expansion is only valid over a sector in the complex plane and that on reaching the 
boundary of such a sector, it acquires a discontinuous term and yet another on moving to an adjacent sector. This is 
known as the Stokes phenomenon \cite{din73,kow09}, while the sectors over which a complete asymptotic expansion is uniform 
and their boundaries are called Stokes sectors and lines, respectively. As a consequence, it is simply not sufficient to 
give an asymptotic expansion without specifying the Stokes sector or line where it is uniform. Therefore, we require the
following definition. 

\begin{definition}
An asymptotic form is defined here as being composed of: (1) a complete asymptotic expansion, which not only possesses
all terms in a dominant asymptotic power series such as $S(z)$ above, but also all the terms in each subdominant 
asymptotic series, should they exist, and (2) the range of values or the sector/ray in the complex plane over which the 
argument of the variable in the series is valid.
\end{definition}

We now truncate the asymptotic series $S(z)$ at $N$ terms, which will be referred to as the truncation parameter. Moreover, 
in the resulting infinite series we express $c_k(1)$ in terms of the Riemann zeta function via Eq.\ (\ref{fourteen}) and 
replace the latter by its Dirichlet series form. Consequently, we find that 
\begin{eqnarray}
S(z) = z \sum_{k=1}^{N-1} \frac{(-1)^k}{(2z)^{2k}}\, \Gamma(2k-1)\,c_k(1) - 2 z \sum_{n=1}^{\infty}
\sum_{k=N}^{\infty} \frac{(-1)^k}{(2 \pi n z)^{2k}}\; \Gamma(2k-1) \;.
\label{eighteen}\end{eqnarray}
The series over $k$ in the second term on the rhs of the above equivalence is an example of a Type I generalized terminant. 
Terminants were first introduced by Dingle in Ch.\ 22 of Ref.\ \cite{din73}, who discovered that a great number of 
special functions in mathematical physics possessed asymptotic expansions in which the coefficients could eventually be 
approximated by coefficients with gamma function growth, viz.\ $\Gamma(k \!+\! \alpha)$. In the case of a Type I terminant 
the coefficient is also accompanied by a phase factor of $(-1)^k$, while a Type II terminant does not possess such a phase 
factor in its summand. Subsequently, terminants were generalized in Ref.\ \cite{kow09} to series where the coefficients 
were given by $\Gamma(pk/q \!+\! \alpha)$ with $p$ and $q$ positively real, instead of $\Gamma(k \!+\! \alpha)$.
  
In Ch.\ 10 of Ref.\ \cite{kow09} the notation $S^{I}_{p,q}(N,z^{\beta})$ was used to denote a Type I generalized terminant.
These are given by
\begin{eqnarray}
S^{I}_{p,q}\left(N,z^{\beta} \right) = \sum_{k=N}^{\infty} (-1)^k\, \Gamma(pk+q) z^{\beta\,k} \;\;.
\label{nineteen}\end{eqnarray}
As a result, we can express Eq.\ (\ref{eighteen}) as 
\begin{eqnarray}
S(z) = z \sum_{k=1}^{N-1} \frac{(-1)^k}{(2z)^{2k}}\, \Gamma(2k-1)\,c_k(1)  
- 2 z \sum_{n=1}^{\infty} S^{I}_{2,-1}\left(N,(1/2n \pi z)^2 \right)  \;.
\label{twenty}\end{eqnarray}  
Specifically, in the above equation we see that $\beta \!=\!2$, while the power variable is equal to $1/2n \pi z$. Although it is 
stated in Ref.\ \cite{kow09} that both $p$ and $q$ have to be positive and real, it is actually the value of $N \!+\!q/p$ that 
appears in the regularized value of a generalized terminant as we shall see soon. This means that as long as the real part of 
this quantity is greater than zero, then the form for the regularized value of the series given in Ref.\ \cite{kow09}
will still apply. Alternatively, the problem can be avoided by replacing $k$ by $k \!+\! 1$ in the series, in which case 
the "new" value of $q$ becomes unity. In other words, since $S^I_{2,-1}\bigl(N,z^2\bigr)  \!=\!-z^2 S^I_{2,1}\bigl(N-1,z^2 \bigr)$,
we can apply the result in Ref.\ \cite{kow09} to $S^I_{2,1}\bigl( N-1,z^2 \bigr)$, thereby avoiding the fact that $q$ is negative
in $S^{I}_{2,-1}\bigl(N,z^2\bigr)$.

According to Rule A in Ch.\ 1 of Ref.\ \cite{din73}, Stokes lines occur whenever the terms in an asymptotic series are 
determined by those phases or arguments for which successive late terms are homogeneous in phase and are all
of the same sign. In the case of the generalized terminant in Eq.\ (\ref{twenty}), Dingle's rule means that Stokes lines occur 
whenever ${\rm arg}\left( -1/z^2 \right)\!=\! 2l \pi$, for $l$, an integer. Under this condition all the terms in either
$S^{I}_{2,-1}\bigl(N,1/z^2)$ or $S^I_{2,1}\bigl(N,1/z^2)$ are all positively real. Because of the arbitrariness of $l$, we 
can replace -1 in this equation by $\exp(-i\pi)$. Then we find that the Stokes lines for $S(z)$ occur whenever
${\rm arg}\, z = -(l+ 1/2) \pi$, i.e. at half integer multiples of $\pi$. 
    
In Ref.\ \cite{kow09} the concept of a primary Stokes sector was introduced, which was necessary for indicating the
Stokes sector over which an asymptotic expansion does not possess a Stokes discontinuity such as Equivalence\ (\ref{twenty}).
It is also necessary for defining asymptotic forms since two functions can have the same complete asymptotic 
expansion, but then they would be valid over different sectors or rays. On encountering the Stokes lines at the boundaries 
of the primary Stokes sector, jump discontinuities need to be evaluated according to the Stokes phenomenon. 
Then as a secondary Stokes sector is encountered either in a clockwise or anti-clockwise direction from the primary 
Stokes sector, still more Stokes discontinuities need to be evaluated. The choice of a primary Stokes sector is arbitrary, 
but we shall choose it to be the Stokes sector lying in the principal branch of the complex plane, since 
most asymptotic expansions are derived under the condition that the variable lies initially in the principal 
branch of the complex plane. 

Before we can investigate the regularization of the asymptotic series, $S(z)$, we require the following lemma:
\begin{lemma}
Regularization of the power/Taylor series for the logarithmic function, viz.\ $\log \bigl( 1+z \bigr)$, yields 
\begin{eqnarray}
\sum_{k=1}^{\infty} \frac{ (-1)^{k+1}} {k} \; z^k \begin{cases} \equiv \ln(1+z) \quad, & \quad \Re\, z \leq -1\quad, \cr
=\ln(1+z) \quad, & \quad \Re\, z > -1 \quad . \end{cases}
\label{twentya}\end{eqnarray}
\end{lemma}
{\bfseries Proof}. There is no need to present the proof for the above result as this has been already carried out in Refs.\ 
\cite{kow11b} and \cite{kow11c} via the integration of the geometric series. The regions in the complex plane where the logarithmic 
series is divergent and convergent are also determined in these references. As was found for the arctan function, the 
logarithmic series is absolutely convergent within the unit disk and either conditionally convergent or divergent outside it. 
The logarithmic series is usually cited as the classic example of a conditionally convergent for positive real values of $z$ 
in the above series as described on p.\ 18 of Ref.\ \cite{whi73}. Moreover, like the geometric and binomial series, the 
regularized value of the above series is the same value as when the series is convergent. Hence we can replace the equals 
sign in the lemma by the less stringent equivalence symbol, which means in turn that 
\begin{eqnarray}
\sum_{k=1}^{\infty} \frac{(-z)^k}{k} \equiv -\ln(1+z) \;, \quad \forall\, z \;.
\label{twentyb}\end{eqnarray}
This completes the proof of the lemma.

We are now in a position to regularize $S(z)$, which means that we can derive the asymptotic forms for $\ln \Gamma(z)$
based on Stirling's approximation. 
\begin{theorem}
Via the regularization of the asymptotic power series given by Eq.\ (\ref{eighteen}), the logarithm of the gamma 
function can be expressed in the following asymptotic forms:
\begin{align} 
\ln \Gamma(z) & =  \Bigl( z- \frac{1}{2} \Bigr) \ln z - z + \frac{1}{2} \,\ln (2 \pi) + 
 z \sum_{k=1}^{N-1} \frac{(-1)^k}{(2z)^{2k}}\, \Gamma(2k-1)\,c_k(1)  
\nonumber\\
& +  \;\; R^{SS}_N(z) + SD^{SS}_M(z)  \;,
\nonumber\\
\label{twentyone}\end{align}
where the remainder $R^{SS}_N(z)$ is given by
\begin{align}
R_N^{SS}(z) & =  \frac{2\,(-1)^{N+1} \,z}{(2\pi z)^{2N}} \int_0^{\infty} dy \; y^{2N-2}\, e^{-y}  
\sum_{n=1}^{\infty} \frac{1}{n^{2N-2} \left( (y/2 \pi z)^2 + n^2 \right)} \;,
\label{twentytwo}\end{align}
and the Stokes discontinuity term $SD_M(z)$ is given by
\begin{align}
SD^{SS}_M(z) & =  -  \lfloor M/2 \rfloor \, \ln \left(-\, e^{\pm 2 i \pi z} \right) - \frac{\left( 1 - (-1)^M \right)}{2} \,
\ln \left( 1- e^{\pm 2i \pi z }\right) \;. 
\label{twentytwoa}\end{align}
The remainder given by Eq.\ (\ref{twentytwo}) is found to be valid for either $(M-1/2) \pi \!<\! \theta \!=\! {\rm arg}\, z \!<\! (M+1/2) \pi$ 
or $-(M+1/2) \pi \!<\! \theta \!<\! -(M-1/2) \pi$, where $M$ is a non-negative integer. However, the Stokes discontinuity term given by Eq.\ 
(\ref{twentytwoa}) has two versions, complex conjugates of one another as evidenced by $\pm$ sign. In this instance the upper-signed version 
of Eq.\ (\ref{twentytwoa}) applies to $(M-1/2) \pi \!<\! \theta \!<\! (M+1/2) \pi$, while the lower-signed version is valid over 
$-(M+1/2) \pi \!<\! \theta  \!<\! -(M-1/2) \pi$.
For the situation along the Stokes lines or rays, where $\theta \!=\! \pm (M+1/2) \pi$, we replace $R^{SS}_N(z)$ and $SD^{SS}_M(z)$
by $R^{SL}_N(z)$ and $SD^{SL}_M(z)$, respectively. Then the remainder is found to be 
\begin{align}
R_N^{SL}(z) & = \frac{2z}{(2 \pi |z|)^{2N-2}} \, P \int_0^{\infty} dy \; y^{2N-2}\, e^{-y} \sum_{n=1}^{\infty} 
\frac{1}{n^{2N-2}(y^2 -4n^2 \pi^2 |z|^2)} \;\;,
\label{twentythree}\end{align}
while the Stokes discontinuity term is given by
\begin{align}
SD^{SL}_{M}(z) & =  (-1)^M \Bigl( \lfloor M/2 \rfloor +\frac{1- (-1)^M}{2} \, \Bigr) 2\pi |z|  - \frac{1}{2} \; \ln \Bigl( 1
- e^{- 2 \pi |z|} \Bigr) \;. 
\label{twentyfour}\end{align}
\end{theorem}
In Eq.\ (\ref{twentythree}) $P$ denotes that the Cauchy principal value must be evaluated.

\begin{remark}
The above results may not appear to be asymptotic, but we shall see this more clearly when we carry out numerical studies, where the 
truncation parameter exceeds the optimal point of truncation.
\end{remark} 

{\bfseries Proof}.  The regularized value of a generalized Type I terminant is derived in Ch.\ 10 of Ref.\ \cite{kow09} and is also discussed
in detail in Ref.\ \cite{kow14b}.
For $-(2M+1) \pi/\beta \!<\! \theta \!= \! {\rm arg}\, z \!<\! -(2M-1) \pi/\beta$, it is given as Equivalence\ (10.21), which is
\begin{align}
S^{I}_{p,q} \left( N, z^{\beta} \right) &\equiv (-1)^N z_M^{\beta(N-1)} p^{-1} \int_0^{\infty} dt\; \frac{t^{N+q/p-1} \,
e^{-t^{1/p}}}{t+z_M^{-\beta}} + 2 i \pi  z_M^{-\beta q/p} \;\;
\nonumber\\
& \times \;\; p^{-1} \sum_{j=1}^M e^{(2j-1)q i \pi/p} \, \exp\Bigl( -z_M^{-\beta/p} e^{(2j-1) i\pi/p} \Bigr) \;\;,
\label{twentyfive}\end{align}
where $z_M \!=\! z \exp(2Mi \pi/\beta)$. The last term in this result represents the Stokes discontinuity term with $M$ being  
the number of Stokes sectors that have been traversed. It arises by summing the residues due to the singularity at $t=-z_M^{-\beta}$ 
in the Cauchy integral in the above result for each traversal of a Stokes line/ray. When $M \!=\!0$, it vanishes and thus, it gives 
the value corresponding to the primary Stokes sector, which, as stated above, has been selected to lie in the principal branch 
of the complex plane. Substituting $t$, $z$, $\beta$, $p$ and $q$ in Equivalence\ (\ref{twentyfive}) by $y^2$, $1/2n \pi z$, 
2, 2 and -1, respectively, we arrive at
\begin{align}
S^{I}_{2,-1} \left( N, (1/2n \pi z)^2 \right) &\equiv \frac{(-1)^N}{(2n\pi z)^{2N-2}} \int_0^{\infty} dy\; \frac{y^{2N-2} \,
e^{-y}}{y^2+4 n^2 \pi^2 z^2} - \frac{1}{2 n z} 
\nonumber\\
& \times \;\; \sum_{j=1}^M (-1)^{M-j}\, \exp\Bigl( -2(-1)^{M-j} n i \pi z \Bigr) \;.
\label{twentyfivea}\end{align}
By introducing the above result into Eq.\ (\ref{twenty}) and carrying out some manipulation, we obtain
\begin{align}
S(z) &  \equiv  z \sum_{k=1}^{N-1} \frac{(-1)^k}{(2z)^{2k}}\, \Gamma(2k-1)\,c_k(1) 
 -  2  \, \Bigl( -\frac{1}{4 \pi^2 z^2} \Bigr) ^{N} \,z \int_0^{\infty} dy \; e^{-y}\, y^{2N-2}
\nonumber\\   
& \times \;\; \sum_{n=1}^{\infty} \frac{1}{n^{2N-2}\, ((y/2\pi z)^2 +n^2)} 
+  \sum_{n=1}^{\infty} \frac{1}{n} \sum_{j=1}^M (-1)^{M-j} \exp \left(2  (-1)^{M-j}\, n i \pi z\right) \;.
\label{twentysix}\end{align} 

Unfortunately, Equivalence\ (\ref{twentysix}) is only a partially regularized result. Because of the summation over $n$, we
are no longer considering one generalized terminant as in Equivalences\ (\ref{twentyfive}) and (\ref{twentyfivea}). Now we
are considering an infinite number of terminants for each Stokes line. Each one of these generalized terminants possesses a Stokes 
discontinuity term as represented by the second term on the rhs of Equivalence\ (\ref{twentyfivea}). Alternatively, we can regard 
each Stokes line as possessing an infinite number of singularities. As a consequence, we obtain another infinite series when summing
all their Stokes discontinuity terms. This series can become divergent. Hence we have to regularize this second series, which would 
not have been necessary if there were only a finite number of generalized terminants. As we shall see, the regularization of the
series over $n$ in the last term of Equivalence\ (\ref{twentysix}) is responsible for providing the logarithmic behaviour including 
the multivaluedness arising from taking the logarithm of the gamma function. That is, the fact that the logarithm is a multivalued 
function means that $\ln \Gamma(z)$ is also a multivalued function. This behaviour does not manifest itself in the other terms 
on the rhs of Equivalence\ (\ref{twentysix}), but in the last term after it has been regularized, which will be seen more clearly 
in the numerical study of Sec.\ 5. 

We shall refer to the last term on the rhs of Eq.\ (\ref{twentysix}) as the Stokes discontinuity term and denote it by $SD_M(z)$. 
It can, however, be simplified drastically by considering the cases of odd and even values of $M$ separately. For $M \!=\! 2J$, the 
Stokes discontinuity term reduces to
\begin{eqnarray}
SD_{2J}(z) =J \sum_{n=1}^{\infty} \frac{1}{n} \; \bigl( \exp(2n i \pi z)- \exp(-2n i \pi z) \bigr) \;\;,
\label{twentyseven}\end{eqnarray}
while for $M \!=\! 2J+1$, we obtain
\begin{eqnarray}
SD_{2J+1}(z) =J \sum_{n=1}^{\infty} \frac{1}{n} \; \bigl( \exp(2n i\pi z)- \exp(-2n i\pi z  ) \bigr) 
+ \sum_{n=1}^{\infty} \frac{1}{n} \; \exp\left( 2n i \pi z \right) \;.
\label{twentyeight}\end{eqnarray}
Since $J= \lfloor M/2 \rfloor$ regardless of whether $M$ is even or odd, we can combine the two preceding equations into one
result, which yields
\begin{eqnarray}
SD_M(z) =\lfloor M/2 \rfloor \sum_{n=1}^{\infty} \frac{1}{n} \; \Bigl( e^{2n i \pi z }- e^{-2n i \pi z} \Bigr) 
+ \frac{\left( 1-(-1)^M \right)}{2} \sum_{n=1}^{\infty} \frac{e^{2n i \pi z}}{n} \;.
\label{twentynine}\end{eqnarray}

It has already been stated that Equivalence\ (\ref{twentysix}) is only partially regularized. That is, in order to obtain the
regularized value of $S(z)$, we need to regularize the Stokes discontinuity term. In fact, all the series in Eq.\ (\ref{twentynine}) 
are variants of the logarithmic series. That is, they can be written as $\sum_{n=1}^{\infty} n^{-1} z^n$ with 
$z \!=\! \exp(\pm 2\pi z i)$. As a consequence, we can use Lemma\ 2.2 to obtain the regularized value of the Stokes discontinuity
term and hence, $S(z)$. 

If we introduce the regularized value of the logarithmic series as given by Equivalence\ (\ref{twentyb}) into $SD_M(z)$, then after a 
little manipulation we arrive at
\begin{eqnarray}
SD_M(z) \equiv -\lfloor M/2 \rfloor \ln \left( - \,e^{2i \pi z}\right) + \frac{\left(1- (-1)^M \right)}{2} \ln 
\left( 1- e^{2 i \pi z}\right) \;. 
\label{thirtyone}\end{eqnarray}  
Furthermore, by replacing the Stokes discontinuity term or the second term on the rhs of Equivalence\ (\ref{twentysix}) by the rhs of 
Equivalence\ (\ref{thirtyone}), we obtain the regularized value of $S(z)$ for $(M-1/2) \pi < \theta <(M+1/2) \pi$. This is
\begin{align}
& S(z)   \equiv  z \sum_{k=1}^{N-1} \frac{(-1)^k}{(2z)^{2k}}\; \Gamma(2k-1)\,c_k(1) 
 - 2  \, \Bigl( -\frac{1}{4 \pi^2 z^2} \Bigr) ^{N} \,z \int_0^{\infty} dy \; e^{-y}\, y^{2N-2}
\nonumber\\   
& \times \;\; \sum_{n=1}^{\infty} \frac{1}{n^{2N-2}\, ((y/2\pi z)^2 +n^2)} +\lfloor M/2 \rfloor \ln \left( -\, e^{2i \pi z} \right) - 
\frac{\left(1- (-1)^M \right)}{2} \ln \left( 1- e^{2 i \pi z }\right) \;.
\label{thirtyonea}\end{align} 
As explained in Refs.\ \cite{kow002}, \cite{kow09}, \cite{kow11} and \cite{kow09a}, the regularized value is a unique quantity,
which means that if an infinite power series possesses two different forms for it, then they are equal to another. This means that
we can equate the rhs of Equivalence\ (\ref{thirtyone}) with the lhs of Equivalence\ (\ref{eleven}). Therefore, we find that 
\begin{align}
 \ln \Gamma(z) & =  \Bigl( z- \frac{1}{2} \Bigr) \ln z -z +\frac{1}{2} \,\ln (2 \pi) + z \sum_{k=1}^{N-1} \frac{(-1)^k}{(2z)^{2k}}\, 
\Gamma(2k-1)\, c_k(1) - 2  \, \Bigl( -\frac{1}{4 \pi^2 z^2} \Bigr) ^{N} 
\nonumber\\
& \times \;\; z \int_0^{\infty} dy \; e^{-y}\, y^{2N-2} \sum_{n=1}^{\infty} \frac{1}{n^{2N-2}\, ((y/2\pi z)^2 +n^2)} 
+\lfloor M/2 \rfloor \ln \left(- \, e^{2 i\pi z}\right) 
\nonumber\\
& - \;\; \frac{\left(1- (-1)^M \right)}{2} \ln \left( 1- e^{2i \pi z}\right) \;,
\label{thirtytwo}\end{align} 
for $(M-1/2) \pi < \theta < (M+1/2) \pi$. 
Note the appearance of the same logarithmic term in the final term when $M$ is odd as in the continuation formula involving $\Omega(z)$ and 
$\Omega(ze^{-i\pi})$ in Eq.\ (\ref{elevenb}).

So far, we have only considered positive values or anti-clockwise rotations of $\theta$. To derive the value of $\ln \Gamma(z)$ 
for the clockwise rotations of $\theta$, we require the regularized value of $S^I_{p,q}(N,z^{\beta})$ 
when $(2M-1)/\beta \!<\! \theta \!<\! (2M+1)/\beta$. This result, which appears as Equivalence\ (10.37) in Ref.\ \cite{kow09} and is
also discussed in Ref.\ \cite{kow14b}, is given by
\begin{align}
S^{I}_{p,q} \left( N, z^{\beta} \right) &\equiv (-1)^N z_{-M}^{\beta(N-1)} p^{-1} \int_0^{\infty} dt\; \frac{t^{N+q/p-1} \,
e^{-t^{1/p}}}{t+z_{-M}^{-\beta}} - 2 \pi i z_{-M}^{-\beta q/p} \;\;
\nonumber\\
& \times \;\; p^{-1} \sum_{j=1}^M e^{-i(2j-1)q \pi/p} \, \exp\Bigl( -z_{-M}^{-\beta/p} e^{-i(2j-1) \pi/p} \Bigr) \;\;,
\label{thirtythree}\end{align}
and $z_{-M} \!=\! z \exp(-2Mi \pi/\beta)$. Note that when $M \!=\!0$, the above result agrees with the $M=0$ case of Equivalence\ 
(\ref{twentyfive}), which is expected since it reduces to the regularized value over the primary Stokes sector. For non-vanishing $M$, 
however, the above result represents the complex conjugate of Equivalence\ (\ref{twentyfive}). By introducing the above result into 
Eq.\ (\ref{twenty}), we obtain
\begin{align}
S(z) & \equiv  z \sum_{k=1}^{N-1} \frac{(-1)^k}{(2z)^{2k}}\, \Gamma(2k-1)\,c_k(1) 
- 2  \, \Bigl( -\frac{1}{4 \pi^2 z^2} \Bigr) ^{N} \,z \int_0^{\infty} dy \; e^{-y}\, y^{2N-2} 
\nonumber\\  
& \times \;\; \sum_{n=1}^{\infty} \frac{1}{n^{2N-2}\, ((y/2\pi z)^2 +n^2)} 
+ \sum_{n=1}^{\infty} \frac{1}{n} \sum_{j=1}^M (-1)^{M-j} \exp \left(-2  (-1)^{M-j}\, n \pi z i \right) \;.
\label{thirtyfour}\end{align}

Once again, this result is partially regularized since the Stokes discontinuity term or the second term on the rhs can become
divergent. Moreover, the only difference between the above result and Equivalence\ (\ref{twentysix}) is that the Stokes discontinuity 
term is now the complex conjugate of Equivalence\ (\ref{twentynine}). Consequently, the same analysis leading to Equivalence\ 
(\ref{thirtyone}) can be employed. Hence for $-(M+1/2) \pi < \theta <-(M+1/2) \pi$, we arrive at  
\begin{align}
S(z)  &  \equiv  z \sum_{k=1}^{N-1} \frac{(-1)^k}{(2z)^{2k}}\; \Gamma(2k-1)\,c_k(1) - 2  \, \Bigl( -\frac{1}{4 \pi^2 z^2} \Bigr)^{N} 
\,z \int_0^{\infty} dy \; e^{-y}\, y^{2N-2} \sum_{n=1}^{\infty} \frac{1}{n^{2N-2}}
\nonumber\\   
& \times \;\; \frac{1}{((y/2\pi z)^2 +n^2)} -\lfloor M/2 \rfloor \ln \left(- \, e^{-2 i\pi z} \right) -\frac{\left(1- (-1)^M \right)}{2} 
\ln \left( 1- e^{-2i \pi z }\right) \;.
\label{thirtyfive}\end{align} 
Equating the rhs of Equivalence\ (\ref{thirtyfive}) with the lhs of Equivalence\ (\ref{eleven}) yields
\begin{align}
 \ln \Gamma(z) & =  \Bigl( z- \frac{1}{2} \Bigr) \ln z -z +\frac{1}{2} \,\ln (2 \pi) + z \sum_{k=1}^{N-1} \frac{(-1)^k}{(2z)^{2k}}\, 
\Gamma(2k-1)\, c_k(1) - 2  \, \Bigl( -\frac{1}{4 \pi^2 z^2} \Bigr) ^{N} 
\nonumber\\
& \times \;\; z \int_0^{\infty} dy \; e^{-y}\, y^{2N-2} \sum_{n=1}^{\infty} \frac{1}{n^{2N-2}\, ((y/2\pi z)^2 +n^2)} 
- \lfloor M/2 \rfloor \ln \left(- \, e^{-2i \pi z} \right)
\nonumber\\
& - \;\; \frac{\left(1- (-1)^M \right)}{2} \ln \left( 1- e^{-2i \pi z}\right) \;.
\label{thirtysix}\end{align} 
As expected, the above result is identical to Equivalence\ (\ref{thirtytwo}) except that the Stokes discontinuity term or those terms
with $M$ in them are complex conjugates. In this case the final logarithmic term for odd values $M$ is the same term that appears in 
the continuation formula involving $\Omega(z)$ and $\Omega(ze^{i\pi})$ in Eq.\ (\ref{elevenb}). Moreover, if we combine Equivalences\ 
(\ref{thirtytwo}) and (\ref{thirtysix}) into one expression, then we obtain the equivalence given in the theorem, viz.\ Equivalence\ 
(\ref{twentyone}) with the remainder and Stokes discontinuity term given by Eqs.\ (\ref{twentytwo}) and (\ref{twentytwoa}) respectively.

The terms on the rhs's of Eqs.\ (\ref{thirtytwo}) and (\ref{thirtysix}) up to the first sum over $k$ are the standard terms or
leading order terms in Stirling's approximation as given by ``Eq."\ (\ref{one}), but with $n$ replaced by $z$. The truncated sum 
over $k$ descends in powers of $1/z^2$ with the first power being $1/z$, while the coefficients are equal to $(-1)^{k+1} 2^{1-2k} 
\zeta(2k)/\pi^{2k}$. Therefore, we find that for all values up to $k\!=\!4$, the coefficients are identical to those appearing in 
No.\ 6.1.41 of Ref.\ \cite{abr65} or ``Eq."\ (\ref{two}) here. The next term is the most fascinating of all the terms since it represents 
an exact value for the remainder, $R^{SS}_N(z)$, a quantity that is at best bounded in standard Poincar$\acute{\rm e}$ asymptotics. 
That is, it represents part of the regularized value of the asymptotic series appearing in Eq.\ (\ref{three}) and is finite, although
as we shall observe in the next section, it diverges for very large values of the truncation parameter $N$. From here on, we shall use 
$R_N(z)$ to denote this finite value, which when added to the other terms on the rhs of of Eq.\ (\ref{thirtytwo}) yields the exact value 
for $\ln \Gamma(z)$. The superscript SS denotes that it represents the remainder for the Stokes sectors as opposed to the superscript SL,
which will used shortly in reference to the remainder for Stokes lines/rays. If we require the series in Eq.\ (\ref{three}), then an 
equivalence symbol will be used when it is equated to $R_N(z)$. Hence we shall no longer regard the remainder as an infinite quantity.  
 
The remaining terms in Eqs.\ (\ref{thirtytwo}) and (\ref{thirtysix}) are only non-zero when $z$ is no longer situated in the primary 
Stokes sector given by $|\theta| < \pi/2$. For $\pi/2 <\theta <3 \pi/2$ or $M=1$, the term with $\lfloor M/2 \rfloor$ in Eq.\ (\ref{thirtytwo})
vanishes, while the final term does not. A similar situation occurs for $-3\pi/2 < \theta <-\pi/2$ and Eq.\ (\ref{thirtysix}).
In the vicinity of the Stokes lines, i.e.\ near the positive or negative imaginary axes, either term is exponentially subdominant for large 
values of $|z|$, but becomes more significant as $\theta$ moves away from the Stokes lines. This is consistent with Rule BA on p.\ 7
of Ref.\ \cite{din73}, which states that relative to its associate (the Stokes discontinuity term), an asymptotic series ($S(z)$) is
dominant where its late terms are of uniform sign, the condition we used to determine the Stokes lines in the first place. However, 
provided $\theta$ is less than $\pi$ and greater than $-\pi$ ($\theta \!=\! \pm \pi$ are known as an anti-Stokes lines), the 
second term remains small when compared with the truncated term as long as $|z|$ is large. This will become more apparent when 
we carry out a numerical study later. This means that the Stokes discontinuity term can be neglected for $|\theta |< \pi$, which is 
why the tilde has been introduced as in ``Eq."\ (\ref{two}).
 
The results presented above are not valid if we wish to determine the regularized value of $S_{2,-1}(N,z^{\beta})$ at a Stokes line 
or ray. In this case we require the regularized value of a generalized Type I terminant or $S^I_{p,q}(N,z^{\beta})$ at the Stokes lines. 
There are two forms appearing in Ref.\ \cite{kow09}, namely Equivalences\ (10.22) and (10.38). They are also discussed in Ref.\ \cite{kow14b}.
However, since they are complex conjugates of one another, they can be combined into one statement, which can be expressed as
\begin{align}
& S^{I}_{p,q} \Bigl( N, z^{\beta}\Bigr) \equiv -|z|^{\beta(N-1)}\, p^{-1} P \int_0^{\infty}dt\; \frac{t^{N+q/p-1}\, e^{-t^{1/p}}}
{t- |z|^{-\beta}} \pm \frac{2 \pi i}{p} \, |z|^{-\beta q/p}
\nonumber\\
& \times \;\; \sum_{j=1}^M e^{\pm 2jq \pi i/p} \exp \Bigl( -|z|^{-\beta/p} e^{\pm 2j \pi i/p}\Bigr) \pm \frac{\pi i}{p} \; 
|z|^{-\beta q/p} \exp\Bigl(-|z|^{-\beta/p} \Bigr) \;,
\label{thirtyseven}\end{align} 
where $P$ denotes that the Cauchy principal value must be evaluated so that we avoid the pole at $t \!=\! |z|^{-\beta}$. The 
upper-signed version in the above equivalence corresponds to $\theta \!=\! -(2M+1) \pi/\beta$, while the lower-signed version corresponds 
to $\theta \!=\! (2M+1)\pi/\beta$. The regularized value of $S_{2,-1}(N,1/(2n\pi z)^2)$ at the Stokes lines of $\theta$ or 
${\rm arg}\, z \!=\! \pm (M+1/2) \pi$ follows by introducing the values for $\beta$, $p$ and $q$ into the above result. Then we find that 
\begin{align}
& S^{I}_{2,-1} \Bigl( N,1/(2n \pi z)^2 \Bigr) \equiv - \Bigl( \frac{1}{2n\pi |z|}\Bigr)^{2N-2} \; P \int_0^{\infty}dy\; 
\frac{y^{2N-2}\, e^{-y}}{y^2- 4 n^2 \pi^2 |z|^2} \; \pm \; \frac{i}{2n |z|}
\nonumber\\
& \times \; \sum_{j=1}^M (-1)^j \exp \Bigl( -2 (-1)^j n \pi |z| \Bigr) \pm \frac{i}{4n |z|}
\exp\Bigl( -2n \pi |z|\Bigr) .
\label{thirtysevenb}\end{align} 
In obtaining this result the substitution, $y =\sqrt{t}$, has again been made. Introducing Equivalence\ (\ref{thirtyeight}) into Eq.\ 
(\ref{twenty}) yields
\begin{align}
S(z)  &\equiv  
\mp i e^{i \theta} \sum_{n=1}^{\infty} \frac{1}{n} \sum_{j=1}^M (-1)^{j} \exp \left(2  (-1)^{j}\, n \pi |z| \right) 
\mp i e^{i \theta} \sum_{n=1}^{\infty} \frac{1}{2n} \exp \left(-2 n\pi |z| \right)
\nonumber\\
& + \; z \sum_{k=1}^{N-1} \frac{(-1)^k}{(2z)^{2k}} \, \Gamma(2k-1)\,c_k(1) + 2  \Bigl( \frac{1}{4 \pi^2 |z|^2} \Bigr)^{N-1} 
z \, P \int_0^{\infty} dy \, e^{-y}\, y^{2N-2}  
\nonumber\\
& \times \; \sum_{n=1}^{\infty} \frac{1}{n^{2N-2}\, (y^2 -4n^2 \pi^2 |z|^2)} \;.
\label{thirtyeight}\end{align} 

Because there is an infinite number of singularities on each Stokes line, we have once more an infinite number of residue contributions to
the regularized value. The main difference between being situated in a Stokes sector and on a Stokes line is that the final residue contribution 
in the latter case is a semi-residue contribution, which cannot be included in the sum over $j$. Instead, this contribution appears separately
and is represented by the second term on the rhs of Equivalence\ (\ref{thirtyeight}). From Lemma\ 2.2, we see that for odd values of $j$ 
in the first term on the rhs the series over $n$ is logarithmically divergent, while all the other series including the sum over $n$ 
in the second term on the rhs are convergent. Therefore, the above result has only been partially regularized as in the Stokes sector case. 
A complete regularization requires that the divergent series on the rhs must be regularized too. Furthermore, all the series are composed of 
real terms, which means that the regularized value must also be  real. In other words, a series composed of real terms cannot 
suddenly acquire imaginary terms as discussed on p.\ 10 of Ref.\ \cite{din73}. In Ref.\ \cite{kow09} this is called the Zwaan-Dingle principle. 
Therefore, whilst we can employ Lemma\ 2.2 to determine the regularized values of the series in Equivalence\ (\ref{thirtyeight}), we must 
introduce the $\Re$ function to ensure that only the real part is evaluated. In so doing, discontinuities may appear in $\ln \Gamma(z)$. These 
discontinuities, however, will only arise in the first sum on the rhs of Equivalence\ (\ref{thirtyeight}) since the second sum yields a regularized 
value that is real. That is, the discontinuities will only apply when $M \geq 1$ for $\theta =\pm (M+1/2) \pi$. Although there is a Stokes 
discontinuity term for $\theta = \pm \pi/2$, $\ln \Gamma(z)$ will not be discontinuous there, which demonstrates that Stokes lines can be 
fictitious. After applying Lemma\ 2.2 and carrying out some algebra, we find that the first two terms on the rhs of Equivalence\ 
(\ref{thirtyeight}) reduce to
\begin{align}
& \pm i e^{i \theta} \sum_{n=1}^{\infty} \frac{1}{n} \sum_{j=1}^M (-1)^{j} \exp \left(2  (-1)^{j}\, n \pi |z| \right) 
\pm i e^{i \theta} \sum_{n=1}^{\infty} \frac{1}{2n} \exp \left(-2 n\pi |z| \right)
\nonumber\\
& \equiv \;\; (-1)^M \Bigl(  \frac{M}{2}  + \frac{1-(-1)^M}{2} \Bigr) 2 \pi |z|-
\frac{1}{2} \; \ln \! \left( 1- e^{-2 \pi |z|} \right) \; . 
\label{thirtynine}\end{align}
In obtaining this result the following identity has been used
\begin{eqnarray}
\Re \ln \Bigl(\frac{1-e^{-2|z|}}{1-e^{2|z|}} \Bigr)=\Re \ln \Bigl( -e^{-2 |z|}\Bigr) = -2 |z| \;\;.
\label{thirtyninea}\end{eqnarray}

Although Eq.\ (\ref{thirtynine}) is valid for any integer value of $M$, the second term in the big brackets on the rhs contributes 
only when $M$ is odd. Moreover, the contribution due to the infinite singularities at a specific Stokes line is represented by the second 
term on the rhs, while the contributions from the preceding Stokes lines is represented by the first term. We shall refer to these terms
as the Stokes discontinuity term. Since it only applies to Stokes lines, we denote it by $SD_M^{-}(z)$. Introducing the above 
equivalence into the partially regularized result given by Equivalence\ (\ref{thirtyeight}) yields
\begin{align}
S(z)  &  \equiv  \; z \sum_{k=1}^{N-1} \frac{(-1)^k}{(2z)^{2k}} \, \Gamma(2k-1)\,c_k(1) + 2  \Bigl( \frac{1}{4 \pi^2 |z|^2} \Bigr)^{N-1} 
z \, P \int_0^{\infty} dy \, e^{-y}\, y^{2N-2} \sum_{n=1}^{\infty} \frac{1}{n^{2N-2}}
\nonumber\\
& \times \; \frac{1}{(y^2 -4n^2 \pi^2 |z|^2)}  +  (-1)^M \Bigl(  \lfloor M/2 \rfloor  + \frac{1-(-1)^M}{2} \Bigr) 2 \pi |z|-
\frac{1}{2} \; \ln \! \left( 1- e^{-2 \pi |z|} \right) \; . 
\label{forty}\end{align} 
As stated previously, the regularized value is unique for all values of $z$. Hence the rhs of the above result is equal to the lhs 
of Equivalence\ (\ref{eleven}), which means in turn that
\begin{align}
\ln \Gamma(z) & = \Bigl( z- \frac{1}{2} \Bigr) \ln z -z + \frac{1}{2} \,\ln (2 \pi) 
+ z \sum_{k=1}^{N-1} \frac{(-1)^k}{(2z)^{2k}} \, \Gamma(2k-1)\,c_k(1) + 2  \Bigl( \frac{1}{4 \pi^2 |z|^2} \Bigr)^{N-1} 
\nonumber\\
& \times \;\; z \, P \int_0^{\infty} dy \, e^{-y}\, y^{2N-2} \sum_{n=1}^{\infty} \frac{1}{n^{2N-2}(y^2 -4n^2 \pi^2 |z|^2)} 
+ (-1)^M \Bigl( \lfloor M/2 \rfloor 
\nonumber\\
& + \;\; \frac{1-(-1)^M}{2} \Bigr)  \, 2 \pi |z| -  \frac{1}{2} \; \ln \! \left( 1- e^{-2 \pi |z|} \right)   \Bigr] \;.
\label{fortyone}\end{align} 
where the upper- and lower-signed versions become one form that is only valid for $\theta \!=\! \pm (M+1/2) \pi$.
Hence we have arrived at the result in Thm.\ 2.1 with the remainder given by Eq.\ (\ref{twentythree}) and
the Stokes discontinuity term given by Eq.\ (\ref{twentyfour}).
This completes the proof of Theorem\ 2.1. 

It should be stressed that the remainder in Theorem\ 2.1 is conceptually different from the remainder term discussed in 
standard Poincar$\acute{\rm e}$ asymptotics. In the latter prescription the remainder is occasionally bounded well before
the optimal point of truncation, but more frequently, it is left open with the introduction of either the Landau gauge symbol 
$\mathcal{O}()$ or $+ \dots$.  That is, an expression like Eq.\ (\ref{twentyone}) would typically be expressed as
\begin{eqnarray}
\ln \Gamma(z)  =  \Bigl( z- \frac{1}{2} \Bigr) \ln z - z + \frac{1}{2} \,\ln (2 \pi) - \frac{c_1(1)}{4z}+ \frac{c_2(1)}{8 z^3}
-\frac{3c_3(1)}{8 z^5} + \LARGE{\mathcal{O}} \Bigl(\frac{1}{z^7}\Bigr) \;.  
\label{fortytwo}\end{eqnarray}
Moreover, by introducing $c_1(1) \!=\!-1/3$, $c_2(1) \!=\! -1/45$, $c_3(1) \!=\!-2/945$ and $c_4(1) \!=\! -1/4725$, into the 
above result, we obtain the Stirling approximation or  ``Eq."\ (\ref{two}). For real values of $z$ the above result is referred 
to as a large $z$ or $z \to \infty$ expansion with the limit point situated at infinity. For $z$ complex, it becomes a large $|z|$ 
expansion. In those cases, where the Landau gauge symbol is dropped, a tilde often replaces the equals sign. Whichever case is 
considered, it means that the later terms in the truncated power series or remainder are neglected due to their eventual divergence. 
On the other hand, the remainder term given by  Eqs.\ (\ref{twentytwo}) or (\ref{twentythree}) is valid for any value of the 
truncation parameter $N$. Although the truncated power series still diverges as $N \to \infty$, when it is combined with the 
corresponding value of its remainder, it yields the same value as when the truncation parameter is equal to unity. That is, the 
divergence in the truncated series as $N$ increases is counterbalanced by the remainder diverging in the opposite sense. We shall
observe this behaviour in the following section.

We can show that the truncated series can be cancelled by the remainder leaving the remainder that one obtains when 
$N \!=\! 1$. To accomplish this, the following identity is required
\begin{eqnarray}
\frac{1}{n^k\, (n+a)} = \sum_{j=1}^k \frac{(-1)^{j+1}}{a^j\, n^{k-j+1}} + \frac{(-1)^k}{a^k\,(n+a)}\;.
\label{fortythree}\end{eqnarray}
Substituting $a$, $n$ and $k$ in the identity by $(y/2\pi z)^2$, $n^2$ and $N \!-\! 1$, respectively, and introducing the rhs into the
expression for the remainder given by Eq.\ (\ref{twentytwo}), one obtains
\begin{align}
R^{SS}_N(z) & = -2z \Bigl( - \frac{1}{4 \pi^2 z^2} \Bigr)^N  \int_0^{\infty} dy \; y^{2N-2}\, e^{-y} \sum_{n=1}^{\infty} 
\left( \sum_{j=1}^{N-1} \frac{(-1)^{j+1}}{(y/2\pi z)^{2j}\,n^{2N-2j}} \right.
\nonumber\\
& + \left. \;\; \frac{(-1)^{N+1}}{(y/2\pi z)^{2N-2}} \; \frac{1}{n^2+ y^2/4 \pi^2 z^2} \right) \;.
\label{fortythreea}\end{align}

In the first term with the summation over $j$ we replace the sum over $n$ by the zeta function. Consequently, we can introduce Eq.\ 
(\ref{fourteen}), which gives
\begin{align}
R^{SS}_N(z) & = z \int_0^{\infty} dy \; e^{-y} \left( \sum_{j=1}^{N-1} \frac{(-1)^{N-j+1}}{(2 z)^{2N-2j}} \, c_{N-j}(1) \, 
y^{2N-2j-2} \right.
\nonumber\\
& + \left. \;\; 2\sum_{n=1}^{\infty}\frac{(-1)^{N+1}}{(2\pi z)^2} \; \frac{1}{n^2+ y^2/4 \pi^2 z^2} \right) \;.
\label{fortyfour}\end{align}
Setting $j \!=\!N \!-\!j$ and then performing the integration in the first on the rhs of the above result, we arrive at
\begin{eqnarray}
R^{SS}_N(z) = -z \sum_{j=1}^{N-1} \frac{(-1)^j}{(2z)^j} \; c_j(1)\, \Gamma(2j-1) + \frac{1}{2 \pi^2 z} \int_0^{\infty} dy\;
\sum_{n=1}^{\infty} \frac{e^{-y}}{n^2+ y^2/4 \pi^2 z^2} \;.
\label{fortyfive}\end{eqnarray}
Hence we see that the remainder is composed of the truncated series, which cancels the truncated series in Eq.\ (\ref{twentyone}),
and the $N \!=\! 1$ value for Eq.\ (\ref{twentytwo}) or $R^{SS}_1(z)$. Moreover, with the aid of No.\ 1.421(4) in Ref.\ \cite{gra94} we 
obtain the following integral representation for $R^{SS}_1(z)$:
\begin{eqnarray}
R^{SS}_1(z) = \frac{1}{2} \int_0^{\infty} dy\; y^{-1} \,e^{-y} \Bigl( {\rm cth}\Bigl( \frac{y}{2z} \Bigr) - \frac{2z}{y} \Bigr) \;.
\label{fortysix}\end{eqnarray}
Introducing the above result into Eq.\ (\ref{twentyone}) yields 
\begin{eqnarray} 
\ln \Gamma(z)  =  \Bigl( z- \frac{1}{2} \Bigr) \ln z - z + \frac{1}{2} \,\ln (2 \pi) + 
\frac{1}{2} \int_0^{\infty} dy\, y^{-1} e^{-y} \Bigl( {\rm cth}\Bigl( \frac{y}{2z} \Bigr) - \frac{2z}{y} \Bigr) 
+ SD^{SS}_M(z)  .
\label{fortyseven}\end{eqnarray}
If we differentiate Eq.\ (\ref{fortyseven}) w.r.t.\ $z$ and note that
\begin{eqnarray}
\frac{d}{dz} \; {\rm cth} \Bigl( \frac{y}{2z}\Bigr)= - \frac{y}{z} \frac{d}{dy} \; {\rm cth} \Bigl( \frac{y}{2z}\Bigr)\;\;,
\label{fortyeight}\end{eqnarray}
then after integrating by parts we arrive at
\begin{align}
\psi(z) & = \ln z -\frac{1}{2z} + \int_0^{\infty} dy\; e^{-y} \Bigl( \frac{1}{2z} \; {\rm cth}\Bigl( \frac{y}{2z}\Bigr)
- \frac{1}{y} \Bigr)  
\nonumber\\
& \mp  \;\; \pi i \Bigl( 2 \lfloor M/2 \rfloor - \frac{1-(-1)^M}{e^{\mp 2 i \pi z }-1}\Bigr) \;.
\label{fortynine}\end{align}
In the above equation $\psi(z)$ represents the digamma function or $d\ln \Gamma(z)/dz$. Although Eq.\ (\ref{fortynine}) has been
obtained by regularizing an asymptotic series, it agrees with No.\ 3.554(4) in Ref.\ \cite{gra94} when $\Re\, z \!>\!0$ in
the principal branch of the complex plane, i.e.\ for $|\theta|< \pi/2$. Outside this sector, however, the Stokes discontinuity term 
appears making Eq.\ (\ref{fortynine}) more general than the result in Ref.\ \cite{gra94}. Furthermore, if we set $z$ equal to unity,
then we see that the Stokes discontinuity or the second term on the rhs of Eq.\ (\ref{fortynine}) vanishes and we are left with
\begin{eqnarray}
\gamma= \frac{1}{2} +\int_0^{\infty} dy \; e^{-y} \Bigl( \frac{1}{2} \; {\rm cth} \Bigl( \frac{y}{2} \Bigr) - \frac{1}{y}
\Bigr) \;,
\label{fifty}\end{eqnarray}
where $\gamma$ denotes Euler's constant. 

In Ref.\ \cite{kow10} it is found that
\begin{eqnarray}
\gamma= \sum_{k=1}^{\infty} \frac{(-1)^{k+1}}{k}\; A_k \;,
\label{fiftyone}\end{eqnarray}
where $A_0 \!=\!1$ and the remaining $A_k$ are positive fractions given by
\begin{eqnarray}
A_k = \frac{(-1)^k}{k!} \int_0^1 dt \; \frac{\Gamma(k+t-1)}{\Gamma(t-1)} \;.
\label{fiftytwo}\end{eqnarray}
These coefficients are referred to as the reciprocal logarithm numbers, but according to p.\ 137 of Ref.\ \cite{ape08}, their 
absolute values are known as either the Gregory or Cauchy numbers. They also satisfy the following recurrence relation:
\begin{eqnarray}
A_k = \sum_{j=0}^{k-1} \frac{(-1)^{k-j+1}}{k-j+1} \; A_j \;,
\label{fiftythree}\end{eqnarray}
whereupon we see that $A_1=1/2$, $A_2= -1/12$, $A_3=1/24$, etc. Eq.\ (\ref{fiftyone}) is referred to as Hurst's formula in 
Ref.\ \cite{kow10}, although since then it has been revealed that Kluyver calculated several of the leading terms \cite{kly24}.
As a consequence, the integral in Eq.\ (\ref{fifty}) can be written as
\begin{eqnarray}
\int_0^{\infty} dy \; e^{-y} \Bigl( \frac{1}{2} \; {\rm cth} \Bigl( \frac{y}{2} \Bigr) - \frac{1}{y} \Bigr)  = \sum_{k=2}^{\infty} 
\frac{(-1)^{k+1}}{k}\, A_k \;.
\label{fiftyfour}\end{eqnarray}

The analysis resulting in Eqs.\ (\ref{fortyfour})-(\ref{fortysix}) can be applied when $z$ lies on a Stokes line, but now
we require Eq.\ (\ref{twentythree}). Once again, we use the identity given by Eq.\ (\ref{fortythree}) with the same values
for $n$ and $k$, but with $a$ set equal to $-1/4\pi^2 |z|^2$. We also replace $x$ by $ix$ in No.\ 1.421(4) of Ref.\
\cite{gra94} so that it becomes
\begin{eqnarray}
\cot x = \frac{1}{x} - \frac{2x}{\pi^2} \sum_{k=1}^{\infty} \frac{1}{k^2-x^2/\pi^2} \;.
\label{fiftyfive}\end{eqnarray}
As a result, the remainder along the Stokes lines of $\theta =\pm (M+1/2) \pi$ is found to be given by
\begin{align}
R^{SL}_N(z)& = -z \sum_{k=1}^{N-1} \frac{\Gamma(2k-1)}{(2 |z|)^{2k}}\; c_k(1) \mp \frac{(-1)^M \, i}{2} \; P \int_{0}^{\infty} dy 
\; y^{-1}\, e^{-y} 
\nonumber\\
& \times \;\; \Bigl( \cot \Bigl( \frac{y}{2 |z|}\Bigr)- \frac{2 |z|}{y} \Bigr) \;. 
\label{fiftysix}\end{align}
Therefore, when the above result is introduced into Eq.\ (\ref{twentyone}), the terms involving the truncated series cancel
each other and we are left with the integral, which represents $R_1(z)$ for $z$ lying on a Stokes line.

\section{Numerical Analysis}
In the previous section it was shown that despite being composed of an asymptotic series and thus, being divergent in certain 
regions of the complex plane, Stirling's approximation in its entirety can be regularized, thereby yielding an equation for 
the logarithm of the gamma function as given by Eq.\ (\ref{twentyone}). With the aid of the material in Ref.\ \cite{kow09}, 
exact expressions were derived for the remainder of this famous approximation and the Stokes discontinuity term over 
all Stokes sectors and lines. Although these results were proved, one still cannot be sure that these results will yield exact values 
for the gamma function unless an effective numerical study is carried out. This is because proofs in standard Poincar$\acute{\rm e}$ 
asymptotics for the most part invoke such symbols as $\sim$, $\approx$, $\mathcal{O}()$, $+ \dots$, $\leq$, and $\geq$. The introduction 
of these ``tools" obscures the most interesting and challenging problem in asymptotics --- the behaviour/evaluation of the 
remainder. Consequently, asymptotic results such as ``Eq."\ (\ref{two}) are often limited in accuracy and are only valid over 
a narrow, but vague, range. These hardly represent the ideal qualities of an exact science as mathematics and as a result, standard 
Poincar$\acute{\rm e}$ asymptotics has often been subjected to much criticism and ridicule, particularly from pure mathematicians. 
Moreover, when it is realized that the remainder possesses an infinity in certain regions of the complex plane, neglecting the late 
terms in a proof is simply invalid. In such cases the proof can only be regarded as fallacious. On the other hand, a numerical study, 
if performed with proper attention to detail and without bias, is often more valuable and incisive, a fact that is generally overlooked 
by the asymptotics community. In addition, it has the potential to uncover new properties about the original function.

It should also be noted that one does not require a vast number of values of $z$ to carry out an effective numerical investigation. 
We have seen that the results in Theorem\ 2.1 only change form according to specific sectors or rays in the complex plane. Within 
each sector or on each line the results behave uniformly with respect to $z$, i.e., there are no singularities that exclude particular 
values of $z$. Hence only a few values of $|z|$ are necessary for conducting a proper numerical study. In fact, to demonstrate that the 
remainder behaves according to standard Poincar$\acute{\rm e}$ asymptotics, we need only two values of $|z|$: a relatively large one, 
where it is valid to truncate the asymptotic series in Eq.\ (\ref{sixteen}), and a small one, where truncation breaks down. Then the more
important issue is to consider a range of values for both the truncation parameter $N$ and the argument or phase of $z$, viz. ${\rm arg}\,z$ 
or $\theta$ as it has been denoted in this work. Furthermore, selecting extremely large/small values of $|z|$ may result in overflow or underflow 
problems in the numerical computations, thereby creating the misleading impression that the results in Theorem\ 2.1 are not correct rather 
than being a deficiency of the computing system. Since the variable in the asymptotic series is $1/(2 n \pi z)^2$ with $n$ ranging from 
unity to infinity, which follows once the Dirichlet series form for the Riemann zeta function is introduced into Eq.\ (\ref{fourteen}), 
a value of 3 represents a ``large value" for $|z|$, courtesy of the $2n \pi$ factor. As for a small value of $|z|$, a value of 1/10 will 
suffice, which we shall see is sufficiently small to ensure that there is no optimal point of truncation for small values of $n$. 

The optimal point of truncation can be estimated by realizing that it occurs when the successive terms in an asymptotic 
series begin to dominate the preceding terms. Specifically, it occurs at the value of $k$ where the $k+1$-th term is greater
than the $k$-th term in $S(z)$ or Eq.\ (\ref{sixteen}). This is given by
\begin{eqnarray}
\left| \frac{2k(2k-1)}{(2z)^2} \; \frac{c_{k+1}(1)}{c_k(1)} \right|=  \left| \frac{2k(2k-1)}{(2 \pi z)^2} \; \frac{\zeta(2k+2)}{\zeta(2k)} 
\right|  \approx 1 \;.
\label{fiftyseven}\end{eqnarray} 
The ratio of the Riemann zeta functions in the above estimate is close to unity. Consequently, we see that the optimal point $N_{OP}$ 
occurs when it is roughly equal to $\pi |z|$. For $|z| \!=\! 3$, this means that the optimal point of truncation will
be in the vicinity of $N_{OP}=10$, while for $|z| \!=\! 1/10$, there is no optimal point of truncation, i.e. $N_{OP}=0$, since the 
truncation parameter must be greater than unity. In this instance the first or leading term of the asymptotic series will yield 
the ``closest" value to the actual value of $\ln \Gamma(z)$, but it will be by no means accurate. It is for such values of $|z|$ 
where no optimal point of truncation exists that standard Poincar$\acute{\rm e}$ asymptotics breaks down. On the other hand, the larger 
$N_{OP}$ is, the more accurate truncation of the asymptotic series becomes.

Typically, when a software package such as Mathematica \cite{wol92} is used to determine values for a special function such as 
the gamma function, it does this only over the principal branch of the complex plane for its variable. This means that for the 
proposed numerical study ${\rm arg}\,z$ or $\theta$ must line in the interval $(-\pi,\pi]$. As a result, the numerical verification 
of Eq.\ (\ref{twentyone}) will only involve three Stokes sectors covering the principal branch of the complex plane, viz.\ 
$-3 \pi/2 < \theta < -\pi/2$, $-\pi/2 < \theta <\pi/2$, and $\pi/2 < \theta < 3 \pi/2$, and the two Stokes lines at 
$\theta = \pm \pi/2$. That is, we can only test the $M=0$ and $M= \pm 1$ results in Theorem\ 2.1. If we denote the truncated sum 
in Eq.\ (\ref{twentyone}) by $TS_N(z)$ such that
\begin{eqnarray}
TS_N(z) = z \sum_{k=1}^{N-1} \frac{(-1)^k}{(2z)^{2k}}\; \Gamma(2k-1)\,c_k(1) \;,
\label{fiftysevena}\end{eqnarray}
then the results that we need to verify over the principal branch for $z$ can be expressed as
\begin{eqnarray}
\ln \Gamma(z) = \begin{cases}  F(z) + TS_N(z) + R^{SS}_N(z) + SD^{SS,U}_{1}(z) \;,\quad   & \pi/2 < \theta \leq \pi\;, \cr
F(z) + TS_N(z) +R^{SL}_N(z) +SD^{SL}_{0}(z) \;, \quad & \theta = \pi/2 \;, \cr
F(z) + TS_N(z) +R^{SS}_N(z) \;, \quad & -\pi/2 < \theta < \pi/2 \;, \cr
F(z) + TS_N(z) +R^{SL}_N(z) + SD^{SL}_{0}(z) \;, \quad & \theta = -\pi/2 \;, \cr
F(z) + TS_N(z) + R^{SS}_N(z) + SD^{SS,L}_{1}(z) \;, \quad & -\pi < \theta < -\pi/2 \;,
\end{cases} 
\label{fiftyeight}\end{eqnarray}
where $F(z)$ represents the leading terms in Eq.\ (\ref{twentyone}). Specifically, $F(z)$ is given by
\begin{eqnarray}
F(z) = \Bigl( z- \frac{1}{2} \Bigr) \ln z- z+\frac{1}{2} \, \ln( 2\pi)\;.
\label{fiftynine}\end{eqnarray}
When compared with ``Eq." (\ref{two}), we see that the leading terms are basically the terms in Stirling's approximation for the gamma function. 
In the above results the superscripts of U and L have been introduced into the Stokes discontinuity terms for the Stokes sectors to indicate 
the upper- and lower-signed versions of Eq.\ (\ref{twentytwoa}). Although equal to zero, we shall refer to the Stokes discontinuity term for the 
third asymptotic form as $SD^{SS}_0(z)$. 

If we put $N \!=\!4$ in the central result of Eq.\ (\ref{fiftyeight}) and neglect the final term or remainder, then we arrive at Eq.\ 
(\ref{fortytwo}). Moreover, in Eq.\ (\ref{fiftyeight}) the remainder terms are given by
\begin{align}
R^{SS}_N(z) & =  \frac{2\,(-1)^{N+1} \,z}{(2\pi z)^{2N-2}} \sum_{n=1}^{\infty} \frac{1}{n^{2N-2}}\int_0^{\infty} dy \; 
\frac{y^{2N-2}\, e^{-y}}{\left( y^2 + 4 \pi^2 n^2 z^2 \right)} \;,
\label{sixty}\end{align}
and 
\begin{align}
R^{SL}_N(z) & =  \frac{2\,z}{(2\pi |z|)^{2N-2}} \sum_{n=1}^{\infty} \frac{1}{n^{2N-2}} P \int_0^{\infty} dy \; 
\frac{y^{2N-2}\, e^{-y}}{\left( y^2 - 4 \pi^2 n^2 |z|^2 \right)} \;,
\label{sixtya}\end{align}
while the Stokes discontinuity terms are given by
\begin{eqnarray}
SD^{SS}_{1}(z) = -\ln \Bigl( 1- e^{\pm 2\pi z i} \Bigr) \;,
\label{sixtyone}\end{eqnarray}
and
\begin{eqnarray}
SD^{SL}_{0}(z) = - \frac{1}{2} \ln \Bigl( 1- e^{- 2\pi |z| } \Bigr) \;.
\label{sixtytwo}\end{eqnarray}

In order to proceed with the numerical investigation, we need to consider the results over the Stokes sectors separately
from those applicable to the Stokes lines at $ \theta \!=\! \pm \pi/2$.  This is because: (1) the latter involve
the evaluation of the Cauchy principal value and (2) the Stokes discontinuity terms possess a factor of $1/2$
compared with zero when $|\theta|<\pi/2$ and unity when $|\theta|> \pi/2$. Therefore, to obtain values of $\ln \Gamma(z)$
using the above results, we shall require two different programs or modules: one, where the standard numerical integration 
routine called NIntegrate in Mathematica is invoked and another, where the NIntegrate routine is required to evaluate the Cauchy 
principal value and half the Stokes discontinuity term. As we shall see shortly, the second module is far more computationally 
intensive and thus, takes much longer to execute. 

\begin{figure}
\begin{center}
\includegraphics{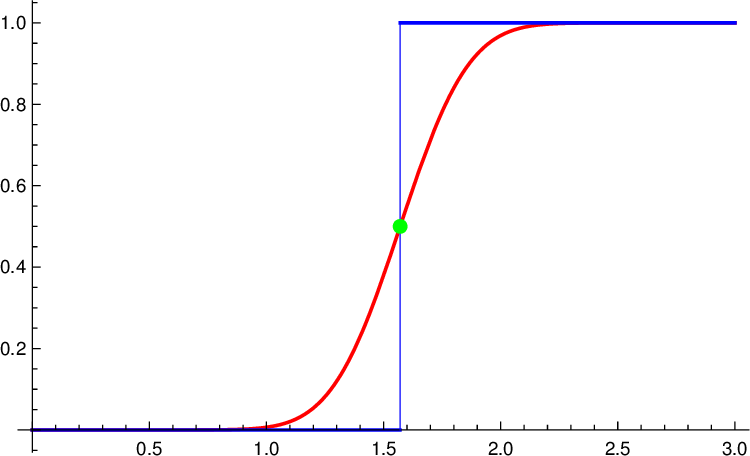}
\end{center}
\caption{Graph of the Stokes multiplier $S^{+}$ given by ``Eq." (\ref{sixtyseven}) with $|z| =3$ versus $\theta$.}
\label{figone}
\end{figure}

When $\theta > 0$, we can combine the Stokes discontinuity terms into one expression, which we denote by $SD^{+}(z)$. Hence
the discontinuity terms can be expressed as
\begin{eqnarray}
SD^{+}(z) = - S^{+} \ln \Bigl( 1- e^{2 \pi i z} \Bigr) ,
\label{sixtythree}\end{eqnarray} 
where the factor $S^{+}$ is given by
\begin{eqnarray}
S^{+} = \begin{cases} 
1 \;, & \; \pi/2 < \theta \leq \pi \;, \cr
1/2 \;, & \theta = \pi/2 \;, \cr
0 \; , & -\pi/2 < \theta < \pi/2 \;, 
\end{cases}
\label{sixtyfour}\end{eqnarray}
Similarly, we can denote the Stokes discontinuity terms in the lower half of the principal branch by $SD^{-}(z)$ and express
them in terms of another factor $S^{-}$. Hence we arrive at
\begin{eqnarray}
SD^{-}(z) = - S^{-} \ln \Bigl( 1- e^{-2 \pi i z} \Bigr) \;,
\label{sixtyfive}\end{eqnarray} 
with $S^{-}$ given by
\begin{eqnarray}
S^{-} = \begin{cases} 
0 \;, & \; -\pi/2 < \theta < \pi/2 \;, \cr
1/2 \;, & \theta = -\pi/2 \;, \cr
1 \; , & -\pi < \theta < -\pi/2 \;, 
\end{cases}  
\label{sixtysix}\end{eqnarray}
In the literature $S^{+}$ and $S^{-}$ are known as Stokes multipliers \cite{par01}. From the preceding analysis we see that 
they are discontinuous, which is in accordance with the conventional view of the Stokes phenomenon. However, as a result of
Ref.\ \cite{ber89}, an alternative view of Stokes phenomenon has arisen where these multipliers are no longer regarded 
as discontinuous step-functions, but experience a smooth and rapid transition from zero to unity, equalling $1/2$ when $z$ 
lies on a Stokes line. This has become known as Stokes smoothing, although it should be really be called Berry 
smoothing of the Stokes phenomenon since at no stage did Stokes ever regard the multipliers as being smooth \cite{sto04}. 
According to the approximate theory initially developed by Berry and then made more ``rigorous" by Olver \cite{olv90}, 
the multiplier can be expressed in terms of the error function ${\rm erf}(z)$. Shortly afterwards, Berry \cite{ber91} 
and Paris and Wood \cite{par92} derived an approximate form for the Stokes multipliers of $\ln \Gamma(z)$. These are given 
by 
\begin{eqnarray}
S^{\pm}(z) \sim  \frac{1}{2} \pm \frac{1}{2} \, {\rm erf} \Bigl((\theta \pm \pi/2) \sqrt{\pi |z|} \Bigr)\; .
\label{sixtyseven}\end{eqnarray}
A graph of this result for $|z| =3$ versus $\theta$ is presented in Fig.\ \ref{figone}. Here we see that for $\theta <1$,
the Stokes multiplier is virtually zero, while for $\theta >2$, it is almost equal to unity, which is consistent with the
conventional view of the Stokes phenomenon. In between, however, the smoothing as postulated by Berry and Olver is expected 
to occur with the greatest deviation from the original step-function occurring in the vicinity of the Stokes line at $\theta = 
\pi/2$. Therefore, the results in Theorem\ 2.1 are not expected to yield accurate values of $\ln \Gamma(z)$, especially for 
$13\pi/32 < \theta <17\pi/32$, if smoothing occurs. Nonetheless, the adherents of this theory have not to this day provided 
one numerical demonstration as to whether ``Stokes smoothing" actually does occur or whether the conventional view still holds. 
Here, however, we can establish whichever view is correct by evaluating $\ln \Gamma(z)$ for deviations of $\theta$ above and 
below $\pi/2$ so that they lie in the range of $(13\pi/32,17\pi/32)$. That is, where the smoothing is expected to be at its 
most pronounced. If the above results for the Stokes discontinuity terms are unable to provide exact values of $\ln \Gamma(z)$, 
then we know that the conventional view of the Stokes phenomenon is not valid and that smoothing is a viable alternative.  
It should also be noted that the above result for the Stokes multiplier only applies to large values of $|z|$, while for $|z|<5$
there is no accurate representation. Thus, the concept has only ever been applied to large values of $|z|$, where truncating the 
asymptotic expansion after a few terms will yield very accurate values for $\ln \Gamma(z)$. The numerical studies presented here
are aimed at small values of $|z|$, where it is no longer possible to obtain accurate values of the function and where the
Stokes discontinuity term cannot be disregarded in order to obtain exact values of $\ln \Gamma(z)$.

Before we carry out the investigation into Stokes smoothing, we need to show that the remainder terms in Eq.\ 
(\ref{fiftyeight}) do in fact behave typically for an asymptotic expansion. That is, we need to show that for 
large values of $|z|$, the remainder can be neglected to yield accurate, but still approximate, values of $\ln \Gamma(z)$ 
up to and not very far from the optimal point of truncation, while for small values of $|z|$, it is simply not valid
to neglect the remainder. For this demonstration we do not require the Stokes discontinuity terms. Therefore, we 
shall concentrate on the asymptotic form for $|\theta| <\pi/2$. This includes $\theta=0$, which is the simplest case to 
study because it does not involve complex arithmetic. Later, when we study the values produced for all the Stokes sectors 
and on the Stokes lines in the principal branch, we shall consider non-zero values of $\theta$.
  
From Eq.\ (\ref{sixty}) we see that the evaluation of the remainder involves two computationally intensive tasks. The first 
is the infinite sum over the integers $n$, which has arisen because there is an infinite number of singularities lying on
each Stokes line. The second issue is the numerical integration of the exponential integral. The latter can be avoided by 
expressing the integral in terms of the incomplete gamma function. That is, by decomposing the denominator into 
partial fractions and introducing No.\ 3.383(10) from Ref.\ \cite{gra94}, we find that the remainder can be expressed as
\begin{align}
R_N^{SS}(z) & =\frac{ \Gamma(2N-1)}{2 \pi i} \sum_{n=1}^{\infty} \frac{1}{n} \Bigl( 
e^{-2 \pi n z i} \, \Gamma(2-2N, -2\pi n z i) 
\nonumber\\
& - \;\; e^{2 \pi n z i} \, \Gamma(2-2N, 2 \pi n z i ) \Bigr) \; ,
\label{sixtyeight}\end{align}
where $|\theta| < \pi/2$. Although not explicitly given, the above result can be inferred by combining Eqs.\ (4.3), (4.10)
and (4.11) in Ref.\ \cite{par92}. 

Eq.\ (\ref{sixtyeight}) is only of use if we have a mathematical software package that has the capability of evaluating the 
incomplete gamma function to very high precision. Fortunately, Mathematica \cite{wol92} does have this capability, but we shall 
also consider Eq.\ (\ref{sixty}) when evaluating $\ln \Gamma(z)$ since it only requires a numerical integration routine, which is 
more easily constructed when one does not have access to a mathematical software package. Furthermore, Eq.\ (\ref{sixtya}) 
represents the continuation of Eq.\ (\ref{sixty}) to the Stokes lines, while we do not know how Eq.\ (\ref{sixtyeight}) can be 
extended beyond the primary Stokes sector. 
  
\begin{table}
\small
\centering
\begin{tabular}{|c|c|c|} \hline
$N$ & Quantity & Value  \\ \hline
& $F(3)$ &       0.66546925487494697026844282871193190148012386819465 \\ \hline
& TS &           0.02777777777777777777777777777777777777777777777777 \\
2 & $R^{SS}_2(3)$ &   -0.0000998520927794385973038298896926468609453577911 \\
& Total &        0.69314718055994530944891677660001703239695628818129 \\ \hline
& TS &           0.02767489711934156378600823045267489711934156378600 \\
3 & $R^{SS}_3(3)$ &   3.028565656775362781062293505563582854900697488 $\times 10^{-6}$ \\
& Total &        0.69314718055994530941723212145811236218232033267815 \\ \hline
& TS &           0.02767789100093626842598036013673873756178282927254 \\
5 & $R^{SS}_5(3)$ &   3.468406207072280893260950592903359343689305700 $ \times 10^{-8}$ \\
& Total &        0.69314718055994530941723212145817656807550013436025 \\ \hline
& TS &           0.02767792490305420773799675002229807193246527808017 \\
9 & $R^{SS}_9(3)$ &   7.819441314107925427239465946629109880854200684 $\times 10^{-10}$  \\
& Total &        0.69314718055994530941723212145817656807550013436025  \\ \hline
& TS &           0.02767792629413478403268961923401255626286988423163 \\
10 & $R^{SS}_{10}(3)$ & -6.091364448839003264877678896674936180660392163 $ \times 10^{-10}$ \\
& Total &        0.69314718055994530941723212145817656807550013436025 \\ \hline
& TS &           0.02767792509609780488374471454379088783178430267081 \\
11 & $R^{SS}_{11}(3)$ &  5.889005342650445782024537787635919634947854418 $\times 10^{-10}$ \\
& Total &        0.69314718055994530941723212145817656807550013436025 \\ \hline
& TS &           0.02767792637739909405287985684200177174299855050138 \\
12 & $R^{SS}_{12}(3)$ & -6.924007549040905640957571051476222843357840714 $\times 10^{-10}$ \\
& Total &        0.69314718055994530941723212145817656807550013436025 \\ \hline
& TS &           0.02767792256451067899110318500595297397521206328541 \\
15 & $R^{SS}_{15}(3)$ &  3.120487660157686107740291692620164202880179513  $\times 10^{-9}$ \\
& Total &        0.69314718055994530941723212145817656807550013436025 \\ \hline
& TS &           0.02767853067626590045913727597678865028583493409103 \\
20 & $R^{SS}_{20}(3)$ & -6.0499126756131034798323054398369045866792543896 $\times 10^{-7}$ \\
& Total &        0.69314718055994530941723212145817656807550013436025 \\ \hline
& TS &           41.2834736138079254966213754129139774958755379575621 \\
30 & $R^{SS}_{30}(3)$ &  -41.255795688122927157472586120167732829280161691396 \\
& Total &        0.69314718055994530941723212145817656807550013436025 \\ \hline
& TS &           6.0039864088710184849557428450939638222762809177 $\times 10^{25}$ \\
50 & $R^{SS}_{50}(3)$ &  -6.003986408871018484955742842326171253776447002 $\times 10^{25}$\\
& Total &        0.69314718055994530941723212145817656807550013436025 \\ \hline
& $\ln \Gamma(3)$ & 0.69314718055994530941723212145817656807550013436025 \\ \hline
\end{tabular}
\normalsize
\caption{Determination of $\ln\Gamma(3)$ via Eq.\ (\ref{sixtyeight}) for various values of the truncation parameter, $N$}
\label{tab1}
\end{table}

The first program presented in the appendix is a Mathematica module that evaluates $\ln \Gamma(z)$ using Eq.\ (\ref{sixtyeight})
for the remainder. Instead of evaluating the sum over $n$ to infinity, the module calculates all the terms
to the value assigned to the variable called limit. In order to achieve great accuracy, limit has been set equal to
$10^5$. That is, we are essentially truncating the series over $n$ to $10^5$, which will also be the case when comparing
the results obtained via Eqs.\ (\ref{sixty}) and (\ref{sixtya}). Although this limit will be sufficient for our purposes,
more astonishing results can be obtained by raising this value higher, but it comes at the expense of computational
time, which will be discussed in the presentation of the numerical data.  

Table\ \ref{tab1} displays a small sample of the results obtained by running the first module in the appendix on
both a Sony VAIO laptop with 2 Gb of RAM and a SunFire X4600M2 server with 64 Gb of RAM for $z \!=\!3$ and
various values of the truncation parameter, $N$. The surprising feature of these calculations was that the laptop 
was able to perform an individual calculation significantly faster than the SunFire alpha server, although the 
latter was to carry a much greater number of calculations simultaneously. Whilst the values were often printed out to more 
than 50 decimal places in full form, we should regard 50 decimal places as the upper limit since the module calculates all the 
values to 50 decimal places. This is due to N[,50] appearing in each print statement, so that as the values are printed 
out, they often appear with the suffix `50. 

Another feature of the results is that they are all real, which is expected since $\ln \Gamma(z)$ is real and positive
for all real values of $z$ greater than unity. In actual fact, Mathematica did print out a tiny imaginary part
with each value, but in nearly all cases it was zero to the first 50 or even more decimal places. Hence the imaginary
contributions have been discarded. In certain situations, which will be discussed later, there were values of zero, but with 
a suffix less than 50. The appearance of these tiny imaginary values is a good sign because it gives an indication of 
the numerical error involved in the calculations.

The first column in Table\ \ref{tab1} displays the values of the truncation parameter $N$, which begin with $N \!=\!2$ and
and end with $N \!=\! 50$, the latter well beyond the optimal point of truncation or $N_{OP}$. Four consecutive values between 
$N \!=\! 9$ and $N \!=\! 12$ have been presented because it is expected that one of these values will be the optimal point of 
truncation according to the discussion below ``Eq."\ (\ref{sixtyseven}). The second row in the table gives the value of $F(z)$ 
for $z \!=\! 3$. That is, this value represents Stirling's approximation to $\ln \Gamma(z)$, which we see is close the
actual value of $\ln \Gamma(3)$ or $\ln 2$ appearing in the bottom row of Table\ \ref{tab1}. Since these values are invariant,
they only appear once in the table. To avoid any possible problems arising out of the multivaluedness of the Log function
in Mathematica \cite{wol92}, the modulus and argument of $z$ are separate inputs in the first Mathematica module in the appendix.
It should also be noted that there has been no rounding-off introduced into any of the results in the table.

Although far more values of the truncation parameter were considered, the eleven results in Table\ \ref{tab1} are deemed
sufficient to demonstrate that the remainder given by Eq.\ (\ref{twentytwo}) possesses the accepted behaviour of the
remainder in standard Poincar$\acute{\rm e}$ asymptotics. For each value of $N$ there are three rows. The first row 
labelled TS and represented by e1 in the first Mathematica module in the appendix is the value of the truncated sum 
given by Eq.\ (\ref{sixty}). The second value labelled $R^{SS}_N(3)$ and represented by rem in the same Mathematica 
module is the value of the remainder with limit in the Do loop set to $10^5$ as mentioned previously. The third row 
labelled Sum displays the sum of $F(3)$, the truncated sum and the remainder. Therefore, according to Thm.\ 2.1, we expect 
this sum to give the value of $\ln \Gamma(z)$ for all values of the truncation parameter.

Most of the calculations performed by the SunFire server took between 5 and 8 hours to execute, whilst for the laptop they were 
generally a few hours shorter. The longest calculation was the final one in the table, which took 8.5 hours to complete, while
the shortest calculation turned out to be the first one, which took only 5.5 hours. There is, however, a method of speeding 
up the calculations so that they take a little over a minute to execute on the laptop. This can be accomplished by introducing 
into the Mathematica module the symbol N[expr,50], which attempts to give a numerical value for the expression expr to a 
precision of 50 digits. Therefore, whenever a value such as e2 or e3 is evaluated in the Mathematica module, the rhs should be 
wrapped around N[expr, 50] with expr equal to the quantity on the rhs. It may, however, affect the accuracy of the results, in 
particular the total value. Nevertheless, it had no effect on any of the results given for $N \!=\! 12$, but for those situations 
where the remainder begins to diverge, e.g. the final calculation, the summed value will only be accurate to a reduced number of 
decimal places. For example, the value obtained by summing all the contributions for $N \!=\! 50$ in the table was found to 
be accurate to the first 24 decimal places of $\ln \Gamma(3)$, whereas the result displayed in the table agrees with 
$\ln \Gamma(3)$ to over 50 decimal places. 
 
From the table we see that for $N \!=\! 2$, the truncated sum is found to equal $0.027\,777\cdots$, while the remainder 
$R^{SS}_2(3)$ is found to equal $-9.98\,529 \cdots \times 10^{-5}$. When these values are summed with $F(3)$, they yield 
a value that agrees with the value of $\ln \Gamma(3)$ to 19 decimal places, which is by no means as accurate as any of 
the other summed results in the table including the sum for $N \!=\! 50$, i.e. well beyond $N_{OP}$. In fact,
the remainder is of the order of $10^{25}$, which means that the first 25 places must be cancelled by the first 25 places
of the truncated sum in order to yield the decimal fraction for $\ln \Gamma(3)$. The reason why the $N \!=\! 2$ result
is not as accurate as the other results is that the factor of $n^{2N-2}$ in the denominator of Eq.\ (\ref{sixty}) still
has an effect on the calculation of the remainder for the small values of $N$ such 1 or 2. For these cases limit needs
to be increased substantially in order to make the remainder significantly more accurate. In other words, the neglected 
terms when limit is set to a large value can become relatively smaller as the truncation parameter increases.

It can also be seen that the remainder is smallest in magnitude when $N \!=\! 11$. Therefore, the point of optimal 
truncation for $z \!=\!3$ occurs at $N_{OP} \!=\! 11$, which compares favourably with our estimate of $N_{OP} \!=\! 10$ 
below ``Eq." (\ref{fiftyseven}). At $N_{OP}$ the sum of the values only differs from the actual 
value of $\ln \Gamma(3)$ at the fifty-third decimal place. We also see that for those values of the truncation parameter 
in the vicinity of $N_{OP}$ that there is little deterioration in the accuracy. For $N \!=\! 30$ and
$N \!=\! 50$, i.e. well past the optimal point of truncation, the remainder dominates, whereas for all the 
other calculations, it is extremely small, which is consistent with standard Poincar$\acute{\rm e}$ asymptotics.
As a consequence, for all but the last two calculations, $F(3)$ or Stirling's approximation represents the dominant 
contribution to $\ln \Gamma(3)$. In the last two calculations, both the truncated sum and remainder dominate, but the 
divergence in one is countered by the divergence in the other. Hence for $N \!=\! 50$, we see that both the remainder
and truncated sum are of the order of $10^{25}$, which means that at least the first the 26 decimal places of these
quantities must cancel each other so that Stirling's approximation becomes the dominant contribution again. This cancellation
of decimal places can only be achieved by a proper regularization of the remainder in an asymptotic expansion. It
was also responsible for creating an imaginary term that was zero to a reduced number of decimal places, viz.\ 23 places
instead of the 50 places described above.
  
We now turn our attention to the evaluation of the remainder when $z \!=\! 1/10$. This means that we are considering a small
value of $|z|$ in a large $|z|$ asymptotic expansion, which is unheard of in standard Poincar$\acute{\rm e}$ asymptotics.
It should also be pointed out such a value has yet to be tested by those claiming to carry out hyperasymptotic
evaluations of asymptotic expansions \cite{ber90,ber91,ber91a,ber89,par01}. In particular, Paris carries out a hyperasymptotic
evaluation of $\ln \Gamma(z)$ at the end of Ref.\ \cite{par11} using awkward Hadamard expansions for $\Omega(z)$. 
Depending on the number of levels he chooses, he displays results that are accurate at best to $10^{-45}$ for real values
of $z$. However, the results appearing in Table\ \ref{tab1} are far more accurate than his results despite the fact that the 
results in Table\ \ref{tab1} have been obtained for $z \!=\! 3$, whereas Paris put $z \!=\!8$. According to ``Eq."\ 
(\ref{fiftyseven}) $N_{OP}$ would equal 25 for this value of $z$. With such a high value for $N_{OP}$ we would easily obtain 
much greater accuracy than the results in Table\ \ref{tab1} without even having to consider a very large value for limit. The value of $z$ 
chosen by Paris is simply too large in order to gain any meaningful understanding as to whether his analysis has resulted in a 
marked improvement to standard Poincar$\acute{\rm e}$ asymptotics. Moreover, the main aim of hyperasymptotics should not be to 
extend the accuracy in the regions where an asymptotic expansion is already very strong, but to be able to obtain results in 
the regions where standard Poincar$\acute{\rm e}$ asymptotics breaks down. 

\begin{table}
\centering
\begin{tabular}{|c|c|c|} \hline
$N$ & Quantity & Value \\ \hline
& $F(1/10)$ &             1.73997257040229101538752631827936332290183806908929 \\ \hline 
& TS &                    0.83333333333333333333333333333333333333333333333333 \\
2 & $R^{SS}_1(1/10)$ &         -0.3205932520014175333654707340213015524648898035500\\
& Total &                 2.25271265173420681535538891759139510377028159887258 \\ \hline
& TS &                    -1.94444444444444444444444444444444444444444444444444 \\
3 & $R^{SS}_3(1/10)$ &          2.457184525776359388926618212330380430142089389324197 \\
& Total &                 2.252712651734205959869700086165299308599483016422822 \\ \hline
& TS &                    -5874.96031746031746031746031746031746031746031746031  \\
5 & $R^{SS}_5(1/10)$ &         5875.473057541649375261942492788406592113174013182747 \\
& Total &                 2.252712651734205959869701646368495118615533791559062 \\ \hline
& TS &                    -2.94867419474845489858725152842799901623431035195 $\times 10^{13}$ \\
9 & $R^{SS}_9(1/10)$ &         2.948674194748506172595384719922447233767119265136 $\times 10^{13}$ \\
& Total &                 2.25271265173420595986970164636849511861562722229495  \\ \hline
& TS &                    -3.60868558918311609670918346035346984255501011055 $\times 10^{31}$ \\
15 & $R^{SS}_{15}(1/10)$ &     3.60868558918311609670918346035352111656314330204 $\times 10^{31}$ \\
& Total &                 2.252712651734205959869701646368495118615627222294953  \\ \hline
& $\ln \Gamma(1/10)$ &    2.252712651734205959869701646368495118615626380692264 \\ \hline
\end{tabular}
\caption{Determination of $\ln\Gamma(1/10)$ via Eq.\ (\ref{sixtyeight}) for various values of the truncation parameter, $N$}
\label{tab2}
\end{table}

Table\ \ref{tab2} presents a much smaller sample of the results for $z \!=\! 1/10$ in the central asymptotic form in Eq.\ (\ref{fiftyeight})
with $R_N^{SS}(z)$ given by Eq.\ (\ref{sixtyeight}). As in the previous table the second row displays the value of Stirling's
approximation for $\ln \Gamma(1/10)$, whose value appears in the bottom row of the table. As expected, $F(1/10)$  represents a 
major contribution to $\ln \Gamma(1/10)$, but it is by no means close or accurate. If we were to consider the truncated series and 
$F(1/10)$ as an approximation to $\ln \Gamma(1/10)$ as in standard Poincar$\acute{\rm e}$ asymptotics, then the best possible 
approximation is the $N \!=\! 1$ case, but it too would not be regarded as accurate. The other values of the truncation parameter are 
even worse since the truncated series diverges far more rapidly than for $z \!=\!3$. Therefore, for this value of $z$, there is no 
optimal point of truncation, i.e. the remainder does not attain a minimum value before beginning to diverge. Because of this, 
the remainder diverges far more rapidly for low values of $N$ and thus, a greater cancellation of decimal places occurs than
in the previous table. Consequently, the total values for each value of the truncation parameter in Table\ \ref{tab2} are generally 
not as accurate as those in the previous table, the exception being the very low values of $N$ such as $N \!=\! 2$. As described 
previously, these values of $N$ are affected by the size of limit. Nevertheless, we would not have been able to achieve these results 
if the remainder had not undergone a correct regularization.

Now we assume that we do not have access to a powerful routine such as Gamma[N,z] in Mathematica to evaluate the incomplete gamma 
function. In this instance we need to create a program based on Eq.\ (\ref{sixty}), which means in turn that we require a numerical
integration routine. The second program in the appendix presents a module that employs the NIntegrate routine in Mathematica.  
Consequently, the remainder is being evaluated by another method. This is despite the fact that the only difference between the two 
programs occurs in the Do loop in which the NIntegrate routine appears. Therefore, the second program can be used as a check on the 
results generated by the first one. In the version appearing in the appendix, both the precision and accuracy goals have been set to
30, meaning that we are seeking thirty figure accuracy. To achieve such accuracy, the working precision must be set to a much higher 
value. Thus, WorkingPrecision has been set to 60, which actually yields more accurate values than those specified by the accuracy and 
precision goals. One can select higher values for all these options, but it comes at the expense of the computing time. The integrand 
used in the NIntegrate routine is called Intgrd. It appears outside the module in the first line of the program and is essentially the 
integrand in the integral given in Eq.\ (\ref{twentytwo}). Moreover, the calculated quantities in the program are all printed out to 
25 decimal places, which means that we should not expect greater accuracy than this, although in practice the results are often more 
accurate than specified. Because the leading term $F(z)$ or the Stirling approximation is the dominant term for $\ln \Gamma(z)$, we 
expect that for $|z| \!=\! 3$ and the truncation parameter lower than the optimal point of truncation, i.e. $N \!<\! N_{OP}$, the 
results will be accurate to at least 25 decimal places if the third asymptotic form in Eq.\ (\ref{fiftyeight}) is correct.
 
\begin{table}
\small
\centering
\begin{tabular}{|c|c|c|c|} \hline
$N$ & $\theta$ & Quantity  &  Value \\ \hline
& &$F(z)$ &           -0.4989796572635297888298191958 + 1.38469038925736775278152701797$\, i$ \\    
& & TS &              0.02250122010739348973103804218 - 0.01622952279862886741804878275$\, i$ \\
4 & $\pi/5$ & $R^{SS}_{4}(z)$ &   1.10075483159775932614119$ \times 10^{-7}$ - 2.3053919203213893994435144$\times 10^{-7}\,i$ \\
& & $SD^{SS}_0(z)$   &  0 \\ 
& & Total &           -0.4764783270806531393228485395 + 1.36846063591954685322453829086$\, i$ \\
& & $\ln \Gamma(z)$ & -0.4764783270806531393228485395 + 1.36846063591954685322453829086$\, i$ \\ \hline  
& &$F(z)$ &           -2.2031482244785443427155424131 + 1.30339979253231970320178992528$\, i$ \\    
& & TS &              0.01399322119549584041027516121 - 0.02405320763867771412867472768$\, i$ \\
12 & $\pi/3$ & $R^{SS}_{12}(z)$ &   8.0074867806568638247593$ \times 10^{-10}$ - 9.415241580707277926560568$\times 10^{-10}\,i $  \\
& & $SD^{SS}_0(z)$   &  0 \\ 
& & Total &           -2.1891550024822998242395808694 + 1.27934658395211783100238740494$\, i$ \\
& & $\ln \Gamma(z)$ & -2.1891550024822998242395808694 + 1.27934658395211783100238740494$\, i$ \\ \hline   
& &$F(z)$ &           -3.5024573297231710758524262838 + 0.51402564699132619931253447636$\, i$ \\  
& & TS &              0.00624853540876701851029121754 - 0.02716316476357736737433024841$\, i$ \\ 
6 & $3\pi/7$ & $R^{SS}_{6}(z)$ &  5.01448516214019677825215$ \times 10^{-9}$ + 1.423983690056893098624148$\times 10^{-8}\,i$ \\
& & $SD^{SS}_0(z)$   &  0 \\ 
& & Total &           -3.4962087892999188952019382880 + 0.48686249646758573250713521418$\, i$ \\
& & $\ln \Gamma(z)$ & -3.4962087892999188952019382880 + 0.48686249646758573250713521418$\, i$ \\ \hline 
& &$F(z)$ &           -5.2196841368342001927923955103 - 3.93258796345066902756928238018$\, i$ \\ 
& & TS &              -0.0139932183105014143643078934 - 0.02405320763867771412867472768$\, i$ \\ 
7 & $2\pi/3$ & $R^{SS}_{7}(z)$ &  -3.68574310411165365028398$\times 10^{-9}$ - 2.18129024012924399653960$\times 10^{-9}\,i$ \\ 
& & $SD^{SS,U}_{1}(z)$   &              -8.13752781094718217957452$\times 10^{-8}$ + 0$\,i$ \\ 
& & Total &           -5.2336774402057228207401788498 - 3.95664117203087089976868490052$\, i$ \\
& & $\ln \Gamma(z)$ & -5.2336774402057228207401788498 - 3.95664117203087089976868490052$\, i$ \\ \hline
& &$F(z)$ &           -5.0042973272271249734301225781 - 5.40748566131578852467236682820$\, i$ \\ 
& & TS &              -0.0182755216757311328100261387 - 0.02093456705773778718274769556$\, i$ \\ 
10 & $8\pi/11$ & $R^{SS}_{10}(z)$ &  -8.5614036344993217453630$\times 10^{-10}$ + 1.9101719530299621368772$\times 10^{-10}\,i$ \\
& & $SD^{SS,U}_{1}(z)$   &                 6.34447088671232496941997$\times 10^{-7}$ + 1.43565508702212979167964$\times 10^{-7}\,i $ \\
& & Total &           -5.0225722153119077984575839494 - 5.42842008461700041433913914211$\, i$ \\ 
& & $\ln \Gamma(z)$ & -5.0225722153119077984575839494 - 5.42842008461700041433913914211$\, i$ \\ \hline
& &$F(z)$ &           -1.3917914586609929543190298089 - 10.2463903684759447097075243173$\,i $ \\ 
& & TS &              -0.0273105469633224052584334852 - 0.00452525983710157150438833607$\, i$ \\ 
5 & $18\pi/19$ & $R^{SS}_{5}(z)$ & -4.90851551222138258253736$\times 10^{-9}$ -  3.45838244247219429147328$\times 10^{-8}\,i$ \\
& & $SD^{SS,U}_{1}(z)$   &   0.0443596569684061732573638356 + 0.01194376356035533612936953659$\, i$\\
& & Total &           -1.3747423535644246985414820411 - 10.2389718993365153698044860315$\, i$ \\
& & $\ln \Gamma(z)$ & -1.3747423535644246985414820411 - 10.2389718993365153698044860315$\, i$ \\ \hline 
\end{tabular}
\normalsize
\caption{Determination of $\ln\Gamma(z)$ via Eq.\ (\ref{fiftyeight}) with $|z| \!=\! 3$ and various values of the truncation
parameter, $N$ and argument, $\theta$}
\label{tab3}
\end{table}

Unlike the results in the previous tables we now consider complex values of $z$ by letting $\theta$ take values within the 
principal branch of the complex plane except for $\pm \pi/2$. Table\ \ref{tab3} presents a very small sample of the the results 
obtained by running the second program in the appendix with $|z| \!=\!3$. Although both positive and negative values of $\theta$
were considered, only positive values appear in the table, while negative values of $\theta$ will be presented when we discuss 
$|z| \!=\! 1/10$. Furthermore, the calculations generally took between 3 and 5 hours with the Sony VAIO laptop, while those utilising 
the Sunfire server took between 5.5 and 8.5 hours. Occasionally, calculations with the latter machine took longer to complete when 
$\theta$ was close to a Stokes line, which will be discussed later. 

From the table it can be seen that there is a set of six results for each couple of the truncation parameter and $\theta$. The first 
of these is Stirling's approximation or $F(z)$, where $z \!=\! 3 \exp(i \theta)$. Because the truncation parameter is either less 
than or close to $N_{OP}$, we see that $F(z)$ is quite close to the actual value of $\ln \Gamma(3 \exp(i \theta))$, which appears 
as the bottom value of each set.  However, it is by no means very accurate. The second value denoted TS represents the value of
the truncated sum given by Eq.\ (\ref{fiftysevena}). If the second value is added to the corresponding value of $F(z)$, then it 
yields a better or far more accurate value to $\ln \Gamma(z)$ than Stirling's approximation. This represents standard asymptotic
procedure, but it is, of course, dependent upon the remainder being very small and the Stokes discontinuity term being negligible.
The next value in each set of six results is the remainder calculated via Eq.\ (\ref{sixty}). From the table we see that the remainder
is at most of the order of $10^{-7}$, which means that it does not represent a significant contribution to $\ln \Gamma(3 \exp(i \theta))$. 
Of course, this is not the case when the truncation parameter exceeds $N_{OP}$ as we observed in the last couple of
calculations in Table\ \ref{tab1}.   

The fourth value in each set or calculation of $\ln \Gamma(z)$ is the Stokes discontinuity term. As indicated earlier, this term 
is zero when $|\theta| < \pi/2$. Therefore, for the first three sets of values in the table, the fourth value is given as zero. For 
$\theta > \pi/2$, the Stokes discontinuity term is given by the upper-signed version of Eq.\ (\ref{sixtyone}), which is denoted by 
$SD^{SS,U}_1(z)$ since $M=1$. This term is subdominant to the truncated series when $\theta$ is close to the Stokes line at $\pi/2$, 
but becomes more dominant as $\theta$ increases to $\pi$. It should be noted that midway in the Stokes sector or at $\theta \!=\!\pi$, 
which is known as an anti-Stokes line, the Stokes discontinuity term is equally as dominant as the sum of the truncated series and its 
remainder. We see this occurring with  $\theta = 18 \pi/19$ in the table. For this calculation the truncated sum and remainder 
is of the same order as $SD^{SS,U}_1 \left(3 \exp(18 \pi i/19) \right)$. In fact, in magnitude the latter is larger. Beyond 
$\theta = \pi$, the Stokes discontinuity is expected to dominate the truncated series and its remainder \cite{hea62}. All the 
behaviour mentioned here only occurs if $|z| \! \gg \! 1$.

The fifth value in each set or calculation of $\ln \Gamma(z)$, which is denoted by Total, represents the sum of the four preceding values, 
whilst the value immediately below, as stated before, is the actual value of the special function obtained by using Mathematica's 
intrinsic routine called LogGamma[z]. In all cases we see that both the real and imaginary parts of the totals agree exactly 
to well over 25 decimal places with those obtained via LogGamma[z], which is well within the accuracy and precision goals specified 
in the second program in the appendix. Once again, we have obtained exact values of $\ln \Gamma(z)$ from the asymptotic forms in 
Thm.\ 2.1. Note also that although the remainder is very small for each set of values in Table\ \ref{tab3} and the Stokes discontinuity 
term is very small in the fourth and fifth sets, they are still necessary in order to achieve the remarkable agreement between 
the Total values and those obtained via LogGamma[z].       

\begin{table}
\small
\centering
\begin{tabular}{|c|c|c|c|} \hline
$N$ & $\theta$ & Quantity & Value \\ \hline
& &$F(z)$ &                  1.75803888205251701300823152720 + 0.38158365834299627447460123156$\, i$\\    
& & TS &                     0.72168783648703220563643597562 - 2.36111111111111111111111111111$\, i$ \\
3 & $-\pi/6$ & $R^{SS}_{3}(z)$ &  -0.2230295240392980035338083054 + 2.52524252152237336263247340087$\, i$ \\
& & $SD^{SS}_0(z)$   &  0 \\ 
& & Total &                  2.25669719450025121511085919742 + 0.54571506875425852599596352132$\, i$ \\
& & $\ln \Gamma(z)$ &        2.25669719450025121511085784624 + 0.54571506875425852599596430142$\, i$ \\ \hline  
& &$F(z)$ &                  1.82912062061888235707377761234 + 0.82415547097542000625296951892$\, i$ \\   
& & TS &                     -1.92303004436996354502288$\times 10^{8}$ - 642086.002046734488604705084952$\, i$ \\
7 & $-4\pi/11$ & $R^{SS}_{7}(z)$ & 1.923030048814754111406220$\times 10^{8}$ + 642086.366490882117931701673660$\, i$ \\
& & $SD^{SS}_0(z)$   &  0 \\ 
& & Total &                  2.27359967725721630854389163795 + 1.18859961860474700284167763994$\, i$ \\
& & $\ln \Gamma(z)$ &        2.27359967725721630854389163795 + 1.18859961860474700284167763994$\, i $ \\ \hline
& &$F(z)$ &                  1.88648341970221940135996338478 + 1.03535606610194782214347998160$\, i$ \\ 
& & TS &                     2.87562548020794239198561$\times 10^{13}$ - 7.0718880105443602759020497$\times 10^{12}\,i$ \\ 
9 & $-6\pi/13$ & $R^{SS}_{9}(z)$ & -2.8756254802079022596675$\times 10^{13}$ + 7.0718880105448298576076792$\times 10^{12}\, i$ \\
& & $SD^{SS}_0(z)$   &  0 \\ 
& & Total &                  2.28780660084741914752819484319 + 1.50493777173150666351075080995$\, i$ \\
& & $\ln \Gamma(z)$ &        2.28780660084741914752819484319 + 1.50493777173150666351075080994$\, i $\\ \hline
& &$F(z)$ &                  1.93811875120925961146815019100 +  1.18372127170949939121742184779$\, i$ \\ 
& & TS &                     -6.2877629092633776151775$\times 10^{10}$ + 1.3406164598876901506339999$\times 10^{10}\,i$ \\
8 & $-8\pi/15$ & $R^{SS}_{8}(z)$ & 6.28776290922350273134009$\times 10^{10}$ - 1.3406164598401660891969575$\times 10^{10}\, i$ \\
& & $SD^{SS,L}_{1}(z)$   &              0.76110557640259383178972540915 + 0.07527936657383153773373307240$\, i$ \\ 
& & Total &                  2.30047548923720786540859926314 + 1.73424125265375514384575434070$\, i $ \\
& & $\ln \Gamma(z)$ &        2.30047548923720786540859926314 + 1.73424125265375514384575434070$\, i $ \\ \hline
& &$F(z)$ &                  2.13715100099092628642763759338 + 1.57823338674728063558799447261$\, i $ \\ 
& & TS &                      -590945.19319121626731651297966 + 599359.202691478098295355725315$\, i $ \\ 
6 & $-3\pi/4$ & $R^{SS}_{6}(z)$ & 590944.712330718986993690696285 - 599358.955111729858535640849047$\, i $ \\
& & $SD^{SS,L}_{1}(z)$   &              0.68679805984095965121150997224 + 0.579703018063676729767943942736$\, i $ \\  
& & Total &                  2.34308856355156311535576891244 + 2.405516153050717080232206634297$\, i $ \\ 
& & $\ln \Gamma(z)$ &        2.34308856355156311535576891244 + 2.405516153050717080232206634297$\, i$ \\ \hline
& &$F(z)$ &                  2.33668492162243351553206801970 + 1.825916904516483614559881067115$\, i$ \\  
& & TS &                     -42.600558891527544536579217000 + 64.60897639111337349406319763878$\, i $ \\ 
4 & $-15\pi/16$ & $R^{SS}_{4}(z)$ & 42.0897905773173511704765793260 - 64.54765565501832133510270286263$\, i $ \\
& & $SD^{SS,L}_{1}(z)$   &              0.54123306366541416118208181725 + 1.072657474660830843519039814447$\, i $ \\
& & Total &                  2.36714967107765431061151216215 + 2.959895115272366617039415657712$\, i $ \\ 
& & $\ln \Gamma(z)$ &        2.36714967107765431061151216215 + 2.959895115272366617039415657712$\, i$ \\ \hline 
\end{tabular}
\normalsize
\caption{Determination of $\ln\Gamma(z)$ via Eq.\ (\ref{fiftyeight}) with $|z| \!=\! 1/10$ and various values of
the truncation parameter, $N$ and argument, $\theta$}
\label{tab4}
\end{table}

Table\ \ref{tab4} presents yet another small sample of the results obtained by running the second program in the appendix, 
but on this occasion, $|z| \!=\! 1/10$. Whilst positive values of $\theta$ were also studied, only negative ones appear 
in the table. Therefore, when the Stokes discontinuity term is non-zero, $SD^{SS,L}_{1}(\exp(i\theta)/10)$  in the fifth
asymptotic form of Eq.\ (\ref{fiftyeight}) appears in the table. The calculations had similar CPU times to those in the
previous table or the $|z| \!=\! 3$ case. For each calculation of $\ln \Gamma(\exp(i \theta)/10)$, there are again six values. 
However, since there is no optimal point of truncation, the results in Table\ \ref{tab4} are radically different from those 
in the previous table. Now the values are dominated primarily by the truncated sum and its regularized remainder. In fact, 
for $N \!>\! 3$, the values of the truncated sum and its remainder eclipse all the other quantities, which results in a problem 
for the final values of $\ln \Gamma(\exp(i \theta)/10)$ or those denoted by Total. Because the remainder and truncated sum 
dominate in different directions, a cancellation of many decimal places occurs. This puts pressure on the accuracy of the 
total values. E.g., for $N \!=\! 9$ and $\theta \!=\! -6 \pi/13$, both the truncated sum and the regularized remainder are 
of the order of $10^{13}$, which results in the loss of thirteen decimal places when they are summed together. Yet, the accuracy 
and precision goals have been set to 30 in the second program. Hence the sum of the truncated series and the regularized remainder 
should only be accurate to only 17 decimal places, which limits the accuracy of the total values. Fortunately, it can be seen that 
the total values agree with the values of $\ln \Gamma(z)$ immediately below them to 28 decimal places. This means that despite the 
fact that the precision and accuracy goals were set to 30, the total values were ultimately far more accurate than expected because 
the working precision had been set to a much higher value, viz.\ 60.  

From the table we see that although Stirling's approximation or $F(\exp(i\theta)/10)$ provides a substantial contribution to 
$\ln \Gamma(\exp(i \theta)/10)$, it is now inaccurate. The truncated sum is capable of improving the accuracy slightly 
for small values of the truncation parameter. For example, we see that when the truncated sum is added to $F(z)$ for $N \!=\! 3$ 
and $\theta \!=\! -\pi/6$, the real part is closer to the real part of $\ln \Gamma(\exp(-i\pi/6)/10)$, but the imaginary part is 
even more inaccurate. In actual fact, all the results in the table are dominated by the truncated sum and its regularized remainder, 
but because these quantities act against each other, their sum is not as significant as Stirling's approximation. That is,
one can no longer neglect the remainder, which means that standard Poincar$\acute{\rm e}$ asymptotics has broken down or is
useless. In order to obtain the exact value of $\ln \Gamma(\exp(-i\pi/6)/10)$ via Eq.\ (\ref{fiftyeight}), we require the remainder 
so that it counterbalances the truncated sum. This counterbalance will only occur if a correct or proper regularization of the later 
terms in the series given by Eq.\ (\ref{sixteen}) has been performed. Therefore, in the table we see that when the regularized 
value of the remainder is added to Stirling's approximation and the truncated sum for $|\theta | \! < \! \pi/2$, we obtain the exact 
values of $\ln \Gamma(\exp(i \theta)/10)$. For $\theta \!<\! -\pi/2$, we also need to sum the Stokes discontinuity term given 
by the lower-signed version of Eq.\ (\ref{sixtyone}). Moreover, it can be seen that $SD^{SS,L}_{1}(z)$ is greater than the sum of
the truncated series and the regularized remainder, which demonstrates the importance of the Stokes discontinuity term 
outside the primary Stokes sector. That is, it represents a crucial contribution for obtaining the values of $\ln \Gamma(z)$ 
when $|\theta| \! >\! \pi/2$.

So far, we have managed to verify the asymptotic forms in Eq.\ (\ref{fiftyeight}) pertaining to the Stokes sectors. Now we turn to the 
asymptotic forms for the two Stokes lines situated within the principal branch or the second and fourth asymptotic forms in Eq.\ 
(\ref{fiftyeight}). Since $\theta$ is fixed in both these asymptotic forms, we see immediately that the Stokes discontinuity term 
will only depend upon the magnitude of $z$, which means in turn that the Stokes discontinuity term along the Stokes lines is real. 
Furthermore, since $TS_N(z)$ depends only on odd powers of $z$ according to Eq.\ (\ref{fiftyseven}), $TS_N(z)$ and consequently, 
$R^{SL,}_N(z)$ must be imaginary along both Stokes lines. This behaviour is consistent with Rule D given in Ch.\ 1 of Ref.\cite{din73}, 
which states that on crossing a Stokes line, an asymptotic series generates a discontinuity in form that is $\pi/2$ out of phase with 
the series on the Stokes line.  

\begin{table}
\small
\centering
\begin{tabular}{|c|c|c|} \hline
$N$ & Quantity & Value  \\ \hline
& $F(3\exp(i\pi/2))$ &                  -4.3427565915140719616112579569 - 0.4895612973931192354299251350522$\, i$ \\
&  $SD^{SL}_0(3\exp(i\pi/2))$   &             3.256206078642828367679816468$\times 10^{-9}$ $\;\;\;\;\;\;\;\;\;\;\;\;\;\;\;\;\;\;\;\;\;\;\;\;\;\;\;
\;\;\;\;\;\;\;\;\;\;\;\;\;\;\;\;\;\;\;\;\;\;\;\;\;$ \\  
&  Combined &                -4.3427565882578658829684295892 - 0.4895612973931192354299251350522$\, i$ \\ \hline
&  TS &  $\;\;\;\; 0 \;\;\;\;\;\;\;\;\;\;$ 
\\
1 & $R^{SL}_{1}(3\exp(i \pi/2))$  &  $\;\;\;\;\;\;\;\;\;\;\;\;\;\;\;\;\;\;\;\;\;\;\;\;\;\;\;\;\;\;\;\;\;\;\;\;\;\;\;\;\;\;\;\;\;\;\;\;\;\;$ 
- 0.0278840894653691199321777792256$\, i$ \\
&  Total &                  -4.3427565882578658829684295892 - 0.5174453868584883553621029142779$\, i$ \\ \hline
&  TS &  $\;\;\;\;\;\;\;\;\;\;\;\;\;\;\;\;\;\;\;\;\;\;\;\;\;\;\;\;\;\;\;\;\;\;\;\;\;\;\;\;\;\;\;\;\;\;\;\;\;\;\;\;\;\;\;0$ 
- 0.0278842394252900781377131527007$\, i$ \\ 
6 & $R^{SL}_{6}(3\exp(i \pi/2))$  &  $\;\;\;\;\;\;\;\;\;\;\;\;\;\;\;\;\;\;\;\;\;\;\;\;\;\;\;\;\;\;\;\;\;\;\;\;\;\;\;\;\;\;\;\;\;\;\;\;\;\;\;\;\;\;\;0$
- 1.8907874105339892863379255$\times 10^{-8}\, i$ \\
&  Total &                 -4.3427565882578658829684295892 - 0.51744555572628341890753115113225$\, i$ \\ \hline
&  TS &  $\;\;\;\;\;\;\;\;\;\;\;\;\;\;\;\;\;\;\;\;\;\;\;\;\;\;\;\;\;\;\;\;\;\;\;\;\;\;\;\;\;\;\;\;\;\;\;\;\;\;\;\;\;\;\;0$ 
- 0.0278842563298976281594154202028$\, i$ \\
9 & $R^{SL}_{9}(3\exp(i\pi/2))$  &  $\;\;\;\;\;\;\;\;\;\;\;\;\;\;\;\;\;\;\;\;\;\;\;\;\;\;\;\;\;\;\;\;\;\;\;\;\;\;\;\;\;\;\;\;\;\;\;\;\;\;\;\;\;\;\;0$
+ 3.2562060786428283676798164$\times 10^{-9}\, i$ \\
&  Total &                 -4.3427565882578658829684295892 - 0.51744555572628341890753115113225$\, i$ \\ \hline
&  TS &  $\;\;\;\;\;\;\;\;\;\;\;\;\;\;\;\;\;\;\;\;\;\;\;\;\;\;\;\;\;\;\;\;\;\;\;\;\;\;\;\;\;\;\;\;\;\;\;\;\;\;\;\;\;\;\;0$ 
- 0.0278842691899612112195938305035$\, i$ \\
15 & $R^{SL}_{15}(3\exp(i \pi/2))$  &  $\;\;\;\;\;\;\;\;\;\;\;\;\;\;\;\;\;\;\;\;\;\;\;\;\;\;\;\;\;\;\;\;\;\;\;\;\;\;\;\;\;\;\;\;\;\;\;\;\;\;\;\;\;\;\;0$
+ 1.0856797027741987814423624$\times 10^{-8}\, i $ \\
&  Total &                 -4.3427565882578658829684295892 - 0.51744555572628341890753115113225$\, i$ \\ \hline
&  TS &  $\;\;\;\;\;\;\;\;\;\;\;\;\;\;\;\;\;\;\;\;\;\;\;\;\;\;\;\;\;\;\;\;\;\;\;\;\;\;\;\;\;\;\;\;\;\;\;\;\;\;\;\;\;\;\;0$ 
- 0.0278853616586139260908931195458$\, i$ \\
20 & $R^{SL}_{20}(3\exp(i \pi/2))$  &  $\;\;\;\;\;\;\;\;\;\;\;\;\;\;\;\;\;\;\;\;\;\;\;\;\;\;\;\;\;\;\;\;\;\;\;\;\;\;\;\;\;\;\;\;\;\;\;\;\;\;\;\;\;\;\;0$
+ 1.1033254497426132871034659$\times 10^{-6}\,i $ \\
&  Total &                 -4.3427565882578658829684295892 - 0.51744555572628341890753115113225$\, i$ \\ \hline
&  TS &  $\;\;\;\;\;\;\;\;\;\;\;\;\;\;\;\;\;\;\;\;\;\;\;\;\;\;\;\;\;\;\;\;\;\;\;\;\;\;\;\;\;\;\;\;\;\;\;\;\;\;\;\;\;\;\;0$ 
- 52.072356609356813219352046137393$\, i $ \\
30 & $R^{SL}_{30}(3\exp(i \pi/2))$  &  $\;\;\;\;\;\;\;\;\;\;\;\;\;\;\;\;\;\;\;\;\;\;\;\;\;\;\;\;\;\;\;\;\;\;\;\;\;\;\;\;\;\;\;\;\;\;\;\;\;\;\;\;\;\;\;0$
+ 52.044472351023649035874440121314$\, i$ \\
&  Total &                 -4.3427565882578658829684295892 - 0.51744555572628341890753115113225$\, i$ \\ \hline
&  TS &  $\;\;\;\;\;\;\;\;\;\;\;\;\;\;\;\;\;\;\;\;\;\;\;\;\;\;\;\;\;\;\;\;\;\;\;\;\;\;\;\;\;\;\;\;\;\;\;\;\;\;\;\;\;\;\;0$ 
- 6.4908409843349435181620453$\times 10^{25}\, i$ \\
50 & $R^{SL}_{50}(3\exp(i \pi/2))$  &  $\;\;\;\;\;\;\;\;\;\;\;\;\;\;\;\;\;\;\;\;\;\;\;\;\;\;\;\;\;\;\;\;\;\;\;\;\;\;\;\;\;\;\;\;\;\;\;\;\;\;\;\;\;\;\;0$
+ 6.4908409843349435181620453$\times 10^{25}\, i$ \\
&  Total &                 -4.3427565882578658829684295892 - 0.51744555572628341890753115113225$\, i$ \\ \hline
&  $\ln \Gamma(3\exp(i \pi/2))$ &      -4.3427565882578658829684295892 - 0.51744555572268341890753115113225$\, i$ \\ \hline
\end{tabular}
\normalsize
\caption{Determination of $\ln\Gamma \! \left(3\exp(i \pi/2)\right)$ via Eq.\ (\ref{fiftyeight}) for various values of the truncation
parameter, $N$}
\label{tab5}
\end{table}

The third program in the appendix presents an implementation of the second and fourth asymptotic forms of Eq.\ (\ref{fiftyeight}) in
Mathematica. When compared with the previous programs, one can see that the Do loop is now very different because of the inclusion 
of a Which statement. This statement is necessary because the singularity in the Cauchy principal value integral in the remainder 
$R^{SL}_N(z)$ given by Eq.\ (\ref{sixtya}) alters with each value of k in the Do loop. Since the integral has been divided into 
smaller intervals, one needs to determine the interval in which the singularity is situated and then divide that interval into 
two intervals with the singularity acting as the upper and lower limits of the resulting integrals. As the range of the Cauchy principal 
value integral has been split into seven intervals, there are seven possibilities or conditions where an interval can be divided. Hence
the Which statement has been designed to consider all these possibilities. In addition, when calling NIntegrate, one is required to 
introduce the option Method--$>$PrincipalValue to avoid convergence problems. Moreover, WorkingPrecision has been extended to a 
value of 80. So, whilst Mathematica does indicate the level of accuracy for each result printed out by the program, which is for the 
most part 25 decimal places, in reality the results are expected to be more accurate. Unfortunately, extending the working precision 
also means there is a significant increase in the amount of processing. Consequently, the third program in the appendix takes much longer 
to execute than any of the preceding programs. To allow for this, the program was converted into a script program and run mainly on 
the SunFire server, thereby allowing many different values of the truncation parameter to be executed simultaneously.  

Table\ \ref{tab5} presents a sample of the results generated by running the third program in the appendix with the variable modz
set equal to 3. Although both Stokes lines were considered by putting the variable theta in the program equal to $\pm {\rm Pi}/2$, 
only the results for the positive value of theta have been presented to save space. Results generated for the other Stokes 
line will be presented when we consider modz equal to 1/10 shortly. The calculations took much longer for larger values of the 
truncation parameter, ranging from 26 hrs for $N \!=\!1$ to 47.5 hrs for $N \!= \! 50$. Because the values of $F(3\exp(i\pi/2))$ and 
$SD^{SL}_0(3\exp(i \pi/2))$ are independent of the truncation parameter, they only appear once at the top of the table, while their sum 
appears immediately below them in the row labelled Combined. As mentioned previously, the Stokes discontinuity term is real,
whereas the truncated sum and regularized value of the remainder are  imaginary. Therefore, the real part of the value in the 
Combined row represents the real part of $\ln \Gamma(3\exp(i \pi/2))$, which can be checked by comparing it with the real part of the value 
appearing at the bottom of the table. This means that the Stokes discontinuity term is responsible for correcting the real part of 
Stirling's approximation on a Stokes line. On the other hand, the imaginary part of $\ln \Gamma(3\exp(i\pi/2))$ can only be calculated 
exactly by a correct or proper regularization of the asymptotic series in Eq.\ (\ref{sixteen}), which led to Eq.\ (\ref{sixtya}).   
With regard to the imaginary part of $\ln \Gamma(3\exp(i \pi/2))$ in the table, the last decimal figure was printed out as a 6 instead
of a 5. This is because the accuracy was set to 25 decimal places due to the statement N[e5,25]. Since more than 25 figures appear
in the table, this statement should have been modified to consider a higher level of accuracy. Thus, we should only be concerned when
there is no agreement for less than 25 decimal places. The redundant places have been introduced to indicate that the results in the
Total column have been computed via a different approach from the intrinsic LogGamma function in Mathematica given at the bottom of the
table. Therefore, we should expect variation to occur at some stage, but outside the specified level of accuracy. 

As expected, from the table it can be seen that the regularized value of the remainder decreases steadily until the truncation
parameter hits the optimal point of truncation, viz.\ $N_{OP}=11$, before it begins to increase. Note, however, that the imaginary
part of the Total value for $N \!=\! 1$ is only accurate to 6 decimal places with the imaginary part of $\ln \Gamma(3 \exp(i \pi/2))$.
As discussed previously, this arises because the power of $n$ in the denominator of $R^{SL}_1(z)$ is zero for $N$ or TP in the thrid
program in the appendix is equal to unity. Consequently, we do not find that as k becomes large in the Do loop, e2 decreases significantly 
when it is divided by k$^{\wedge}$(2TP-2). From the vast number of results obtained in running the program, it is found that the 
optimal point of truncation occurs when $N \!=\! 11$, while $R_{11}(3 \exp(i \pi/2))$ has a magnitude of the order of $10^{-11}$ even 
though the calculation is not displayed in Table\ \ref{tab5}. Beyond the optimal point of truncation or for $N > 11$, the magnitude of the 
regularized value of the remainder increases steadily so that its magnitude is of the order of $10^{-6}$ when $N \!=\!20$. By the time 
we reach $N \!=\! 30$, both the truncated series and regularized value of the remainder dominate the calculation of $\ln \Gamma(3 \exp(i\pi/2))$, 
but as we found previously, they act against or cancel each other, thereby yielding the extra small imaginary part required to make the 
imaginary part in the Combined row agree with that for $\ln \Gamma(3 \exp(i\pi/2))$. In fact, the most surprising result in the table 
is the last run or $N \!=\! 50$ result because at least 25 decimal places need to cancel before we end up with the regularized value 
for the entire asymptotic series. As mentioned previously, the cancellation of these decimal places puts pressure on the accuracy and precision 
goals, which have been set to 30 in the third program in the appendix. Fortunately, because WorkingPrecision was set to 80, it appears that 
the neglected terms in setting the Do loop to a limit of $10^5$ have negligible effect. Consequently, the remainder has been evaluated to a 
much greater accuracy than specified by the accuracy and precision goals in the program. Hence the Total value for $N \!=\! 50$  is as 
accurate as the actual value of $\ln \Gamma(3 \exp(i \pi/2))$.   

\begin{table}
\small
\centering
\begin{tabular}{|c|c|c|} \hline
$N$ & Quantity &  Value \\ \hline
& $F(\exp(-i\pi/2)/10)$ &                  1.91315144702220592186619329458 + 1.11565667269685287801745999128$\, i$\\ 
&  $SD^{SL}_0(\exp(-i\pi/2)/10)$   &            0.381235865406433806218304673501 + 0$\,i\;\;\;\;\;\;\;\;\;\;\;\;\;\;\;\;\;\;\;\;\;\;\;\;\;\;\;
\;\;\;\;\;\;\;\;\;\;\;\;\;\;\;\;\;\;\;\;\;\;\;\;\;$ \\  
&  Combined &                  2.29438731242863972808449796808 + 1.11565667269685287801745999128$\, i$ \\ \hline
&  TS &  $\;\;\;\;\;\;\;\;\;\;\;\;\;\;\;\;\;\;\;\;\;\;\;\;\;\;\;\;\;\;\;\;\;\;\;\;\;\;\;\;\;\;\;\;\;\;\;\;\;\;\;\;\;\;\;0$ 
+ 3.61111111111111111111111111111$\, i$\\
3 & $R^{SL}_{3}(\exp(-i\pi/2)/10)$  &  $\;\;\;\;\;\;\;\;\;\;\;\;\;\;\;\;\;\;\;\;\;\;\;\;\;\;\;\;\;\;\;\;\;\;\;\;\;\;\;\;\;\;\;\;\;\;\;\;\;\;\;\;\;\;\;0$
- 3.09864851659634765254576003084$\, i$ \\ 
&  Total &                  2.29438731242863972808449796808 + 1.62811926721161633658281107155$\, i$ \\ \hline
&  TS &  $\;\;\;\;\;\;\;\;\;\;\;\;\;\;\;\;\;\;\;\;\;\;\;\;\;\;\;\;\;\;\;\;\;\;\;\;\;\;\;\;\;\;\;\;\;\;\;\;\;\;\;\;\;\;\;0$ 
+ 6035.35714285714285714285714285$\, i$ \\  
5 & $R^{SL}_{5}(\exp(-i\pi/2)/10)$  &  $\;\;\;\;\;\;\;\;\;\;\;\;\;\;\;\;\;\;\;\;\;\;\;\;\;\;\;\;\;\;\;\;\;\;\;\;\;\;\;\;\;\;\;\;\;\;\;\;\;\;\;\;\;\;\;0$
- 6034.844680262628093684291790216$\, i$ \\
&  Total &                  2.29438731242863972808449796808 + 1.62811926721161633658281263176$\, i$ \\ \hline
&  TS &  $\;\;\;\;\;\;\;\;\;\;\;\;\;\;\;\;\;\;\;\;\;\;\;\;\;\;\;\;\;\;\;\;\;\;\;\;\;\;\;\;\;\;\;\;\;\;\;\;\;\;\;\;\;\;\;0$ 
+ 1.92600477951585451585451$\times 10^8 \,i$ \\
7 & $R^{SL}_{7}(\exp(-i\pi/2)/10)$  &  $\;\;\;\;\;\;\;\;\;\;\;\;\;\;\;\;\;\;\;\;\;\;\;\;\;\;\;\;\;\;\;\;\;\;\;\;\;\;\;\;\;\;\;\;\;\;\;\;\;\;\;\;\;\;\;0$
- 1.926004774391228570706881$\times 10^8 \,i$ \\
&  Total &                  2.29438731242863972808449796808 + 1.62811926721161633658281263176$\, i$ \\ \hline
&  TS &  $\;\;\;\;\;\;\;\;\;\;\;\;\;\;\;\;\;\;\;\;\;\;\;\;\;\;\;\;\;\;\;\;\;\;\;\;\;\;\;\;\;\;\;\;\;\;\;\;\;\;\;\;\;\;\;0$ 
+ 1.7994052185642409074011$\times 10^{16}\, i$ \\
10 & $R^{SL}_{10}(\exp(-i\pi/2)/10)$  &  $\;\;\;\;\;\;\;\;\;\;\;\;\;\;\;\;\;\;\;\;\;\;\;\;\;\;\;\;\;\;\;\;\;\;\;\;\;\;\;\;\;\;\;\;\;\;\;\;\;\;\;\;\;\;\;0$
- 1.79940521856424085615493$\times 10^{16}\, i$\\
&  Total &                  2.29438731242863972808449796808 + 1.62811926721161633658281263176$\, i $ \\ \hline
&  TS &  $\;\;\;\;\;\;\;\;\;\;\;\;\;\;\;\;\;\;\;\;\;\;\;\;\;\;\;\;\;\;\;\;\;\;\;\;\;\;\;\;\;\;\;\;\;\;\;\;\;\;\;\;\;\;\;0$ 
+ 1.5698245268960591367284$\times 10^{25}\, i$ \\
13 & $R^{SL}_{13}(\exp(-i\pi/2)/10)$  &  $\;\;\;\;\;\;\;\;\;\;\;\;\;\;\;\;\;\;\;\;\;\;\;\;\;\;\;\;\;\;\;\;\;\;\;\;\;\;\;\;\;\;\;\;\;\;\;\;\;\;\;\;\;\;\;0$
- 1.56982452689605913672845$\times 10^{25}\, i$ \\
&  Total &                  2.29438731242863972808449796808 + 1.62811926721161633658281263176$\, i $ \\ \hline
&  TS &  $\;\;\;\;\;\;\;\;\;\;\;\;\;\;\;\;\;\;\;\;\;\;\;\;\;\;\;\;\;\;\;\;\;\;\;\;\;\;\;\;\;\;\;\;\;\;\;\;\;\;\;\;\;\;\;0$ 
+ 1.2373076433716581258003$\times 10^{52}\, i$ \\
21 & $R^{SL}_{21}(\exp(-i\pi/2)/10)$  &  $\;\;\;\;\;\;\;\;\;\;\;\;\;\;\;\;\;\;\;\;\;\;\;\;\;\;\;\;\;\;\;\;\;\;\;\;\;\;\;\;\;\;\;\;\;\;\;\;\;\;\;\;\;\;\;0$
- 1.2373076433716581258003$\times 10^{52}\, i$ \\                                              
&  Total &                  2.29438731242863972808449796808 + 1.62811926721161633658281263176$\, i $ \\ \hline
&  TS &  $\;\;\;\;\;\;\;\;\;\;\;\;\;\;\;\;\;\;\;\;\;\;\;\;\;\;\;\;\;\;\;\;\;\;\;\;\;\;\;\;\;\;\;\;\;\;\;\;\;\;\;\;\;\;\;0$ 
+ 5.3585692058787768499098$\times 10^{66}\, i$ \\
25 & $R^{SL}_{25}(\exp(-i\pi/2)/10)$  &  $\;\;\;\;\;\;\;\;\;\;\;\;\;\;\;\;\;\;\;\;\;\;\;\;\;\;\;\;\;\;\;\;\;\;\;\;\;\;\;\;\;\;\;\;\;\;\;\;\;\;\;\;\;\;\;0$
- 5.3585692058787768499098$\times 10^{66}\, i$ \\
&  Total &                  2.29438731242863972808449796808 - 6.729048843994533$\, i \;\;\;\;\;\;\;\;\;\;\;\;\;\;\;\;\;\;\;\;\;\;\;\;\;\;\;$ \\ \hline
&  $\ln \Gamma(\exp(-i\pi/2)/10)$ &        2.29438731242863972808449796808 + 1.628119267211616336582812631761$\, i $ \\ \hline
\end{tabular}
\normalsize
\caption{Determination of $\ln\Gamma\!\left(\exp(-i\pi/2)/10\right)$ via Eq.\ (\ref{fiftyeight}) for various values of the truncation 
parameter, $N$}
\label{tab6}
\end{table}

Table\ \ref{tab6} presents another small sample of the results obtained from running the third program in the appendix with the variables 
modz and theta set equal to 1/10 and $\pm {\rm Pi}/2$, respectively. Only the results for the alternative Stokes line at 
$\theta = -\pi/2$ are displayed. Again, the values obtained from the Stirling approximation and the Stokes discontinuity term appear 
at the top of the table. Immediately below these values is their combined sum, which represents the real part 
of $\ln \Gamma(\exp(-i\pi/2)/10)$. As in the previous table these results are only valid to about 25 decimal places because of the
accuracy and precision goals set in the numerical integration of the remainder. Nevertheless, because the working precision was set equal 
to 80, the accuracy of the results can be extended beyond 25 decimal places. Unlike the previous table, however, the Stokes discontinuity 
term is no longer negligible compared with the value of the Stirling approximation, and thus it cannot be neglected to yield a reasonable 
approximation for the real part of $\ln \Gamma(\exp(-\pi/2)/10)$. So whilst the asymptotic series is at peak exponential dominance to the 
Stokes discontinuity term on the Stokes line, which corresponds to Rule B in Ch.\ 1 of Ref.\ \cite{din73}, it does not mean that the 
latter is always negligible to the asymptotic series after regularization. We also see that both the truncated series and the regularized
value of the remainder diverge far more rapidly than in the previous table, which is due to the fact that there is no optimal point 
of truncation. In fact, the divergence is so rapid that by the time $N = 13$, we require 25 decimal places to be cancelled when
combining the regularized value of the remainder with the value of the truncated series. As a consequence, the most interesting 
result in the table is the final result for $N \!=\! 25$. In this instance both the truncated series and the regularized value of the 
remainder are of the order of $10^{66}$. Hence at least 66 decimal places need to be cancelled and as the remainder had the 
WorkingPrecision set equal to 80, it means that it will only possess about 14 decimal figure accuracy after cancellation with the 
truncated series. Furthermore, because the calculation proceeds well outside the accuracy and precision goals, the resulting value 
for the sum of both quantities will be incorrect, which turns out to be the case. Note also that since the real part of the Total 
value is not affected by the calculation of the remainder, it still yields the correct value for the real part of 
$\ln \Gamma(\exp(-i\pi/2)/10)$. We have thus seen the danger in selecting a value greater than 10 for the truncation parameter 
when there is no optimal point of truncation. If we say that the optimal point of truncation is zero in such cases as Table\ 
\ref{tab6}, then we see that exactification is extremely difficult to achieve when $N \gg 0$. To have any hope of success, 
we would have to increase the working precision and the accuracy and precision goals to unreasonable levels. 

So far, we have not seen any evidence whatsoever of the purported smoothing of the Stokes phenomenon as postulated by Berry and Olver 
in Refs.\ \cite{ber89} and \cite{olv90} and supported by Paris, Kaminski and Wood in Refs.\ \cite{par01}, \cite{par11} and \cite{par95}.
As indicated earlier, smoothing implies that there is no discontinuity in the vicinity of a Stokes line, whereas we have been
able to obtain exact values of $\ln \Gamma(z)$ near Stokes lines assuming the existence of a discontinuity. Because such smoothing occurs 
very rapidly in the vicinity of Stokes lines, it could, perhaps, be argued that the preceding analysis has not investigated the asymptotic 
behaviour of $\ln \Gamma(z)$ sufficiently close to the Stokes lines. If such a rapid transition does occur, then it means that we 
have still not exactified the Stokes approximation in the vicinity of the Stokes lines. From Fig.\ \ref{figone}, which represents the situation 
for $|z| \! =\! 3$, smoothing of the Stokes phenomenon is expected to be most pronounced for $\theta$ lying between $13 \pi/32$ and 
$19 \pi/32$. For these values of $\theta$ the Stokes multiplier should experience a rapid transition from 0 to 1. That is, 
the Stokes multiplier is expected to be quite close to a value of a 1/2 for small values of $\delta$, where $\theta= \pi (1/2+\delta)$ and 
$|\delta|< 3/32$. On the other hand, if the conventional view of the Stokes phenomenon still holds, then the Stokes multiplier $S^{+}$ will 
equal unity for $0<\delta <1$ and zero for $-1< \delta <0$. In other words, according to the conventional view, the Stokes multiplier 
remains a step function as indicated by Eq.\ (\ref{sixtyfour}). Hence there is a small region of positive and negative values of $\delta$, 
where one of the views can be disproved. In particular, by introducing very small values of $\delta$ such that $\theta$ always lies between 
$13 \pi/32$ and $19 \pi/32$ into the respective asymptotic forms in Eq.\ (\ref{sixtyeight}), we should not obtain exact values of 
$\ln \Gamma(z)$ if smoothing occurs because the Stokes multiplier should be close to 1/2 instead of toggling between zero and unity 
according to the sign of $\delta$. Therefore, to complete this section, we shall examine the results obtained from running the second 
program in the appendix by introducing small perturbations on either side of the Stokes line at $\theta \!=\! \pi/2$ with $|z| \!=\! 3$.    

\begin{table}
\small
\centering
\begin{tabular}{|c|c|c|} \hline
$\delta$ & Method & Value \\ \hline
1/10 &   LogGamma[z]  &         -5.1085546405054331385771175 - 2.43504864133618239587613036$\, i$ \\
& $SD^{SS,U}_1(z)$ &             0.0000000146924137960847328 + 0.00000000724920978735477097$\, i$  \\ 
& 1st AF &                      -5.1085546405054331385771175 - 2.43504864133618239587613036$\, i$ \\ \hline
-1/10 &  LogGamma[z]  &         -3.1156770612855851062960250 + 0.79152717486178700663566144$\, i$ \\
& 3rd AF &                      -3.1156770612855851062960250 + 0.79152717486178700663566144$\, i$ \\ \hline
1/100 &   LogGamma[z]  &        -4.4448078360199294879676721 - 0.68426539470619315579497619$\, i$ \\ 
& $SD^{SS,U}_1(z)$ &             0.0000000054543808883397577 - 0.00000000366845661861183983$\, i$ \\
& 1st AF &                      -4.4448078360199294879676721 - 0.68426539470619315579497619$\, i$ \\ \hline
-1/100 &   LogGamma[z]  &       -4.2360547825638102221663061 - 0.35681003461125834209091866$\, i$ \\
& 3rd AF &                      -4.2360547825638102221663061 - 0.35681003461125834209091866$\, i$ \\ \hline
1/1000 &   LogGamma[z]  &       -4.3531757575591613140088085 - 0.53385166100905755261595669$\, i$ \\ 
& $SD^{SS,U}_1(z)$ &             0.0000000065016016472424544 - 0.00000000038545945628149871$\, i$ \\
& 1st AF &                      -4.3531757575591613140088085 - 0.53385166100905755261595669$\, i$ \\ \hline
-1/1000 &   LogGamma[z]  &      -4.3322909095906129602545969 - 0.50110130347126170951651903$\, i$ \\
& 3rd AF &                      -4.3322909095906129602545969 - 0.50110130347126170951651903$\, i$ \\ \hline
1/10000 &   LogGamma[z]  &      -4.3438006028809735966127763 - 0.51908338527968766540121412$\, i$ \\ 
& $SD^{SS,U}_1(z)$ &             0.0000000065123040290213875 - 0.00000000003856476898298508$\, i$ \\
& 1st AF &                      -4.3438006028809735966127763 - 0.51908338527968766540121412$\, i$ \\ \hline
-1/10000 &   LogGamma[z]  &     -4.3417121085407199183370966 - 0.51580834470414165478538635$\, i$ \\
& 3rd AF &                      -4.3417121085407199183370966 - 0.51580834470414165478538635$\, i$ \\ \hline
1/20000 &   LogGamma[z]  &      -4.3438006028809735966127763 - 0.51908338527968766540121412$\, i$ \\ 
& $SD^{SS,U}_1(z)$ &             0.0000000065123851251757157 - 0.00000000001928245580002624$\, i$ \\
& 1st AF &                      -4.3438006028809735966127763 - 0.51908338527968766540121412$\, i$ \\ \hline
-1/20000 &   LogGamma[z]  &     -4.3422344065179726897501879 - 0.51662687288967352139359494$\, i$ \\
& 3rd AF &                      -4.3422344065179726897501879 - 0.51662687288967352139359494$\, i$ \\ \hline
\end{tabular}
\normalsize
\caption{Evaluation of $\ln\Gamma \! \left(3\exp(i(1/2+\delta)\pi ) \right)$ via Eq.\ (\ref{fiftyeight}) for 
various values of $\delta$}
\label{tab7}
\end{table}

Table\ \ref{tab7} presents another small sample of the results obtained by running the second program in the appendix for 
$|z| \!=\! 3$ and various values of $\delta$, where $\theta \!=\! (1/2 \!+\! \delta)\pi$. Besides using different values values 
of $\delta$, the code was also run for different values of the truncation parameter except for small values so as to eliminate 
the problem of not being to calculate the remainder accurately when the power of $n$ in the denominator of Eq.\ (\ref{sixty}) 
is too small. For each positive value of $\delta$ there are three rows of values, while for each negative value there are only 
two rows. This is because the Stokes discontinuity term is zero for negative values of $\delta$. Hence there is no need to 
display it. The first row for each value of $\delta$ represents the value obtained by using the LogGamma routine in Mathematica 
and is denoted by the row called LogGamma[z] in the Method column. Depending upon whether $\delta$ is positive or not,
the second row presents the Stokes discontinuity term. In general, this term was found to possess real and imaginary parts
of the order of $10^{-8}$ or a couple of orders lower. The next value for each $\delta$ is labelled either 1st AF 
or 3rd AF in the Method column. This means that depending on the sign of $\delta$ the first or third asymptotic form in Eq.\ 
(\ref{fiftyeight}) was used to calculate the value of $\ln \Gamma(z)$. The values of the truncated sum, the regularized value 
of the remainder and the Stirling approximation are not displayed due to limited space. 

It should also be noted that when $|\delta|$ is extremely small, e.g.\ $10^{-5}$, NIntegrate experiences convergence problems 
because for such values of $\delta$, it is carrying out the integration too close to the singularities lying on the Stokes line. 
In particular, for $\delta = 10^{-5}$, the program printed out a value of $\ln \Gamma(z)$ that agreed with the actual value 
to 25 decimal places for the real part, but in the case of the imaginary part the results only agreed to 18 decimal places.  
Although this calculation is not presented in the table, it still represents a degree of success since the imaginary part of the 
Stokes discontinuity term is of the order of $10^{-12}$. That is, the Stokes discontinuity term had to be correct to the first six
decimal places in order to yield the imaginary part of $\ln \Gamma(z)$ for $\delta \!=\! 10^{-5}$.

With the exception of the first value of $\delta$, which reflects the situation as the error function begins to veer away from
the step-function, we expect for all other values of $\delta$ that the Stokes multiplier $S^{+}$ is close to 1/2 according
to Fig.\ \ref{figone} assuming that smoothing occurs. For those results where $\delta >0$, this means
that the first asymptotic form with only half the Stokes discontinuity term should be a far more accurate approximation
to the actual value of $\ln \Gamma(z)$ than the entire asymptotic form in Eq.\ (\ref{fiftyeight}). However, we see the opposite 
where the first asymptotic form yields the exact value of $\ln \Gamma(z)$ for all values of $\delta$ despite the fact that the
Stokes discontinuity term has no effect on the first nine decimal places. For $\delta < 0$, if smoothing does occur, then the 
third asymptotic form in Eq.\ (\ref{fiftyeight}) should not yield exact values of $\ln \Gamma(z)$ because it is missing almost 
half the Stokes discontinuity term. Yet, we see the opposite that the third asymptotic form yields exact values of $\ln \Gamma(z)$ 
for all negative values of $\delta$ in the table. Thus, we have seen clearly that there is no smoothing of the Stokes phenomenon 
occurring the vicinity of the Stokes line at $\theta = \pi/2$ as postulated in Refs. \cite{ber89} and \cite{olv90}. For further
discussion of this issue, the reader should consult Ref.\ \cite{kow14b}. 

An explanation as to why there is no smoothing of the Stokes phenomenon is given in Sec.\ 6.1 of Ref.\ \cite{kow09}, 
where it is shown that the form for the Stokes multiplier as proposed by Berry and Olver is based on applying standard 
asymptotic techniques to its integral form. Ref.\ \cite{olv90}, which is regarded as a ``rigorous proof" that smoothing
occurs in the vicinity of a Stokes line, is based on truncating at a few orders the expression obtained by 
the application of Laplace's method. Consequently, this rigorous analysis is littered with the Landau gauge symbol $O()$,
$+ \cdots$ and tildes. Since only the lowest order terms are retained, Olver obtains the error function result given
by ``Eq." (\ref{sixtyseven}). However, the neglected terms are not only divergent; they are extremely difficult to 
regularize. If they could be regularized, then they would produce the necessary corrections to turn the error function 
in Fig.\ \ref{figone} into the step-function as postulated in the conventional view of the Stokes phenomenon. This also
vindicates the statement made at the beginning of this section that a proper numerical study can be far more powerful than 
a ``proof" in asymptotics. 

Since we have seen that the asymptotic forms in Thm.\ 2.1 give exact values of $\ln \Gamma(z)$ for all values of $z$,
we can differentiate both sides, thereby obtaining asymptotic forms for the digamma function, $\psi(z)$. This represents
another advantage in deriving complete asymptotic expansions over standard Poincar${\acute{\rm e}}$ asymptotics because according 
to p.\ 153 of Ref.\ \cite{whi73}, it is generally not permissible to differentiate an asymptotic expansion when adopting 
the latter prescription. Therefore, we arrive at the following theorem:

\begin{theorem}
The digamma function $\psi(z)$ possesses the following asymptotic forms:
\begin{align} 
\psi(z) & =  \ln z -  \frac{1}{2z} - \sum_{k=1}^{N-1} \frac{(-1)^k}{(2z)^{2k}} \, \Gamma(2k)\, c_k(1)  
+  R\psi^{SS}_N(z) + SD\psi^{SS}_M(z)  \;,
\label{sixeightone}\end{align}
where the remainder $R\psi^{SS}_N(z)$ is given by
\begin{align}
R\psi_N^{SS}(z) & =  \Bigl( \frac{3-2N}{z}\Bigr)  R^{SS}_N(z) + \frac{4\,(-1)^{N}}{ (2\pi z)^{2N-4}} \int_0^{\infty} dy \; y^{2N-2}\, e^{-y}  
\nonumber\\
& \times \;\;  \sum_{n=1}^{\infty} \frac{1}{n^{2N-4} \left( y^2 + 4 \pi^2 n^2 z^2 \right)^2} \;,
\label{sixeighttwo}\end{align}
$R^{SS}_N(z)$ is given by Eq.\ (\ref{twentytwo}), and the Stokes discontinuity term $SD\psi_M(z)$ is given by
\begin{align}
SD\psi^{SS}_M(z) & =  \Bigl( \mp 2 \lfloor M/2 \rfloor   \mp \frac{\left( 1 - (-1)^M \right)}{ (1- e^{\mp 2i \pi z}) }\Bigr) \pi i \;. 
\label{sixeightthree}\end{align} 
Eq.\ (\ref{sixeighttwo}) is valid for either $(M-1/2) \pi \!<\! \theta \!=\! {\rm arg}\, z \!<\! (M+1/2) \pi$ 
or $-(M+1/2) \pi \!<\! \theta \!<\! -(M-1/2) \pi$, where $M$ is a non-negative integer. As in Thm.\ 2.1, the Stokes discontinuity term 
or Eq.\ (\ref{sixeightthree}) has two versions, which are complex conjugates of one another. The upper-signed version is valid for 
$(M-1/2) \pi \!<\! \theta \!<\! (M+1/2) \pi$, while the lower-signed version is valid over $-(M+1/2) \pi \!<\! \theta  \!<\! -(M-1/2) \pi$.
For the situation along the Stokes lines or rays, where $\theta \!=\! \pm (M+1/2) \pi$, we replace $R\psi^{SS}_N(z)$ and $SD\psi^{SS}_M(z)$
by $R\psi^{SL}_N(z)$ and $SD\psi^{SL}_M(z)$, respectively. Then the remainder for these particular values of $\theta$ is 
\begin{align}
R\psi_N^{SL}(z) & = \Bigl(\frac{3-2N}{z}\Bigr) R^{SL}_N(z) +\frac{4(-1)^N }{(2 \pi |z|)^{2N-4}} \, P \int_0^{\infty} dy \; y^{2N-2}\, e^{-y} 
\nonumber\\
& \times \;\;\sum_{n=1}^{\infty} \frac{1}{n^{2N-4}(y^2 -4n^2 \pi^2 |z|^2)^2} \;\;,
\label{sixeightfour}\end{align}
where $R^{SL}_N(z)$ is given by Eq.\ (\ref{twentythree}). On the other hand, the Stokes discontinuity term is given by
\begin{align}
SD\psi^{SL}_{M}(z) & =  \pm \Bigl( (-1)^M -1 -2\lfloor M/2 \rfloor  +  \frac{(-1)^M }{ e^{2 \pi |z|} -1} \Bigr) \pi i \;. 
\label{sixeightfive}\end{align}
\end{theorem}

{\bfseries Proof}. The above results have been obtained by taking the derivative of the various results in Thm.\ 2.1 and carrying
out some elementary algebraic manipulation, which is left as an exercise for the reader. The results along the Stokes lines have
been obtained by writing $|z|$ as $z \exp(\mp (M+1/2) \pi i)$. This completes the proof.

From Thm.\ 3.1 we see that the remainder is considerably more complicated than for the asymptotic forms of $\ln \Gamma(z)$. However,
the second and third programs in the appendix can be adapted to evaluate the remainder in Thm.\ 3.1. Moreover, it can be seen why it is often
not permissible to take the derivative of the resulting truncated asymptotic expansion from the application of the Poincar${\acute{\rm e}}$ 
prescription. This is because taking the derivative of the neglected remainder can yield terms that are comparable to the leading order
terms. Finally, as a result of Thm.\ 3.1, we can continue to take derivatives of the digamma function, thereby obtaining asymptotic forms 
for $\psi^{(n)}(z)$, which in turn equals $(-1)^{n+1} n! \, \zeta(n+1,z)$, where $\zeta(n,z)$ is Hurwitz zeta function. 

\section{Mellin-Barnes Regularization}
In the preceding section we were able to exactify Stirling's approximation by carrying out a spectacular numerical 
study of the asymptotic forms given by Eq.\ (\ref{fiftyeight}), which in turn were derived from the results given in
Theorem\ 2.1. There are, however, two drawbacks with the numerical study in the previous section. The first is 
that we need to place an upper limit on the infinite sums appearing in the expressions for the regularized value of 
the remainder of the asymptotic series in Eq.\ (\ref{sixteen}). Thus, in order to evaluate the remainders, we effectively 
truncated the sums, although in this instance they were convergent and not asymptotic as in the case of Eq.\ (\ref{sixteen}).
In order to evaluate the regularized values extremely accurately, an upper limit of $10^5$ in the sum over $n$ was chosen. This 
results in the second drawback, which is the considerable computation required to calculate the regularized values of the remainder. 
The second drawback arises from the fact in order to obtain extremely accurate values of the remainder, either $10^5$ calls to 
the NIntegrate routine in Mathematica in the second and third programs in the appendix or $2 \times 10^5$ calls to the routine 
for the incomplete gamma function in the first program had to be made. Therefore, much time was expended in obtaining all the 
results presented in the preceding section. Moreover, we do not wish to truncate any result here so that we can once and 
for all dispel any doubt that we are evaluating an approximation. If the infinite sum over $n$ in the regularized values of 
the remainder can be replaced by a single result, then there will be a huge reduction in the execution times since only one 
call to the NIntegrate routine would be required.

The problems mentioned above were similar to those encountered when carrying out a numerical investigation of
the Borel-summed forms of the complete asymptotic expansion for a particular value of the generalized Euler-Jacobi 
series, viz.\ $S_3(a) =\sum_{n=1}^{\infty}\exp(-a n^3)$, in Ref.\ \cite{kow95}. There, the Borel-summed regularized 
values of the remainder were even more problematic than those presented in the previous sections because they not
only involved an infinite sum, but the integrals were also two-dimensional. To avoid dealing with such computationally 
intractable situations, the technique of Mellin-Barnes regularization was devised. Since then, the 
technique has been developed further and applied to various problems in Refs.\ \cite{kow002}, \cite{kow09}-\cite{kow11c} 
and \cite{kow11}. In this section we aim to apply the technique to the asymptotic series in Eq.\ (\ref{sixteen}), 
thereby obtaining an another set of regularized asymptotic forms yielding $\ln \Gamma(z)$. As we shall see, these forms are very 
unlike the corresponding Borel-summed asymptotic forms and thus, constitute a different method for evaluating
the regularized value.

\begin{theorem} Via the Mellin-Barnes (MB) regularization of the asymptotic series $S(z)$ given by Eq.\ (\ref{twenty}), the
logarithm of the gamma function or $\ln \Gamma(z)$ can be expressed as
\begin{align}
\ln \Gamma(z) & = \Bigl( z- \frac{1}{2} \Bigr) \ln z -z + \frac{1}{2} \,\ln (2 \pi) +
z \sum_{k=1}^{N-1} \frac{(-1)^k}{(2z)^{2k}}\, \Gamma(2k-1)\,c_k(1)
\nonumber\\
& - \;\; 2z \!\!\!\!\! \!\!\!\!\! \int\limits_{\substack{c-i \infty \\ {\rm Max}[N-1,1/2]<c=\Re\, s<N}}^{c+i \infty} 
\!\!\!\!\! \!\!\!\!\!\!\!\!  ds\; \left( \frac{1}{2\pi z} \right)^{2s} \frac{e^{\pm 2 M i \pi s}}{e^{-i \pi s} -e^{i \pi s}} 
\; \zeta(2s) \, \Gamma(2s-1) + S_{MB}(M,z) \;,
\label{sixtyeighta}\end{align}
where 
\begin{align}
S_{MB}(M,z) & = \pm \lfloor M/2 \rfloor \ln \Bigl( -\, e^{- 2i \pi z} \Bigr) - \Bigl(\frac{1-(-1)^M}{2} \Bigr)  
\ln \Bigl( 1- e^{\pm 2i \pi z} \Bigr)\;\;, 
\label{sixtyeightb}\end{align}
for $(\pm M-1)\pi< \theta= {\rm arg}\, z< (\pm M+1)\pi$ and $M \geq 0$, but excluding $\theta$ equal to half-integer values of $\pi$. 
Within each  domain of convergence there are two half-integer cases of $S_{MB}(M,z)$, viz.\ $\theta= (\pm M-1/2) \pi$ and $\theta=
(\pm M +1/2) \pi$. For $\theta = \pm(M-1/2)\pi$, $S_{MB}(M,z)$ is given by 
\begin{align}
S_{MB}(M,z) & = \Bigl[ \Bigl( \frac{(-1)^{M}-1}{2} \Bigr) \ln \Bigl(1 - e^{-2 \pi |z|}\Bigr)  + 2 (-1)^{M+1} \, 
\lfloor M/2 \rfloor \pi \, |z| \Bigr] \;\;,
\label{sixtyeightc}\end{align}
while for $\theta= \pm (M+1/2) \pi$, it is given by
\begin{align}
S_{MB}(M,z) & = \Bigl[ \Bigl( \frac{(-1)^{M}-1}{2} \Bigr) \ln \Bigl(1 - e^{-2 \pi |z|}\Bigr)  + 2 \pi |z| 
\nonumber\\
& \times \;\; \Bigl( (-1)^M \lfloor M/2 \rfloor + \frac{(-1)^M-1}{2} \Bigr) \Bigr] \;\;.
\label{sixtyeightd}\end{align}
\end{theorem}

\begin{remark}
Note that Eqs.\ (\ref{sixtyeightc}) and (\ref{sixtyeightd}) only apply when $M$ is a positive integer, i.e. $S_{MB}(0,z)$ vanishes.
In addition, the logarithmic term in both forms of $S_{MB}(M,z)$ vanishes for even integer values of $M$.
\end{remark}
\begin{remark}
For each value of $M$ the sector over which Eq.\ (\ref{sixtyeightb}) is defined represents the domain of convergence for 
the MB integral in Eq.\ (\ref{sixtyeighta}). Each domain of convergence contains the lines where the extra terms
given by Eqs.\ (\ref{sixtyeightc}) and (\ref{sixtyeightd}) are valid except for the $M\!=\! 0$ domain of convergence. These 
lines do not represent the boundaries of Stokes sectors as observed in Thm.\ 2.1, although they occur at the same locations as 
the Stokes lines in Borel summation. Here they need to be isolated as a result of the MB regularization of $S(z)$ 
since we have already seen that $\ln \Gamma(z)$ possesses jump discontinuities at $\theta= (l+1/2)\pi$, where $l$ can be
any integer. By definition, MB regularization yields an alternative representation of the original function via its asymptotic 
expansion, and relies on the continuity of the function. If the original function possesses discontinuities as in the case of 
$\ln \Gamma(z)$, then the MB-regularized value will not yield the value of the function, although the analysis can be adapted to 
obtain the correct value of the function as described in the proof. 
\end{remark} 
\begin{remark}
Since MB regularization is unable to yield the discontinuities of $\ln \Gamma(z)$, i.e. they do not appear in the MB integral
in Eq.\ (\ref{sixtyeighta}), they can only arise from the infinite series obtained by summing the contributions due to the 
infinite number of singularities lying on the Stokes lines. It is the regularization of this series that leads to the discontinuities
of $\ln \Gamma(z)$. 
\end{remark}
{\bfseries Proof}. The MB regularization of both types of generalized terminants is discussed in Ch.\ 9 of Ref.\ \cite{kow09}. 
Because the series $S(z)$ represents a specific case of a Type I generalized terminant as defined by Eq.\ (\ref{nineteen}), we 
only require the MB-regularized value of a Type I generalized terminant, which is presented as Proposition\ 4 on p.\ 139
of the reference. Hence we find that 
\begin{align}
S^{I}_{p,q} \left( N, z^{\beta} \right) &\equiv \int\limits_{\substack{c-i \infty \\ {\rm Max}[N-1,-q/p]<c=\Re\, 
s<N}}^{c+i \infty} \!\!\!\!\!  ds\; \frac{z^{\beta s}\,e^{\mp 2 M i \pi s}}{e^{-i \pi s} -e^{i \pi s}} \;\Gamma(ps+q) 
\mp \frac{2 \pi i}{p} \;  z^{-\beta q/p}
\nonumber\\
& \times \;\; e^{\mp  q i \pi/p}  \sum_{j=1}^M e^{\pm 2j q i\pi/p} \, \exp\Bigl( -z^{-\beta/p} e^{\pm (2j-1) i\pi/p} 
\Bigr) \;\;,
\label{sixtynine}\end{align}
for $(\pm 2M-1-p/2)\pi/\beta < {\rm arg}\, z < (\pm 2M+1+p/2) \pi/\beta$, $M>0$ and $N> -q/p$.  The last condition is 
required to avoid the poles due to the gamma function in the integrand. For $N< -q/p$, however, we can separate those 
terms in $S^{I}_{p,q}(N, z^{\beta})$ up until the first value of $k$, where $pk+q$ is greater than zero and then 
re-adjust $N$, thereby allowing us to use Equivalence\ (\ref{sixtynine}) again. Substituting $z$, $\beta$, $p$ and $q$ 
in the above equivalence respectively by $1/2n \pi z$, 2, 2 and -1, we find that the MB-regularized value of the 
generalized terminant in Eq.\ (\ref{twenty}) is given by
\begin{align}
S^{I}_{2,-1} \left( N, (1/2n \pi z)^2 \right) &\equiv \int\limits_{\substack{c-i \infty \\ {\rm Max}[N-1,1/2]<
c=\Re\, s<N}}^{c+i \infty} \!\!\!\!\!  ds\; \left( \frac{1}{2n \pi z} \right)^{2s} \frac{e^{\mp 2  M i \pi s}}
{e^{-i \pi s} -e^{i \pi s}} \;\Gamma(2s-1) 
\nonumber\\
& + \;\;  \frac{1}{2n z}  \sum_{j=1}^M (-1)^j \, \exp\Bigl( \pm 2 (-1)^j n i\pi z \Bigr) \;\;,
\label{sixtyninea}\end{align}
where $(\mp M -1) \pi< \theta< (\mp M+1) \pi$. If we compare this result with its Borel-summed analogue given by 
Equivalence\ (\ref{twentyfivea}), we see that a Mellin-Barnes integral has replaced the Cauchy integral in the latter, 
while the Stokes discontinuity term appears to have been retained. However, the second term on the rhs of 
Equivalence\ (\ref{sixtyninea}), whilst resembling the Stokes discontinuity term, is nothing of the sort. This is 
because Equivalence\ (\ref{sixtyninea}) is valid over domains of convergence given by $(\mp M-1 ) \pi < \theta< 
(\mp M+1) \pi$. Furthermore, adjacent domains of convergence overlap each other. E.g., for $M=0$, the above result 
is valid for $-\pi< \theta<\pi$, while for $M=1$, it is valid for $0< \theta< 2 \pi$. That is, both the $M=0$ 
and $M=1$ forms of Equivalence\ (\ref{sixtyninea}) apply over $0 < \theta<\pi$. The extra term on the rhs of Eq.\ 
(\ref{sixtyninea}) arises when the MB-regularized value for one domain of convergence is set equal to the 
MB-regularized value of an adjacent domain of convergence in the overlapping sector. See Ref.\ \cite{kow09} for 
more details. Therefore, although the exponential term is similar to the Stokes discontinuity term in Eq.\ 
(\ref{thirtynine}), the above equivalence possesses no line of discontinuity that acts as a boundary 
between adjacent sectors as we observed in Thm.\ 2.1.

Now we introduce Equivalence\ (\ref{sixtyninea}) into Eq.\ (\ref{twenty}). This yields
\begin{align}
S(z) \; & \equiv  z \sum_{k=1}^{N-1} \frac{(-1)^k}{(2z)^{2k}}\, \Gamma(2k-1)\,c_k(1) -
\sum_{n=1}^{\infty} \frac{1}{n}  \sum_{j=1}^M (-1)^j \, \exp\Bigl( \pm 2 (-1)^j n i\pi z \Bigr) 
\nonumber\\  
& - \;\; 2z \!\!\!\!\! \int\limits_{\substack{c-i \infty \\ {\rm Max}[N-1,1/2]<c=\Re\, s<N}}^{c+i \infty} 
\!\!\!\!\! \!\!\!\!\!\!\!\!  ds\; \left( \frac{1}{2\pi z} \right)^{2s} \frac{e^{\mp 2 M i \pi s}}
{e^{-i \pi s} -e^{i \pi s}} \; \zeta(2s) \, \Gamma(2s-1)  \;\;,
\label{seventy}\end{align}  
where $(\mp M-1)\pi< \theta< (\mp M+1)\pi$ and $M>0$. In the above result the sum over $n$ in the MB integral 
has been replaced by the Riemann zeta function $\zeta(s)$. However, as in the case of Equivalence\ 
(\ref{twentysix}), the sum over $n$ involving the exponential terms  can be divergent. Therefore, we need to 
replace it by its regularized value, which can be determined again by using Lemma\ 2.2. In particular, by 
expressing the difference of the logarithms, $ \ln \left((1-e^{-2i\pi z})/(1-e^{2i \pi z}) \right)$,
as $\ln \left(-\, e^{-2i \pi z}\right)$, we find after a little algebra that
\begin{align}
\sum_{n=1}^{\infty} \frac{1}{n}  \sum_{j=1}^M (-1)^j \, \exp\Bigl( \pm 2 (-1)^j n i\pi z \Bigr) & \equiv
\pm \lfloor M/2 \rfloor \ln \Bigl( -\, e^{-2i \pi z} \Bigr) +\frac{\left(1-(-1)^M \right)}{2}\; 
\nonumber \\
& \times \;\; \ln \Bigl( 1- e^{\mp 2i \pi z} \Bigr) \;\;, 
\label{seventyone}\end{align}  
where $(\mp M -1)\pi < \theta< (\mp M +1) \pi$. Note that the final term contributes only when $M$ is odd. 

There is, however, a problem with the above result. When $\theta = (\mp M -1/2) \pi$ or $\theta=(\mp M+1/2) \pi$, 
the series resulting from summing an infinite number of residue contributions as represented by the second term 
on the rhs of Equivalence\ (\ref{sixtyninea}) is real. Yet the regularized value can become complex as 
witnessed by the first term on the rhs of Equivalence\ (\ref{seventyone}). Furthermore, this term is ambiguous 
because the imaginary part of the logarithm of a negative real number is equal to $(2j+1) \pi$, where $j$ can be 
any integer. Therefore, we need to invoke the Zwaan-Dingle principle \cite{din73}, which states that an initially 
real-valued function cannot suddenly acquire imaginary terms. For $\theta= \pm (M-1/2) \pi$ and $M > 0$, this 
means that we must take the real part of the regularized value given above. As a consequence, we find that 
$\ln \Gamma(z)$ becomes discontinuous for these values of $\theta$, while the sum becomes
\begin{align}
- \sum_{n=1}^{\infty} \frac{1}{n}  \sum_{j=1}^M (-1)^j \, \exp\Bigl( -2 (-1)^{M-j} n \pi |z| \Bigr) & \equiv
\Re \Bigl[ (-1)^M \lfloor M/2 \rfloor \ln \Bigl( -\, e^{- 2 \pi |z|} \Bigr) +\frac{\left((-1)^M -1 \right)}{2}\; 
\nonumber \\
& \times \;\; \ln \Bigl( 1- e^{-2\pi |z|} \Bigr) \Bigr] \;\;. 
\label{seventyonea}\end{align}
In obtaining this result both cases of $M$ being odd and even have been considered separately. In addition, 
the sum yields two separate sums, one in terms of growing exponentials, viz.\ $\exp(2 n \pi |z|)$ and the other 
in terms of decaying exponentials, $\exp(-2n \pi |z|)$. The first of these sums is divergent and must be regularized. 
Therefore, Lemma\ 2.2 has been employed again in arriving at the result on the rhs. However, note that it is still 
necessary to take the real part on the rhs because the logarithm of the first term on the rhs yields a complex value.   

For $\theta=\pm (M+1/2) \pi$, we adopt a similar procedure since the sum over $j$ is now over
$\exp \left( 2(-1)^{M-j} n \pi |z|) \right)$. Then we find that the sum over the residue contributions can be 
expressed as
\begin{align}
- \sum_{n=1}^{\infty} \frac{1}{n}  \sum_{j=1}^M (-1)^j \, \exp\Bigl( 2 (-1)^{M-j} n \pi |z| \Bigr) & \equiv
\Re \Bigl[ \Bigl( (-1)^{M+1} \lfloor M/2 \rfloor+ \frac{1-(-1)^M}{2} \Bigr) \ln \Bigl( -\, e^{- 2 \pi |z|} \Bigr) 
\nonumber\\
& + \;\; \frac{\left((-1)^M-1 \right)}{2}\, \ln \Bigl( 1- e^{-2 \pi |z|} \Bigr) \Bigr] \;\;. 
\label{seventyoneb}\end{align}
In obtaining this result the following identity has been used
\begin{eqnarray}
\Re \Bigl[ \ln \Bigl( 1- e^{2 \pi |z|}\Bigr) \Bigr]  = 2 \pi |z| + \ln \Bigl( 1-e^{-2 \pi |z|}\Bigr) \;\;.
\label{seventyonec}\end{eqnarray}
Consequently, taking the real part of the rhs of Equivalence\ (\ref{seventyoneb}) only applies to the first term,  
which can, in turn, be simplified according to 
\begin{align}
\Re \Bigl[ \ln \Bigl(-e^{\pm 2 \pi |z|} \Bigr) \Bigr] & = \pm  2  \pi |z| \;\;.
\label{seventyoned} \end{align}
Therefore, Eq.\ (\ref{seventyoneb}) reduces to
\begin{align}
- \sum_{n=1}^{\infty} \frac{1}{n}  \sum_{j=1}^M (-1)^j \, \exp\Bigl( 2 (-1)^{M-j} n \pi |z| \Bigr) & \equiv
\frac{\left(  (-1)^{M} -1\right)}{2} \ln \Bigl( 1 - e^{-2\pi |z|} \Bigr) 
\nonumber\\
& +\;\; \Bigl( (-1)^{M} \lfloor M/2 \rfloor  +\frac{(-1)^M -1}{2} \Bigr) 2 \pi |z| \;\;.
 \label{seventyonee}\end{align}

With the aid of Equivalence (\ref{seventyone}), we can express Equivalence\ (\ref{seventy}) as
\begin{align}
S(z)  & \equiv z \sum_{k=1}^{N-1} \frac{(-1)^k}{(2z)^{2k}}\, \Gamma(2k-1)\,c_k(1)
- 2z \!\!\!\!\! \!\!\!\!\! \int\limits_{\substack{c-i \infty \\ {\rm Max}[N-1,1/2]<c=\Re\, s<N}}^{c+i \infty} 
\!\!\!\!\! \!\!\!\!\!\!\!\!  ds\; \left( \frac{1}{2\pi z} \right)^{2s} \frac{e^{\mp 2 M i \pi s}}
{e^{-i \pi s} -e^{i \pi s}} 
\nonumber\\
& \times \;\; \zeta(2s) \, \Gamma(2s-1) + S_{MB}(M,z) \;\;,
\label{seventytwo}\end{align}
where 
\begin{align}
S_{MB}(M,z)& = 
\pm  \lfloor M/2 \rfloor \ln \Bigl( -\, e^{- 2i \pi z} \Bigr) -  \frac{\left(1-(-1)^M \right)}{2}  
\ln \Bigl( 1- e^{\pm 2i \pi z} \Bigr) \;\;.
\label{seventytwoa}\end{align}  
The above result is valid for $(\pm M -1)\pi < \theta< (\pm M +1) \pi$, but excludes half-integer values 
of $\pi$ within this domain. For the specific case of $\theta = \pm (M-1/2) \pi$, we require Equivalence\
(\ref{seventyonea}). Then $S_{MB}(M,z)$ is given by
\begin{align}
S_{MB}(M,z)& = 2 (-1)^{M+1} \lfloor M/2 \rfloor\, \pi |z| +  \frac{\left((-1)^M -1\right)}{2}  
\ln \Bigl( 1- e^{- 2\pi |z|} \Bigr) \;\;.
\label{seventytwob}\end{align}  
On the other hand, for $\theta = \pm (M+1/2) \pi$, we require Equivalence (\ref{seventyonee}), which means
\begin{align}
S_{MB}(M,z)& = 2 \pi |z| \Bigl( (-1)^M \, \lfloor M/2 \rfloor + \frac{(-1)^M-1}{2} \Bigr) +  
\frac{\left((-1)^M-1 \right)}{2} \ln \Bigl( 1- e^{- 2\pi |z|} \Bigr) \;\;.
\label{seventytwoc}\end{align}  

The series $S(z)$ appears on the rhs of Equivalence\ (\ref{ten}), which means in turn that the lhs of 
this statement also represents the regularized value of the series. Since the regularized value is unique, 
the rhs of the Equivalence\ (\ref{seventytwo}) is, therefore, equal to the lhs of Equivalence\ (\ref{ten}). 
Hence we finally arrive at the results given in Theorem 4.1, which completes the proof.

It should be noted that $S_{MB}(0,z)$ vanishes for $\theta =\pm \pi/2$. That is, the regularized value 
of the remainder is given by the MB integral in Eq.\ (\ref{sixtyeighta}), whose domain of convergence 
is $-\pi< \theta <\pi$. As a result, there are no discontinuities when $\theta =\pm \pi/2$ in 
$\ln \Gamma(z)$, although these lines are Stokes lines. Hence we see that Stokes lines/rays of 
discontinuities do not necessarily imply that the original function is discontinuous along them. 

By comparing the results in Theorem\ 4.1 with the Borel-summed results in Theorem\ 2.1, we see that 
not only is the remainder of $S(z)$ different in that it is now expressed in terms of an MB integral, 
but there are no Stokes discontinuity terms arising out of the traversal of Stokes lines. Instead, the 
MB integral is valid over a sector, which represents its domain of convergence. The Stokes lines occurring 
in the Borel-summed results are now situated within the domains of convergence. There are also no 
Stokes multipliers in MB regularization. In particular, there are no discontinuities at 
$\theta = \pm \pi/2$ as indicated in the results of Theorem\ 2.1. That is, they are fictitious, an 
artefact of Borel summation. In fact, discontinuities only occur on Stokes lines if the original function 
possesses singularities on them, which is the case for $\ln \Gamma(z)$ at $\theta = \pm (l+1/2) \pi$, 
where $l>0$. Nevertheless, whilst there are no Stokes discontinuities in the MB-regularized results, in 
order that the MB-regularized results agree with each other where the domains of convergence overlap, 
there are extra logarithmic terms appearing in the MB-regularized value of $\ln \Gamma(z)$ given in 
Theorem\ 4.1, which are similar in form to the Stokes discontinuity terms in Theorem\ 2.1. 

Another feature of the results in Theorem\ 4.1 is that the sum over $n$ in the regularized value 
appearing in Theorem\ 2.1 has vanished. It has effectively been replaced by the Riemann zeta function. 
As a consequence, we only have one integral to evaluate in the remainder of $S(z)$, but at the same time, 
we need to ensure that the software package we use is able to evaluate the zeta function extremely 
accurately. Fortunately, Mathematica \cite{wol92} is capable of doing this via its Zeta routine.

It should also be pointed out that equating forms for the regularized value where they apply in a domain 
of convergence can produce new Mellin transform pairs. For example, if we put $M=0$, Eq.\ (\ref{sixtyeighta}) 
is valid for $-\pi < \theta< \pi$, while if $M = 1$, then it is valid for $0< \theta<2\pi$. Hence we can 
equate the two forms of the regularized value for $\pi< \theta< 2\pi$. As a consequence, we arrive at 
\begin{eqnarray}
\frac{1}{2 \pi i} \!\!\!\!\! \!\!\!\!\! \int\limits_{\substack{c-i \infty \\ {\rm Max}[N-1,1/2]<c=\Re\, 
s<N}}^{c+i \infty} \!\!\!\!\! \!\!\!\!\!\!\!\!  ds\; \left(2 \pi z\right)^{-2s} e^{i \pi s} \, \zeta(2s) 
\, \Gamma(2s-1) = - \frac{1}{4 \pi i z} \, \ln \Bigl( 1 - e^{2i \pi z} \Bigr) \;\;. 
\label{seventythree}\end{eqnarray}
Substituting $\sqrt{\smash[b]{y}}= 2\pi z \exp(-i \pi/2)$ into the above result yields
\begin{eqnarray}
\frac{1}{2 \pi i} \!\!\!\!\! \!\!\!\!\! \int\limits_{\substack{c-i \infty \\ {\rm Max}[N-1,1/2]<c=\Re\, 
s<N}}^{c+i \infty} \!\!\!\!\! \!\!\!\!\!\!\!\!  ds\; y^{-s} \zeta(2s) \, \Gamma(2s-1)
= \frac{1}{2 \sqrt{y} } \, \ln \Bigl( 1 - e^{-\sqrt{y}}  \Bigr) \;\;. 
\label{seventyfour}\end{eqnarray}
The lhs in Eq.\ (\ref{seventyfour}) is now in the form of an inverse Mellin transform, which means 
alternatively that
\begin{eqnarray}
\int_0^{\infty} dz \; z^{s-3/2} \ln \Bigl( 1- e^{-\sqrt{z}} \Bigr)= 2 \, \zeta(2s) \, \Gamma(2s-1) \;\;,
\label{seventyfive}\end{eqnarray}
for $\Re\, s > 1/2$ and $|{\rm arg}\, z|< \pi$.

Although the results in Theorem\ 4.1 have been proven, as in the case of Theorem\ 2.2, we cannot be 
certain that they are indeed valid because it has already been observed in the case of ``Stokes 
smoothing" that proofs in asymptotics are not reliable unless they are validated by an effective 
numerical analysis such as that in Sec.\ 3. Since the results in Theorem\ 2.2 have been validated, we can
use them to establish the validity of the MB-regularized forms given in Theorem\ 4.1. Therefore, in 
the next section we present another numerical analysis in which the results obtained from the 
MB-regularized forms for $\ln \Gamma(z^3)$ in Theorem\ 4.1 are compared with the corresponding 
Borel-summed asymptotic forms in Sec.\ 2.

\section{Further Numerical Analysis}
The numerical analysis presented in Sec.\ 3 was concerned with evaluating $\ln \Gamma(z)$ by using the asymptotic forms derived from
regularizing the asymptotic series $S(z)$ given by Eq.\ (\ref{sixteen}) via Borel summation. In carrying out the analysis, we were 
restricted to considering the particular forms which were valid only over the principal branch of the complex plane for $z$. This 
restriction was due to the fact that Mathematica \cite{wol92} evaluates $\ln \Gamma(z)$ for $z$ lying in the principal branch of 
the complex plane. Consequently, we were unable to determine whether the asymptotic forms in Theorem\ 2.1 for higher/lower Stokes sectors, 
i.e., those where $M \geq 2$, were indeed correct. 

According to the definition of the regularized value \cite{kow95}, \cite{kow002}, \cite{kow09}-\cite{kow11c}, it must be invariant 
irrespective of the method used to evaluate it. Therefore, if we can demonstrate that the MB-regularized asymptotic forms of the regularized value 
yield identical values to the Borel-summed ones presented in Sec.\ 2, then we can be satisfied that the results of the previous section 
are indeed correct, especially for those Stokes sectors and lines not studied in Sec.\ 3. In order to access the higher/lower sectors or lines we 
now consider powers of the variable $z$, viz.\ $z^3$, in $\ln \Gamma(z)$. That is, we are assuming that there is a solution to 
a problem, $f(z)$, which happens to possess the asymptotic forms of $\ln \Gamma(z^3)$. The principal branch is still $(-\pi,\pi]$, 
but Mathematica is only able to evaluate $\ln \Gamma(z^3)$ in the sector where $-\pi/3 < {\rm arg}\,z \leq \pi/3$.
 
From Theorem\ 4.1 we find that there can be two different representations for the regularized value of $\ln \Gamma(z)$ since 
replacing $M$ by either $M-1$ or $M+1$ in Eq.\ (\ref{sixtyeighta}) gives a different form for the regularized value, but which is valid 
over the common half of the sector or domain of convergence for $M \!=\! M$. For example, the upper-signed version of Eq.\ (\ref{sixtyeighta})
is valid only for $\pi<\theta<3\pi$ when $M=2$, while for $M \!=\! 1$ and $M \!=\!3$, it is only valid for $0< \theta< 2 \pi$
and $2\pi < \theta< 4 \pi$, respectively. Thus, the $M=1$ result applies over the bottom half of the domain of convergence for the $M \!=\! 2$ result, 
while the $M \!=\! 3$ result applies over the top half of the domain of convergence for the $M \!=\! 2$ result. Consequently, we are not only in a 
position to evaluate $\ln \Gamma(z)$ for higher/lower arguments of $z$, but we can check the MB-regularized asymptotic forms against each 
other for those values of $M$ where the domains of convergence overlap one another. Moreover, the $M \!=\! 0$ results can be checked with the values 
of $\ln \Gamma(z^3)$ calculated by Mathematica. If this proves to be successful, then we can make another final check to observe whether the 
MB-regularized forms of $\ln \Gamma(z^3)$ yield identical values to those evaluated by the appropriate Borel-summed asymptotic forms in Sec.\ 2. 
Previously, we had no method of checking whether the Borel-summed asymptotic forms for $\ln \Gamma(z)$ outside the principal branch of the 
complex plane were indeed correct. We can now overcome this problem by replacing $z$ by $z^3$ and checking them against the values obtained
from the corresponding MB-regularized forms.

If we make the substitution, $z\!=\!z^3$, in Theorem\ 4.1, then we observe that all the results up and to including $M \!=\! 3$ cover a 
segment within the principal branch of the complex plane for $z$. When $M \!=\! 0$, which is valid for $-\pi/3 < \theta < \pi/3$, we find 
that either the upper- or lower-signed version of Eq.\ (\ref{sixtyeighta}) yields
\begin{eqnarray}
\ln \Gamma \left( z^3 \right)= F \left( z^3 \right)+  TS_N \! \left(z^3 \right)  - 
2z^3 \!\!\!\!\! \!\!\!\!\! \int\limits_{\substack{c-i \infty \\ {\rm Max}[N-1,1/2]<c=\Re\, s<N}}^{c+i \infty} 
\!\!\!\!\! \!\!\!\!\!\!\!\!  ds\; \frac{(1/2\pi z^3)^{2s}}{e^{-i \pi s} -e^{i \pi s}} 
\; \zeta(2s) \, \Gamma(2s-1) \;\;, 
\label{seventysix}\end{eqnarray}
where $TS_N(z)$ has been defined as the truncated part of the asymptotic series $S(z)$ at $N$ as in Eq.\ (\ref{fiftysevena}).
Hence we see that $\ln \Gamma(z^3)$ is composed of the Stirling approximation, $F(z)$, given by Eq.\ (\ref{fiftynine}),
the truncated series $TS_N(z^3)$, and a Mellin-Barnes integral, representing the regularized value of the remainder of $S(z)$ when 
it is truncated at $N$. 

Now, if we put $M=1$ in the upper-signed version of Eq.\ (\ref{sixtyeighta}), then we obtain
\begin{align}
\ln \Gamma \left( z^3 \right)& =F \left( z^3 \right) +  TS_N \! \left( z^3 \right)  - 
2z^3 \!\!\!\!\! \!\!\!\!\! \int\limits_{\substack{c-i \infty \\ {\rm Max}[N-1,1/2]<c=\Re\, s<N}}^{c+i \infty} 
\!\!\!\!\! \!\!\!\!\!\!\!\!  ds\; \frac{(1/2\pi z^3)^{2s}\,e^{2i \pi s}}{e^{-i \pi s} -e^{i \pi s}} 
\; \zeta(2s) \, \Gamma(2s-1) 
\nonumber\\
& - \;\; \ln \left( 1- e^{2 i \pi z^3} \right)\;\;. 
\label{seventyseven}\end{align}
The domain of convergence of the MB integral in this result is $0< \theta< 2\pi/3$, but the above result is not valid for $\theta=\pi/2$ 
according to Thm.\ 4.1. This is because $S_{MB}(M,z^3)$ is discontinuous whenever $\theta =\pm (M \pm 1/2) \pi/3$ except for $M \!=\!0$.
For $\theta =\pi/6$, we can actually use Eq.\ (\ref{sixtyeightc}), but all it does is replace the logarithmic term on the rhs of Eq.\ 
(\ref{seventyseven}) by $\ln \Bigl(1- e^{-2\pi |z|^3} \Bigr)$, again reinforcing the point that there is no discontinuity at $\theta \!=\! \pi/6$. 

When $M \!=\! 1$ in $\theta \!=\! \pm (M \!+\! 1/2) \pi/3$, we have $\theta \!=\! \pm \pi/2$. The upper value of $\theta$ lies in the 
the domain of convergence for $M \!=\! 1$ in the upper-signed version of Eq.\ (\ref{sixtyeighta}). In addition from Thm.\ 4.1, we replace
$S_{MB}(M,z)$ by Eq.\ (\ref{sixtyeightd}) with $z$ replaced by $z^3$ and $M \!=\! 1$. Consequently, we obtain 
\begin{align}
\ln \Gamma \left( z^3 \right)& =F \left( z^3 \right) +  TS_N \! \left( z^3 \right)  - 
2z^3 \!\!\!\!\! \!\!\!\!\! \int\limits_{\substack{c-i \infty \\ {\rm Max}[N-1,1/2]<c=\Re\, s<N}}^{c+i \infty} 
\!\!\!\!\! \!\!\!\!\!\!\!\!  ds\; \frac{(1/2\pi z^3)^{2s}\,e^{2i \pi s}}{e^{-i \pi s} -e^{i \pi s}} 
\; \zeta(2s) \, \Gamma(2s-1) 
\nonumber\\
& - \;\; 2 \pi |z|^3  - \;\; \ln \left( 1- e^{-2 \pi |z|^3} \right)\;\;. 
\label{seventysevenb}\end{align}
One consequence of the penultimate term is that we expect a discontinuity to emerge when the above result for $\ln \Gamma(z^3)$ is programmed as  
a Mathematica module later in this section. In addition, in Eq.\ (\ref{seventyseven}) we can replace $F(z^3)$ and $TS_N(z^3)$ by $F(-i|z|^3)$ and 
$TS_N(-i|z|^3)$, respectively, while $z^3$ in the term with the MB integral can be replaced by $-i |z|^3$.

When compared with the $M \!=\!0$ result for $\ln \Gamma(z^3)$ or Eq.\ (\ref{seventysix}), we see that Eqs.\ (\ref{seventyseven}) and 
(\ref{seventysevenb}) possess an extra term or terms with the MB integral. These are analogous to the Stokes discontinuity term in the Borel-summed 
asymptotic forms in Sec.\ 2, but the major difference here is that the lines of discontinuity are located inside the domains of convergence. 
As a consequence, the asymptotic form is only different on the lines, whereas in the case of Stokes lines, the regularized value is different 
before, on and after them. Moreover, where the domains of convergence overlap, we expect the forms for $\ln \Gamma(z^3)$ to yield identical 
values such as the common region of $0 < \theta <\pi/3$ for both the $M=0$ and $M=1$ results. This is simply not possible
with the Stokes phenomenon because the regions of validity for the asymptotic forms do not overlap. 

For $M \!=\! 2$ the upper-signed version of Eq.\ (\ref{sixtyeighta}) with $z$ replaced by $z^3$ becomes
\begin{align}
\ln \Gamma \left( z^3 \right) &=F \left( z^3 \right) + TS_N \! \left( z^3 \right)  - 
2z^3 \!\!\!\!\! \!\!\!\!\! \int\limits_{\substack{c-i \infty \\ {\rm Max}[N-1,1/2]<c=\Re\, s<N}}^{c+i \infty} 
\!\!\!\!\! \!\!\!\!\!\!\!\!  ds\; \frac{(1/2\pi z^3)^{2s}\,e^{4i \pi s}}{e^{-i \pi s} -e^{i \pi s}} 
\; \zeta(2s) \, \Gamma(2s-1)
\nonumber\\
& + \;\; \ln \left( -\, e^{-2 i \pi z^3}\right) \;\;, 
\label{seventyeight}\end{align}
while for $M \!=\! 3$, the upper-signed version of Eq.\ (\ref{sixtyeighta}) with $z$ replaced by $z^3$ reduces to 
\begin{align}
\ln \Gamma \left( z^3 \right) &= F \left( z^3 \right) +  TS_N \! \left( z^3 \right)  - 
2z^3 \!\!\!\!\! \!\!\!\!\! \int\limits_{\substack{c-i \infty \\ {\rm Max}[N-1,1/2]<c=\Re\, s<N}}^{c+i \infty} 
\!\!\!\!\! \!\!\!\!\!\!\!\!  ds\; \frac{(1/2\pi z^3)^{2s}\,e^{6i \pi s}}{e^{-i \pi s} -e^{i \pi s}} 
\; \zeta(2s) \, \Gamma(2s-1)
\nonumber\\
&  + \;\; \ln \left( -\, e^{-2 i \pi z^3}\right)  - \ln \left( 1- e^{2 i \pi z^3} \right) \;\;.  
\label{seventynine}\end{align}
Eqs.\ (\ref{seventyeight}) and (\ref{seventynine}) are only valid respectively for $\pi/3 < \theta < \pi$ and $2\pi/3 < 
\theta < 4 \pi/3$, except for $\theta \!=\! \pi/2$, $\theta \!=\! 5\pi/6$ and $\theta=7 \pi/6$. The last case can be discarded
since it lies outside of the principal branch for $z$. The above results are similar to Eq.\ (\ref{seventyseven}) 
except that the logarithmic terms are slightly different as a result of Eq.\ (\ref{sixtyeightc}). Note, however, that the $M\!=\!3$ 
result possesses the same logarithmic term as the $M \!=\! 1$ result plus (instead of minus) the extra logarithmic term in the 
$M \!=\! 2$ result.

For $M=1$ in $\theta \!=\! \pm (M \!+\! 1/2)\pi/3$, we used Eq.\ (\ref{sixtyeightd}) to derive the asymptotic form of $\ln \Gamma(z^3)$. However, 
when $\theta \!=\! \pm (M \!-\! 1/2) \pi/3$, $\theta$ can also equal $\pi/2$, but on this occasion $\theta=\pi/2$ applies to the upper-signed 
version of Eq.\ (\ref{sixtyeighta}) with $M \!=\! 2$ in the domain of convergence immediately below it. Furthermore, $S_{MB}(M,z)$ is determined by
putting $M\!=\! 2$ and replacing $z$ by $z^3$ in Eq.\ (\ref{sixtyeightc}). Hence for $M \!=\! 2$ and $\theta \!=\! \pi/2$, we arrive at
\begin{align}
\ln \Gamma \left( z^3 \right) &=F \left( z^3 \right) + TS_N \! \left( z^3 \right)  - 
2z^3 \!\!\!\!\! \!\!\!\!\! \int\limits_{\substack{c-i \infty \\ {\rm Max}[N-1,1/2]<c=\Re\, s<N}}^{c+i \infty} 
\!\!\!\!\! \!\!\!\!\!\!\!\!  ds\; \frac{(1/2\pi z^3)^{2s}\,e^{4i \pi s}}{e^{-i \pi s} -e^{i \pi s}} 
\; \zeta(2s) \, \Gamma(2s-1)
\nonumber\\
& - \;\; 2 \,\pi |z^3| \;\;. 
\label{seventyninea}\end{align}

For $\theta=5 \pi/6$, we have either $M\!=\!3$ when $\theta \!=\!(M-1/2) \pi/3$ or $M\!=\! 2$ when $\theta \!=\! (M+1/2) \pi/3$. In 
the first case we use $S_{MB}(M,z)$ as given by Eq.\ (\ref{sixtyeightc}) with $z \!=\! z^3$ and $M \!=\! 3$. Therefore, $\ln \Gamma(z^3)$ 
via the MB regularization of $S(z)$ becomes
\begin{align}
\ln \Gamma \left( z^3 \right) &=F \left( z^3 \right) + TS_N \! \left( z^3 \right)  - 
2z^3 \!\!\!\!\! \!\!\!\!\! \int\limits_{\substack{c-i \infty \\ {\rm Max}[N-1,1/2]<c=\Re\, s<N}}^{c+i \infty} 
\!\!\!\!\! \!\!\!\!\!\!\!\!  ds\; \frac{(1/2\pi z^3)^{2s}\,e^{6i \pi s}}{e^{-i \pi s} -e^{i \pi s}} 
\; \zeta(2s) \, \Gamma(2s-1)
\nonumber\\
&- \;\; \ln  \Bigl( 1- e^{-2 \pi |z^3|} \Bigr) + 2 \pi |z^3|  \;\;. 
\label{eighty}\end{align} 
For the second case we require Eq.\ (\ref{sixtyeightd}) with $z \!=\! z^3$ and $M \!=\! 2$ to obtain $\ln \Gamma(z^3)$. Then we 
obtain
\begin{align}
\ln \Gamma \left( z^3 \right) &= F \left( z^3 \right) +  TS_N \! \left( z^3 \right)  - 
2z^3 \!\!\!\!\! \!\!\!\!\! \int\limits_{\substack{c-i \infty \\ {\rm Max}[N-1,1/2]<c=\Re\, s<N}}^{c+i \infty} 
\!\!\!\!\! \!\!\!\!\!\!\!\!  ds\; \frac{(1/2\pi z^3)^{2s}\,e^{4i \pi s}}{e^{-i \pi s} -e^{i \pi s}} 
\; \zeta(2s) \, \Gamma(2s-1)
\nonumber\\
&+ \;\; 2 \pi |z^3|  \;\;. 
\label{eightyone}\end{align}

The lower-signed version of Eq.\ (\ref{sixtyeighta}) with $z$ replaced by $z^3$ gives the values of $\ln \Gamma(z^3)$, where the 
domains of convergence for the MB integral cover the negative or lower half of the principal branch of the complex plane. That is, 
the lower-signed version is required for obtaining the values of $\ln \Gamma(z^3)$ when $\theta$ is negative. For $-2\pi/3 < \theta < 0$ 
or $M=1$, we find that $\ln \Gamma(z^3)$ is given by
\begin{align}
\ln \Gamma \left( z^3 \right)& =F \left( z^3 \right) +  TS_N \! \left( z^3 \right)  - 
2z^3 \!\!\!\!\! \!\!\!\!\! \int\limits_{\substack{c-i \infty \\ {\rm Max}[N-1,1/2]<c=\Re\, s<N}}^{c+i \infty} 
\!\!\!\!\! \!\!\!\!\!\!\!\!  ds\; \frac{(1/2\pi z^3)^{2s}\,e^{-2i \pi s}}{e^{-i \pi s} -e^{i \pi s}} 
\; \zeta(2s) \, \Gamma(2s-1) 
\nonumber\\
& + \;\; \ln \left( 1- e^{-2 i \pi z^3} \right)\;\;, 
\label{eightyoneb}\end{align}
while for $-\pi < \theta < -\pi/3$ or $M=2$, one finds that
\begin{align}
\ln \Gamma \left( z^3 \right) &=F \left( z^3 \right) + TS_N \! \left( z^3 \right)  - 
2z^3 \!\!\!\!\! \!\!\!\!\! \int\limits_{\substack{c-i \infty \\ {\rm Max}[N-1,1/2]<c=\Re\, s<N}}^{c+i \infty} 
\!\!\!\!\! \!\!\!\!\!\!\!\!  ds\; \frac{(1/2\pi z^3)^{2s}\,e^{-4i \pi s}}{e^{-i \pi s} -e^{i \pi s}} 
\; \zeta(2s) \, \Gamma(2s-1)
\nonumber\\
& -\;\; \ln \left( -\, e^{-2 i \pi z^3}\right) \;\;. 
\label{eightytwo}\end{align}
Following from our study of the positive values of $\theta$, Eq.\ (\ref{eightyoneb}) is not valid for $\theta=-\pi/2$, while Eq.\  
(\ref{eightytwo}) is not valid for $\theta=-\pi/2 $ and $\theta =-5\pi/6$. 

When $M \!=\! 1$ in $\theta \!=\! -(M+1/2) \pi/3$, we have $\theta \!=\! -\pi/2$. In this case $S_{MB}(M,z)$ is given by Eq.\ 
(\ref{sixtyeightd}) with $z$ replaced by $z^3$ and $M \!=\!1$. Hence we find that
\begin{align}
\ln \Gamma \left( z^3 \right)\left|_{{\rm arg}\,z =-\pi/2}  \right. & =F \left( i |z|^3 \right) +  TS_N \! \left(i |z|^3 \right)  - 
2i |z|^3 \!\!\!\!\! \!\!\!\!\! \int\limits_{\substack{c-i \infty \\ {\rm Max}[N-1,1/2]<c=\Re\, s<N}}^{c+i \infty} 
\!\!\!\!\! \!\!\!\!\!\!\!\!  ds\; \frac{(1/2\pi z^3)^{2s}\,e^{-2i \pi s}}{e^{-i \pi s} -e^{i \pi s}} 
\nonumber\\
& \times \;\; \zeta(2s) \, \Gamma(2s-1)  -  2\pi \, |z|^3 - \ln \left( 1- e^{-2 \pi |z|^3} \right)\;\;. 
\label{eightytwob}\end{align}
The above result represents the complex conjugate of its analog given by Eq.\ (\ref{seventysevenb}). Note that the terms for $S_{MB}(M,z)$ 
or rather $S_{MB}(1,z^3)$ in Eq.\ (\ref{eightytwob}) are identical to the corresponding terms in Eq.\ (\ref{seventysevenb}) since the 
regularized value of the logarithmic series is only real for $\theta \!=\! -\pi/2$ as described in the proof to Thm.\ 4.1.
 
We also have $\theta=-\pi/2$ for $M=2$ in $\theta=-(M-1/2)\pi/3$. In this instance $S_{MB}(M,z)$ is given by Eq.\ (\ref{sixtyeightc})
with $z$ replaced by $z^3$ and $M \!=\!2$. Therefore, we have
\begin{align}
\ln \Gamma \left( z^3 \right) \left|_{{\rm arg} z= -\pi/2} \right. & =F \left( i |z|^3 \right) +  TS_N \! \left( i |z|^3 \right)  - 
2i |z|^3 \!\!\!\!\! \!\!\!\!\! \int\limits_{\substack{c-i \infty \\ {\rm Max}[N-1,1/2]<c=\Re\, s<N}}^{c+i \infty} 
\!\!\!\!\! \!\!\!\!\!\!\!\!  ds\; \frac{(1/2\pi z^3)^{2s}\,e^{-4i \pi s}}{e^{-i \pi s} -e^{i \pi s}} 
\nonumber\\
& \times \;\; \zeta(2s) \, \Gamma(2s-1) - 2 \pi |z|^3 \;\;. 
\label{eightytwoc}\end{align}

The lower-signed version of Eq.\ (\ref{sixtyeighta}) with $M \!=\!2$ also applies to $\theta=-5\pi/6$, but in this case $S_{MB}(2,z^3)$ is
given by Eq.\ (\ref{sixtyeightd}) since $\theta=-(M+1/2)\pi/3$. Then we arrive at
\begin{align}
\ln \Gamma \left( z^3 \right) \left|_{{\rm arg}\,z=-5\pi/6}  \right. & =F \left( -i |z|^3 \right) +  TS_N \! \left( -i|z|^3 \right)  + 
2i |z|^3 \!\!\!\!\! \!\!\!\!\! \int\limits_{\substack{c-i \infty \\ {\rm Max}[N-1,1/2]<c=\Re\, s<N}}^{c+i \infty} 
\!\!\!\!\! \!\!\!\!\!\!\!\!  ds\; \frac{(1/2\pi z^3)^{2s}\,e^{-2i \pi s}}{e^{-i \pi s} -e^{i \pi s}}
\nonumber\\
& \times \;\;   \zeta(2s) \, \Gamma(2s-1) +  2 \pi \, |z|^3 \;\;. 
\label{eightytwod}\end{align}

For $-4 \pi/3 \!<\! \theta \!<\! -2\pi/3$, we can use the $M \!=\!3$ lower-signed version of Eq.\ (\ref{sixtyeighta}) with $z$
replaced by $z^3$ to obtain values of $\ln \Gamma(z^3)$. This yields 
\begin{align}
\ln \Gamma \left( z^3 \right) &= F \left( z^3 \right) +  TS_N \! \left( z^3 \right)  - 
2z^3 \!\!\!\!\! \!\!\!\!\! \int\limits_{\substack{c-i \infty \\ {\rm Max}[N-1,1/2]<c=\Re\, s<N}}^{c+i \infty} 
\!\!\!\!\! \!\!\!\!\!\!\!\!  ds\; \frac{(1/2\pi z^3)^{2s}\,e^{-6i \pi s}}{e^{-i \pi s} -e^{i \pi s}} 
\; \zeta(2s) \, \Gamma(2s-1)
\nonumber\\
&  + \;\; \ln \left( 1- e^{-2 i \pi z^3} \right) - \ln \left( -\, e^{-2 i \pi z^3}\right) \;\;.
\label{eightythree}\end{align} 
The above result, however, is not valid for either $\theta \!=\! -5 \pi/6$ or $\theta \!=\! -7 \pi/6$. For the former case, where $M \!=\! 3$ 
in $\theta \!=\! -(M \!-\! 1/2) \pi/3$, $S_{MB}(M,z)$ is obtained by putting $M \!=\! 3$ and $z \!=\! z^3$ in Eq.\ (\ref{sixtyeightc}). Thus,
we find that 
\begin{align}
\ln \Gamma \left( z^3 \right) &= F \left( z^3 \right) +  TS_N \! \left( z^3 \right)  - 
2z^3 \!\!\!\!\! \!\!\!\!\! \int\limits_{\substack{c-i \infty \\ {\rm Max}[N-1,1/2]<c=\Re\, s<N}}^{c+i \infty} 
\!\!\!\!\! \!\!\!\!\!\!\!\!  ds\; \frac{(1/2\pi z^3)^{2s}\,e^{-6i \pi s}}{e^{-i \pi s} -e^{i \pi s}} 
\; \zeta(2s) \, \Gamma(2s-1)
\nonumber\\
&  - \;\; \ln \left( 1- e^{-2  \pi |z^3|} \right) + 2 \pi |z^3|\;\;.
\label{eightythreea}\end{align}

If we compare the above results with the Borel-summed asymptotic forms given by Eq.\ (\ref{fiftyeight}), then we see that the first
two terms, viz.\ the Stirling approximation or $F(z)$ and the truncated sum $TS_N(z)$, are basically the same, but the remainder
is completely different in that it is now one MB integral rather than an infinite convergent sum of integrals. This is primarily
due to the introduction of the zeta function during MB regularization. In addition, although the logarithmic terms in Eqs.\ 
(\ref{seventyseven}) and (\ref{eightyone}) are similar to those in Eq.\ (\ref{fiftyeight}), they do not appear as discontinuous
quantities with Stokes multipliers. Therefore, we see that Stokes lines are fictitious. They do not imply that the original 
function is necessarily discontinuous along them, only that the Borel-summed asymptotic forms are. For example, although there are 
Stokes lines at $\theta \!=\! \pm \pi/6$, $\ln \Gamma(z^3)$ is not discontinuous along them.  
 
In particular, from Eqs.\ (\ref{seventysix})-(\ref{eightythreea}) we have seen that $\ln \Gamma(z^3)$ is discontinuous whenever 
$\theta= \pm(l+1/2)\pi/3$ for $l$, a positive integer. These values of $\theta$ also represent Stokes lines when $S(z)$ is Borel-summed. 
Consequently, we need to carry out separate numerical investigations: the first will be aimed at showing the agreement between the 
MB-regularized asymptotic forms for $\ln \Gamma(z^3)$ and their Borel-summed counterparts, while the second will deal with the 
the behaviour of $\ln \Gamma(z^3)$ specifically at the Stokes lines/rays. In the first investigation we shall begin by describing how 
to evaluate $\ln \Gamma(z^3)$ via the MB-regularized asymptotic forms. Then we shall be able to compare the results with the Borel-summed 
asymptotic forms with $z$ replaced by $z^3$ in the results of Sec.\ 3. We shall observe that although both MB-regularized 
asymptotic forms are defined at each Stokes line, they give the incorrect values of $\ln \Gamma(z^3)$, which differ by 
discontinuous jumps of $2\pi i$. The second study is aimed at obtaining the correct values of $\ln \Gamma(z^3)$ via the 
Borel-summed and MB-regularized asymptotic forms by the application of the Zwaan-Dingle principle at the Stokes lines.
 
One consequence of MB regularization is that we do not have to truncate the remainder as we did in our previous 
numerical study in Sec.\ 3. That is, we are now in a position to obtain an exact result for the remainder within the limitations
of our computing system, not an approximation as a result of truncating to a large value such as $10^5$. Moreover, we expect to
obtain the regularized value far more quickly than the times reported in Sec.\ 3. Consequently, we can consider even smaller 
values of $|z|$ than we did in Sec.\ 3. 

Since there are no Stokes lines of discontinuity in the preceding results, there are two MB-regularized asymptotic forms that yield 
the values of $\ln \Gamma(z^3)$ for all values of $\theta$ or $\arg\, z$, except when $\theta = k \pi/3$ and $k$ is an integer. This 
means that we can check values of the two different asymptotic forms for the regularized value of $\ln \Gamma(z^3)$, which was not 
possible with the Borel-summed results in Sec.\ 3. As a result of the preceding discussion we, therefore, see that MB regularization 
represents an important technique in (hyper)asymptotics.  

The fourth program presented in the appendix is the Mathematica module called MBloggam. Basically, this code uses the MB-regularized
asymptotic forms given above to evaluate $\ln \Gamma(z^3)$ over the entire principal branch for $z$. First, outside the module there 
appears a statement for Intgrd, in which the integrand of the MB integrals is expressed in the following form:
\begin{eqnarray}
I\left(|z|,\theta,s,M \right)= \left( 2\pi |z|^3 \right)^{-2s} \,\zeta(2s) \, \Gamma(2s-1)\, \left( \frac{e^{i(M\pi-3 \theta) s}}
{e^{-i \pi s} -e^{i \pi s}}\right)  \;\;.
\label{eightyfour}\end{eqnarray}
In this statement $z^3$ in the MB integral has been replaced by $|z|^3 \exp(3i \theta)$ so that the exponential can be combined with
the phase factor associated with the domains of convergence, viz.\ $\exp(M i \pi s)$. This is necessary to ensure that 
the integrand does not diverge at any stage when the combined exponential factor is divided by the exponential term in the 
denominator. For example, as $s \to \infty$, $\exp(M i \pi s)$ diverges either when $\Im \, s < 0$ and $M$ is positive or 
when $\Im\, s >0$ and $M$ is negative. By combining it with $\exp(-3i \theta s)$ and the exponential term in the denominator, 
we avoid this divergence when numerically integrating the MB integral in the specific asymptotic forms between Eqs.\ (\ref{seventysix})
and (\ref{eightythree}), not pertaining to the Stokes lines.

The next three statements outside the module are the expression for: (1) the $c_k$, (2) the Stirling approximation $F(z)$ and (3) the
summand in the truncated sum $TS_N(z)$. These statements are identical to those in the modules when $\ln \Gamma(z)$ was evaluated via the 
Borel-summed asymptotic forms in Sec.\ 3. 

Inside the module the variables s and s1 denote the variables of integration for the upper and lower halves of the MB integral
respectively. Next the variable zcube representing the complex value of $z^3$ is calculated followed by the Stirling approximation, 
which is denoted by e0. Then the truncated sum $TS_N(z^3)$ is evaluated as e1. Because two different asymptotic forms can be used to 
evaluate the regularized value of $\ln \Gamma(z^3)$, two different values of $M$ need to be determined. These are symbolized by M1 and 
M2 and are both initialized to zero. Their true values are determined in the first Which statement. If M2 remains zero, then it means 
that only one of the above equations can be used to evaluate the regularized value. When this occurs, e.g. for $\theta = \pm \pi/3$, 
the module prints out that only one value of $M$ applies. In order to evaluate the entire MB integral in the above equations, two 
separate calls to the NIntegrate routine are made: one corresponding to positive imaginary values or s and the other to negative 
imaginary values or s1. When $M$=M1, these values are represented by e2 and e3, while for $M$=M2, they are represented by 
e5 and e6. The combined values are then multiplied by $-2 z^3$ or rather -2 zcube, to give the total contribution of the MB 
integral to $\ln \Gamma(z^3)$. The next Which statement uses the value of M1 to determine the appropriate logarithmic term that 
needs to be included to obtain the combined value of $\ln \Gamma(z^3)$, while if M2 is non-zero, then the logarithmic term is
evaluated in the last Which statement for the second asymptotic form of $\ln \Gamma(z^3)$. Finally, if $-\pi/3 < \theta \leq \pi/3$, 
then the value of $\ln \Gamma(z^3)$ is evaluated by Mathematica's intrinsic routine LogGamma[z] and printed out with the other results.  

In the previous numerical study we considered both a ``large" and intermediate value of $|z|$. The large value or $|z| \!=\! 3$ 
was chosen primarily because it was deemed sufficiently large in order to observe whether a smoothing of the Stokes phenomenon 
as postulated by Berry \cite{ber89} and Olver \cite{olv90} occurs or not. The intermediate value of $|z|$ or $|z| = 1/10$ was chosen 
because it represented a value where standard Poincar$\acute{\rm e}$ asymptotics breaks down. We shall choose the latter value here 
again for the same reason, but on this occasion the value becomes very small because the variable in the MB integrals in Eqs.\ 
(\ref{seventysix})-(\ref{eightythree}) is actually $2\pi z^3$. This means that we are effectively considering $2\pi \times 10^{-3}$, 
instead of $2\pi/10$ as the magnitude of the variable in Sec.\ 3. Such a small value would be deemed impossible under standard 
Poincar${\acute {\rm e}}$ asymptotics or even by employing the hyperasymptotic methods in Refs.\ \cite{ber90}-\cite{ber91a} and 
\cite{par01}.

There are, however, some issues arising from selecting such a small value. The first is that we expect that both the truncated
series, $TS_N(z)$, and the MB integral in Eqs.\ (\ref{seventysix})-(\ref{eightythree}) to begin to diverge very rapidly for relatively
small values of the truncation parameter such as $N \!=\! 4$. Consequently, there will be a great cancellation of decimal places
when adding the $S_N(z)$ to the MB integral, which in turn means that even though the accuracy and precision goals have been 
set to 30 in all the calls to the NIntegrate routine in MBloggam, we may not necessarily obtain a final value that is accurate
to this level. Although WorkingPrecision has been set higher to 80 to allow for this possibility, there is still no guarantee that
the final value of $\ln \Gamma(z^3)$ will indeed be accurate to 30 decimal places. To overcome this problem, one needs to specify
larger values of AccuracyGoal, PrecisionGoal and WorkingPrecision, but this comes at the expense of computing time.

\begin{table}
\centering
{\small
\begin{tabular}{|c|c|c|c|} \hline
 $\theta$ & $N$ & Quantity  &  Value \\ \hline
& &$ F(z^3)$ &                     4.3666691849467394839681993920 + 0.39773534871318634708519397906$\, i$ \\
& & $TS_3(z^3)$ &                   1.964244428861064224518267$ \times 10^{6}$ - 1.9641265777308664665975342$ \times 10^{6}\, i$ \\
$-\pi/12$ & 3 & MB Int (M1=0) &    -1.964241888183123016922580$ \times 10^{6}$ + 1.964126965801012078012431$ \times 10^{6}\, i$\\
& & $S_{MB}(0,z^3)$  &  0 \\ 
& & Total via M1 &                 6.9073471261543351713515623993 + 0.7858054943246012439632804975$\, i$ \\
& & MB Int (M2=-1) &              -1.964246960282777058252170$ \times 10^{6}$ + 1.964127748979378940448405$ \times 10^{6}\, i$ \\   
& & $S_{MB}(1,z^3)$  &              5.0720996540413295899999675136 - 0.783178366862435974379475023745$\, i $ \\
& & Total via M2    &              6.9073471261543351713515623992 + 0.78580549432460124396328049761$\, i $ \\
& & LogGamma[zcube] &              6.9073471261543351713515623992 + 0.78580549432460124396328049761$\, i$ \\ \hline
& &$ F(z^3)$ &                     4.3790700299188033385944250366 - 1.38001259938612058620816944744$\, i$ \\
& & $TS_4(z^3)$ &                   3.03718072746107697324584$ \times 10^{11}$ - 7.332351579151624817641834$ \times 10^{11}\, i$ \\
$7\pi/24$ & 4 & MB Int (M1=0) &   -3.03718072743578478216056$ \times 10^{11}$ +  7.33235157913793379318864$ \times 10^{11}\, i$ \\
& & $S_{MB}(0,z^3)$ &  0 \\ 
& & Total via M1 &                 6.9082891384461367831946353384 - 2.749115044705401162838355615704$\, i$ \\   
& & MB Int (M2=1) &               -3.03718072748649559827233$ \times 10^{11}$ +  7.33235157914968575273952$ \times10^{11}\, i$ \\    
& & $S_{MB}(1,z^3)$  &              5.0710816111765935339415417904 - 1.175195955088607468412668235518$\, i$ \\
& & Total via M2    &              6.9082891384461367827408810781 - 2.749115044705401165409263040026$\, i$ \\
& & LogGamma[zcube] &              6.9082891384461367827408810777 - 2.749115044705401165409263038133$\, i$ \\ \hline
& &$ F(z^3)$ &                     4.3807239279747234048593708927 - 1.57393791944848641246978433502$\, i$ \\
& & $TS_2(z^3)$ &                  -83.333333333333333333333333333 + 0 $\,i\;\;\;\;\;\;\;\;\;\;\;\;\;\;\;\;\;\;\;\;\;\;\;\;\;\;\;\;\;\;\;\;\;\;\;\;\;\;\;\;\;\;\;\;\;\;\;\;\;$ \\
$\pi/3$ & 2 & MB Int (M1=1) &      80.791062865366238781576251609 + 0 $\,i\;\;\;\;\;\;\;\;\;\;\;\;\;\;\;\;\;\;\;\;\;\;\;\;\;\;\;\;\;\;\;\;\;\;\;\;\;\;\;\;\;\;\;\;\;\;\;\;\;\;$ \\ 
& & Log. Term (M1=1)   &           5.0698798575073995786757215377 - 1.56765473414130682599285904825$\, i$ \\ 
& & Total via M1 &                 6.9083333175150284317780107065 - 3.141592653589793238462643383279$\, i$ \\
& & LogGamma[zcube] &              6.9083333175150284317780107065 - 3.141592653589793238462643383279$\, i$ \\ \hline
& &$ F(z^3)$ &                     4.3671839976260822611773860371 + 0.45442183940812929747019906926$\, i$ \\
& & $TS_5(z^3)$ &                  -5.95238271839508333182790$ \times 10^{17}$ - 7.737535164228715974701668$ \times 10^{11}\, i$ \\
$4\pi/7$ & 5 & MB Int (M1=1) &     5.95238271839508330650668$ \times 10^{17}$ + 7.737535164239864652939712$ \times 10^{11}\, i$ \\
& & $S_{MB}(1,z^3)$   &           5.0674216504983723676332001993 + 2.46643387314754360954113189663$\, i$ \\
& & Total via M1 &                 6.9024828161628933968353749191 + 4.03572353636012752353352062801$\, i$  \\   
& & MB Int (M2=2) &                5.95238271839508335723002$ \times 10^{17}$ + 7.737535164233152240154912$ \times 10^{11}\, i$ \\ 
& & $S_{MB}^{U}(2,z^3)$  &              5.0710816111765935339415417904 - 1.175195955088607468412668235518$\, i$ \\
& & Total via M2    &              6.9024828161628933880358773203 + 4.03572353636012787577347290435$\, i$ \\ \hline
& &$ F(z^3)$ &                     4.3749562509709981184827498273 - 1.05509306570630337542065838646$\, i$ \\
& & $TS_6(z^3)$ &                   8.41751139369492541714725$ \times 10^{23}$ + 5.154919990982005385545807$ \times 10^{17}\, i$ \\
$8\pi/9$ & 6 & MB Int (M1=2) &    -8.41751139369492541714722$ \times 10^{23}$ - 5.154919990982005395943833$ \times 10^{17}\, i$ \\
& & $S_{MB}^{U}(2,z^3)$   &           0.0054413980927026535517822347 - 3.13845106093620344522418073989$\, i$ \\
& & Total via M1 &                 6.9134848732085864689307216827 - 5.23334675905750858730776717724$\, i $ \\
& & MB Int (M2=3) &               -8.41751139369492541714727$ \times 10^{23}$ - 5.154919990982005390723539$ \times 10^{17}\, i$ \\
& & $S_{MB}^{U}(3,z^3)$  &              5.0780394872445415442384564511 - 3.660480464761928202988399777482$\, i$ \\
& & Total via M2    &              6.9134848732729541247531300647 - 5.233346759035781054140026255124$\, i$ \\ \hline
\end{tabular}
}
\normalsize
\caption{Determination of $\ln\Gamma(z^3)$ via the MB-regularized forms with $|z| \!=\! 1/10$ and various values of 
$\theta$ and $N$}
\label{tab8}
\end{table}

Table\ \ref{tab8} presents a very small sample of results obtained by running MBloggam on a Sony VAIO laptop with 2 GB RAM and
Mathematica 9.0 for various values of the truncation parameter and argument of $z$. In general, the calculations took between
150 and 350 CPU seconds, although the results for $\theta \!=\! \pm 11\pi/17$ and $N \!=\! 2$ took about 545 seconds to execute. 
This is considerably faster than the corresponding Borel-summed results in Table\ \ref{tab4}, which not only had to be truncated, 
but employed a larger value of $|z|$. It should also be mentioned that hundreds of results were obtained by creating a script file
out of the fourth program in the appendix and running the resulting file on the SunFire alpha server mentioned in Sec.\ 3. 
Although each calculation was slower than the laptop, the advantage with the latter system is that a far greater number of 
calculations can be performed simultaneously. Nevertheless, the results in the table are indicative of those obtained from the 
SunFire server.

The table presents the results obtained from five calculations, four with positive values of $\theta$ and one with a negative value.
In addition, different values of the truncation parameter have been selected. The first line or row of each calculation represents
the Stirling approximation or $F(z^3)$. Then comes the value of the truncated series, which is denoted by $TS_N(z^3)$. The remainder
denoted by MB Int appears next. As mentioned earlier, because the domains of convergence for the MB integrals overlap one another, two 
different MB integrals can be computed in the evaluation of the remainder of $S_N(z^3)$. The first MB integral is represented 
by M1, while the second is represented by M2. The second MB integral is not evaluated if the value of M2 is zero, which occurs
when $\theta \!=\!l \pi/3$ and $l$ is an integer. When printing out the values of the MB integrals, the values of M1 and M2 are also 
specified by the module and thus, appear in the adjacent column. The values of $N$ and $\theta$ are presented on the row with the value 
of the first MB integral for each calculation. Appearing on the following row after each value of an MB integral is the value of 
$S_{MB}(M,z^3)$. The superscripts U and L denote whether the upper-signed or lower-signed version in Eq.\ (\ref{sixtyeightb}) has
been used in obtaining the displayed value. For example, for M1=0, this term vanishes and hence there is no superscript, while for M1 
or M2 equal to 2, we have $S_{MB}^{U}(1,z^3)=\ln(-\exp(-2i \pi z^3))$ and $S_{MB}^{L}(1,z^3)=-\ln(-\exp(-2i \pi z^3))$.     

The first calculation in Table\ \ref{tab8} presents the results for $\theta \!=\! -\pi/12$ and $N \!=\! 3$. For this situation we 
expect Eqs.\ (\ref{seventysix}) and (\ref{eightyone}) corresponding to M1 \!=\! 0 and M2 \!=\! -1, respectively, to yield the value of 
$\ln \Gamma(\exp(-i \pi/4)/1000)$. The first line gives the Stirling approximation to the value, which we see is substantial, but by no 
means accurate when compared with the actual value evaluated by Mathematica's LogGamma routine presented at the bottom or ninth 
row of the calculation. The second row of the calculation displays the value of the truncated series $TS_3(\exp(-i\pi/4)/1000)$, which is 
already of the order of $10^6$. Therefore, to obtain the actual value of $\ln \Gamma(\exp(-i \pi/4)/1000)$, we require the cancellation of 
at least six decimal figures from this result by the remainder or MB integral. The value of the MB integral determined from Eq.\ 
(\ref{seventysix}) appears on the next row, which is not only of the order of $10^6$, but does in fact cancel the first six decimal 
places of the value of the truncated series. The value of $S_{MB}(0,\exp(-i\pi/4)/1000)$, which is zero in this instance, appears on 
the fourth row of the calculation, while the sum of the Stirling approximation, the truncated sum, the MB integral and 
$S_{MB}(0,\exp(-i\pi/4)/1000)$ appears in the fifth row represented by `Total via M1'. As we can see, the total value agrees with the 
actual value of $\ln \Gamma(\exp(-i \pi/4)/1000)$ to 30 decimal places, well within the accuracy and precision limits specified 
in each call to the NIntegrate routine in the program, especially when it is borne in in mind that six decimal places have been 
lost by summing the MB integral with the truncated series. In this case we have been saved by the fact that the working precision 
was set to a very high value. The sixth row of the first calculation displays the value of the MB integral in Eq.\ (\ref{eightyone}). 
As expected, it agrees with the first six decimal figures of the values for both the truncated sum and the MB integral in Eq.\ 
(\ref{seventysix}) given in the third row. Unlike $S_{MB}(0,z^3)$, however, $S_{MB}(1,z^3)$ is non-vanishing and is presented on 
the seventh row. There it can be seen that the real and imaginary parts of this value are much greater in magnitude than those for 
the Stirling approximation. If this value is summed only with the Stirling approximation, then the resulting value deviates from 
the actual of $\ln \Gamma(\exp(-i \pi/4)/1000)$ far more than either value on its own. However, when the latter value is summed 
with the values of the truncated sum and the MB-regularized remainder, it yields the value of $\ln \Gamma(\exp(-i \pi/4)/1000)$ 
to the 29 decimal places listed in the table despite the fact the first six decimal figures for the truncated sum and MB integral 
cancel each other. 
  
The second calculation in Table\ \ref{tab8} presents the results obtained for $\theta \!=\! 7 \pi/24$ and $N \!=\! 4$. In this instance 
Eqs.\ (\ref{seventysix}) and (\ref{seventyseven}), which means in turn that M1=0 and M2=1 respectively, will be the only asymptotic forms
that are valid for obtaining the MB-regularized value of $\ln \Gamma(\exp(7i \pi/8)/1000)$. Because the overall argument is $7\pi/8$,
the value of $\ln \Gamma(\exp(7i \pi/8)/1000)$ can still be evaluated via the LogGamma routine in Mathematica. Therefore, we can 
check the results obtained by summing all the contributions in Eqs.\ (\ref{seventysix}) and (\ref{seventyseven}) with that evaluated
by Mathematica. As in the first calculation the first row gives the value for the Stirling approximation or $F(z^3)$ as is the case 
for all other calculations in the table. Surprisingly, we see that there is very little variation in the real part of $F(z^3)$,
but the imaginary does display considerable variation over the entire table. Although the truncation parameter for the second calculation 
is one larger than in the first calculation, we see that the values of the truncated series and MB integrals are of the order of $10^{11}$, 
which means that now 11 decimal figures will need to be cancelled in order to arrive at $\ln \Gamma(\exp(7i \pi/8)/1000)$. This places 
pressure on all the calls to the NIntegrate routine to yield accurate values for both the total calculations. If one compares the total 
calculation via M2 with the final or LogGamma value, then one sees that both values agree almost to within the accuracy and precision 
goals set in the MBloggam. However, if one checks the total via M1 with the LogGamma value, then they only agree to about 17 decimal 
places. More interesting is the fact that when WorkingPrecision is set equal 100 rather than 60, the M2 total becomes more accurate by 
agreeing to well past 30 decimal places with the LogGamma value, but the M1 remains as accurate as the result in Table\ \ref{tab8}. This 
means that AccuracyGoal and PrecisionGoal also need to be extended, not just WorkingPrecision. When AccuracyGoal and PrecisionGoal are 
set to 40 with WorkingPrecision at 100 in all the calls to NIntegrate, it is found the total via M1 agrees with the value of 
$\ln \Gamma(\exp(7i \pi/8)/1000)$ to more than thirty decimal figures, but now the calculation takes 635 CPU seconds compared with 
the 301 CPU for the results in Table\ \ref{tab8} and 490 CPU seconds when only WorkingPrecision is altered to 100. If $N \!=\! 2$, 
then both the totals via M1 and M2 are accurate within the original accuracy and precision goals set in the module. In general, it is 
found that one of the calculations is significantly less accurate than the other when $\theta$ is close to one of the limits of the 
domain of convergence of its MB integral. For example, in the second calculation $\theta=7 \pi/24$, which is quite close to the upper 
limit of $\pi/3$ in the domain of convergence for Eq.\ (\ref{seventysix}). Thus, the MB-regularized remainder obtained from the form 
that is closer to the centre of its domain of convergence is frequently the more reliable result.  

The third calculation in Table\ \ref{tab8} presents the results for $\theta \!=\! \pi/3$. As stated previously, only one MB-regularized form, 
viz.\ Eq.\ (\ref{seventyseven}), can be used to evaluate $\ln \Gamma(z^3)$ since this value of $\theta$ is outside the domain of convergence
for Eq.\ (\ref{seventysix}). Because there is one valid MB-regularized asymptotic form in this calculation, only six rows appear in Table\ \ref{tab8},
not nine as in both the previous calculations. In this case MBloggam prints out that M2 is zero, while M1 equals unity, indicating that the 
MB-regularized value has been evaluated via Eq.\ (\ref{seventyseven}). Since the truncation parameter was set equal to 2, it results in the 
lowest magnitudes for both the truncated series and the MB integral of all the calculations presented in the table. Consequently, there is 
only one decimal figure that is cancelled when the value of truncated sum is added to the value of the MB integral. As expected, the total 
agrees with the LogGamma value well within the accuracy and precision goals set in the module. Furthermore, because this calculation involved 
much less computation than the other calculations, it took only 63 CPU seconds to execute. 

If we look closely at the imaginary part of the $\theta=\pi/3$ calculation in Table\ \ref{tab8}, then we see that $\Im (\ln \Gamma(z^3))=
-\pi$. This can be proved by noting that the asymptotic series $S(z^3)$ is composed of real terms when $\theta = k\pi/3$ and
$k$, an integer. Hence the imaginary part of $TS_N(z^3)$ vanishes for all these values. In addition, the imaginary
part of the regularized remainder of the series as given by the MB integral can be shown to vanish by splitting the integral
into two integrals and making the substitutions, $s=c +it$ for the integral in the upper half of the complex plane and $s=c-it$ for
the integral in the lower half of the complex plane. Then it is found that all the factors in one of the integrals are complex
conjugates of all the factors in the other integral when $\theta =\pi/3$. By writing all these factors as sums of a real and imaginary
part, it is found that when they are multiplied out or expanded, the imaginary terms cancel and we are left with a real integral. 
Therefore, the imaginary part of $\ln \Gamma(z^3)$ for $\theta= \pi/3$ becomes
\begin{eqnarray}
\Im \ln \Gamma \Bigl( |z|^3 \exp(i\pi)\Bigr)= \Im F\Bigl(|z|^3 e^{i\pi} \Bigr) - \Im \ln \Bigl(1 - e^{-2i \pi |z|^3}\Bigr)\;.
\label{eightyeighta}\end{eqnarray} 
From Eq.\ (\ref{fiftynine}) we arrive at
\begin{eqnarray}
\Im F \Bigl( |z|^3 e^{i \pi} \Bigr)= -\left( |z|^3 + \frac{1}{2} \right) \pi \;,
\label{eightyeightb}\end{eqnarray}
while the second term on the rhs of Eq.\ (\ref{eightyeighta}) can be expressed as
\begin{eqnarray}
\Im \ln \Bigl( 1- e^{-2i \pi |z|^3} \Bigr)= \Im \ln \Bigl( e^{-i\pi |z|^3} \Bigr) + \Im \ln \Bigl(2 i \sin(\pi |z|^3) \Bigr) \;.
\label{eightyeightc}\end{eqnarray}
Introducing the above results into Eq.\ (\ref{eightyeighta}) yields
\begin{eqnarray}
\Im \ln \Gamma \left( z^3 \right) \Big|_{\theta= \pi/3}= -\pi    \;. 
\label{eightyeightd}\end{eqnarray} 

The two remaining calculations in Table\ \ref{tab8} have been carried out for $\theta \!>\! \pi/3$, which means that we can no longer use 
the LogGamma routine as a check on the results. That is, for these calculations we can only compare the totals via the M1 and M2 asymptotic 
forms with each other. Moreover, when $\theta$ equals $2\pi/3$ or $\pi$, there will only be one MB-regularized form that we can use 
to obtain the regularized value. For these situations we require the Borel-summed regularized values as a check, which will be done
later in this section. 

The fifth calculation presents the calculation for $\theta \!=\! 4\pi/7$ and $N \!=\! 5$, the latter representing the highest value of the truncation
parameter presented so far in the table. As a result, the truncated sum is of the order of $10^{17}$, which means that 17 decimal figures need 
to be cancelled before we can obtain the value of $\ln \Gamma(\exp(12i \pi/7)/1000)$. Because $\theta$ is situated in the domains of convergence of
$(0,2\pi/3)$ and $(\pi/3,\pi)$, we have M1=1 and M2=2. For this calculation Eqs.\ (\ref{seventyseven}) and (\ref{seventyeight}) apply, which is 
interesting because $S_{MB}(M,z^3)$ is very different for both values of $M$ in these asymptotic forms, particularly the imaginary parts. As 
expected, the MB integrals for both asymptotic forms yield the 17 decimal figures in their real parts needed to cancel those in the 
real part of truncated sum or $TS_5(\exp(12i \pi/7)/1000)$. The imaginary parts only result in the cancellation of eleven decimal figures. 
As stated previously, the cancellation of a large number of decimal places places pressure on the accuracy of the totals. Here we see
that the real parts of the totals only agree to 17 decimal figures. Surprisingly, although there were less decimal figures involved in
the cancellation of the imaginary parts, the imaginary parts in the totals agree to the same number of decimal figures as the real
parts.

Since $4\pi/7$ is closer to the upper limit of $2\pi/3$ for the domain of convergence of Eq.\ (\ref{seventyseven}), we expect the 
total obtained via M1 in the table to be the less accurate of the two forms. In actual fact, it turns out that the total obtained
via Eq.\ (\ref{seventyseven}) is more accurate than the total via M2 or Eq.\ (\ref{seventyeight}) by a few extra decimal places. 
This is readily seen by extending WorkingPrecision to 100 and AccuracyGoal and PrecisionGoal to 40 in MBloggam. Although the 
time of execution increases from 172 to 411 CPU seconds, both totals agree to 32 decimal places. Another option for improving
the fourth calculation in the table is to reduce the truncation parameter. 

The final calculation in Table\ \ref{tab8} presents the results for $\theta \!=\! 8\pi/9$ and $N \!=\! 6$, the latter representing the
largest value of the truncation parameter in the entire table. For this value of $\theta$, where M1=2 and M2=3, Eqs.\ (\ref{seventyeight}) 
and (\ref{seventynine}) are the valid MB-regularized asymptotic forms. Because of the high value of the truncation parameter, we see that 
the real parts of the truncated sum and MB integrals are of the order of $10^{23}$. However, the respective imaginary parts are six orders 
down as in the previous calculation, similar to the $\theta \!=\! 4\pi/7$ calculation. This was not observed in the first two calculations 
in Table\ \ref{tab8}, where the orders of the real and and imaginary parts of the truncated sum and MB integrals were found to be identical 
for lower values of the truncation parameter. Because the real parts of the truncated sum and MB integrals in the final calculation are 
so large, we expect that the totals via M1 and M2 to be the least accurate of all the calculations in the table, which is indeed the case 
as both totals only agree with each other to ten decimal figures. Nevertheless, all the calculations including those not displayed in the 
table confirm the validity of Eqs.\ (\ref{seventysix}) to (\ref{eightythree}).

\begin{table}
\centering
{\small
\begin{tabular}{|c|c|c|} \hline
$\theta$ & Quantity  &  Value \\ \hline
$19 \pi/40$ &  Total via M1(=1) & 6.9017797138225092740511474835 - 1.3331484570580039616320161702$\, i$ \\
&              Total via M2(=2) & 6.9017797138225092740511474835 - 1.3331484570580039616320161702$\, i $ \\
$\pi/2$     &  Total via M1(=1) & 6.9014712712081946221027015741 - 1.5702191115306805133718396291$\, i$ \\                
            &  Total via M2(=2) & 6.9014712712081946221027015741 + 4.7129661956489059635534471374$\, i$ \\
$21\pi/40$  &  Total via M1(=1) & 6.9015102177011253639269574406 + 4.4758636443386950563445853873$\, i$ \\
            &  Total via M2(=2) & 6.9015102177011253639269574406 + 4.4758636443386950563445853873$\, i$ \\ \hline
$62\pi/75$  &  Total via M1(=2) & 6.9139890062805982640446908483 + 1.6326576822112902043838680277$\, i$ \\
&              Total via M2(=3) & 6.9139890062805982640446908483 + 1.6326576822112902043838680277$\, i $\\
$5\pi/6$    &  Total via M1(=2) & 6.9140376418225537950565521476 + 1.5702191115306805133718396291$\, i$ \\                
            &  Total via M2(=3) & 6.9140376418225537950565521476 + 1.5702191115306805133718396291$\, i$ \\
$63\pi/75$  &  Total via M1(=2) & 6.9140614934733956410110131643 - 4.7754024883371855095758498194$\, i$ \\
            &  Total via M2(=3) & 6.9140614934733956410110131643 - 4.7754024883371855095758498194$\, i$ \\ \hline
$-12\pi/75$ &  Total via M1(=0) & 6.9077182194043652368591902822 + 1.5085404469087662793902283694$\, i $ \\
&              Total via M2(=-1) & 6.9077182194043652368591902822 + 1.5085404469087662793902283694$\, i $ \\
&              $\ln \Gamma(z^3)$ & 6.9077182194043652368591902822 + 1.5085404469087662793902283694$\, i $ \\
$-\pi/6$    &   Total via M1(=0) & 6.9077544565153742085796268609 + 1.5713735420591127250908037541$\, i $ \\  
            &   Total via M2(=-1) & 6.9077544565153742085796268609 + 1.5713735420591127250908037541$\, i$\\
&              $\ln \Gamma(z^3)$ &  6.9077544565153742085796268609 + 1.5713735420591127250908037541$\, i $ \\
$-13\pi/75$ &   Total via M1(=0) &  6.9077907065971626138255125982 + 1.6342043592171290258017534223$\, i $ \\
            &   Total via M2(=-1) & 6.9077907065971626138255125982 + 1.6342043592171290258017534223$\, i $ \\
&              $\ln \Gamma(z^3)$ &  6.9077907065971626138255125982 + 1.6342043592171290258017534223$\, i $ \\ \hline
\end{tabular}
}
\normalsize
\caption{Determination of $\ln\Gamma(z^3)$ via the MB-regularized forms in the vicinity of the lines of discontinuity given by 
$\theta \!=\! -\pi/6$, $\theta \!=\! \pi/2$ and $\theta \!=\! 5\pi/6$ with $|z| \!=\! 1/10$ and $N \!=\! 5$}
\label{tab8a}
\end{table}
 
Now let us examine when we run the code for values of $\theta$ in the vicinity of the Stokes lines. Although we should not run
the code when $\theta$ corresponds directly to a Stokes line, we shall nonetheless do so since the MB integral in the MB-regularized
results is defined on each Stokes line. In fact, the discontinuities at the Stokes lines are due to the regularization of the series 
that emerges when all the contributions from the infinite number of singularities in the Cauchy integrals along these lines are summed. 
Had we been dealing with a finite number terminants as in Ref.\ \cite{kow09}, we would not have obtained the logarithmically divergent 
series appearing in the proof of Thm.\ 4.1 at all and thus, there would have been no need to regularize them.

Table\ \ref{tab8a} presents the results obtained by running MBloggam in the vicinity of the Stokes lines at $\theta \!=\! \pi/2$, $\theta \!=\! 
5\pi/6$ and $\theta \!=\! -\pi/6$ with $|z| \!=\! 1/10$ and the truncation parameter $N$ set equal to 5. When $\theta \!=\! \pi/2$, the code 
evaluates $\ln \Gamma(z^3)$ via Eqs.\ (\ref{seventyseven}) and (\ref{seventyeight}). These are indicated in the table by $M1 \!=\! 1$ and 
$M2 \!=\!2$, respectively. The first two results in the table give the values of Eqs.\ (\ref{seventyseven}) and (\ref{seventyeight}) for 
$\theta \!=\! 19 \pi/40$ or close to the discontinuity at $\theta =\pi/2$. As expected, we see that both forms of $\ln \Gamma(z^3)$ give 
identical values. At $\theta \!=\! \pi/2$, we find that both forms give different results, but only for the imaginary parts. In fact, there is 
a jump discontinuity of $2\pi i$ between the results with the first form giving a value of $-i \pi/2$ and the second giving a value of $3i\pi/2$. 
The reason why that the differences will occur in multiples of $2\pi i$ is because taking the exponential yields $\Gamma(z^3)$, which does 
not possess discontinuities at all for $\theta \!=\! (M \!-\! 1/2) \pi/3$. That is, the discontinuities arise only as a result of taking the 
logarithm of the gamma function. As far as the results in the table for $\theta \!=\! \pi/2$ are concerned, neither is the correct result. We 
need to evaluate the average of the two results, in which case we find that $\Im \ln \Gamma(z^3)\vert_{\theta =\pi/2} \!=\! \pi/2$.   

The next six results represent the values of $\ln \Gamma(z^3)$ when $\theta$ is even closer to $5 \pi/6$ than the previous case of 
$\theta \!=\! \pi/2$, viz. $\pm \pi/150$ above and below $5 \pi/6$. In this instance Eqs.\ (\ref{seventyeight}) ($M1=2$) and (\ref{seventynine}) 
($M2=3$) are used to evaluate the values of $\ln \Gamma(z^3)$. Once again, we observe that the values obtained via both asymptotic forms give 
identical results to one another before and after $\theta \!=\! 5 \pi/6$. However, for $\theta \!=\! 5\pi/6$, they give identical values for both 
the real and imaginary parts. In fact, although the imaginary parts have the same value of $i\pi/2$, we can see that this is not correct 
because the value of $\ln \Gamma(z^3)$ experiences a jump discontinuity of almost $-2 \pi i$. In this instance it appears that Mathematica 
has chosen the wrong value of the logarithmic terms in Eqs.\ (\ref{seventyeight}) and (\ref{seventynine}), which is discussed on p.\ 564 
of Ref.\ \cite{wol92}. If we acknowledge that there is a jump of $-2 \pi$, then we find that $\Im \ln \Gamma(z^3)\vert_{\theta=5 \pi/6} 
\!=\! -\pi/2$.
   
The final set of results in Table\ \ref{tab8a} are the values of $\ln \Gamma(z^3)$ obtained when $\theta$ is very close to $-\pi/6$, 
actually $\pi/150$ above and below $-\pi/6$, as in the previous set of results. Because $|\theta|< \pi/3$, we are also able to evaluate 
$\ln \Gamma(z^3)$ using LogGamma[z] in Mathematica. Therefore, there are more results for this calculation than for the other calculations
in the table. For this calculation MBloggam uses Eqs.\ (\ref{seventysix}) ($M1=0$) and (\ref{eightyoneb}) ($M2=-1$) to evaluate 
$\ln \Gamma(z^3)$. As in the previous calculations we find that the two versions of $\ln \Gamma(z^3)$ give identical values above and 
below the Stokes line at $\theta=-\pi/6$. Moreover, these results agree with the values obtained via LogGamma[z]. The interesting
point about this calculation, however, is that all three values at the Stokes line $\theta=-\pi/6$ also agree with each other.
This is to be expected for Eq.\ (\ref{seventysix}) since there is no extra logarithmic term. For the second term it turns out
that the logarithmic term is  real at $\theta \!=\!\pm \pi/6$. Therefore, there is no need to restrict the sum over the contributions
or the sum on the lhs of Equivalence (\ref{seventyone}) to real values as for $\theta \!=\! \pm \pi/6$. Thus, we conclude that there is no
discontinuity in $\ln \Gamma(z^3)$ when $\theta \!=\! \pm \pi/6$. Hence the Stokes discontinuities in the Borel-summed regularized values
at these Stokes lines are fictitious.  

There is still one more task to carry out before we can complete this numerical study. Whilst we have verified the MB-regularized asymptotic
forms for $\ln \Gamma(z^3)$, we need to do the same for the Borel-summed asymptotic forms given in Thm.\ 2.1. It was not possible to check 
these results because their regions of validity did not overlap one another as the domains of convergence in the MB-regularized asymptotic
forms do. Now we can use the MB-regularized asymptotic forms to check the validity of the Borel-summed counterparts since the regularized 
value of an asymptotic series is by definition independent of the method used to derive it. Such an investigation can also confirm whether the
MB-regularized asymptotic forms for $\theta \!=\! k \pi/3$, where $k$ is equal to $\pm 1$, $\pm 2$ and $3$, are correct, because we have seen 
that only one MB-regularized asymptotic form applies for these values of $\theta$.              

Before we can conduct the second part of the numerical study, we need to replace $z$ by $z^3$ in Thm.\ 2.1. By adopting the same notation used
to obtain Eqs.\ (\ref{seventysix})-(\ref{eightythree}), we find that the Borel-summed regularized values for $\ln \Gamma(z^3)$ become
\begin{eqnarray}
\ln \Gamma \left( z^3 \right) =F \left( z^3 \right) + TS_N \! \left( z^3 \right) + R_N^{\pm} \!\left( z^3 \right) +
SD^{\pm}_M \!\left( z^3 \right) \;\;,
\label{eightyfive}\end{eqnarray}
where, as before, $TS_N(z^3)$ is the truncated form of $S_N(z^3)$ at $N$, 
\begin{eqnarray}
R_N^{+} \!\left(z^3 \right) = \frac{2(-1)^{N+1}\,z^3}{(2 \pi z^3)^{2N-2}} \int_0^{\infty} dy\, y^{2N-2}\,e^{-y} \sum_{n=1}^{\infty}
\frac{1}{n^{2N-2} \left( y^2 + (2 n \pi z^3)^2 \right)} \;\;,
\label{eightysix}\end{eqnarray}
\begin{eqnarray}
R_N^{-} \!\left(z^3 \right) = \frac{2\,z^3}{(2 \pi |z^3|)^{2N-2}} P \int_0^{\infty} dy\, y^{2N-2}\,e^{-y} \sum_{n=1}^{\infty}
\frac{1}{n^{2N-2} (y^2 - 4n^2 \pi^2 |z^3|^2)} \;\;,
\label{eightyseven}\end{eqnarray}
\begin{eqnarray}
SD_M^{+} \!\left( z^3 \right)= -\lfloor M/2  \rfloor \ln \!\left(-\, e^{\pm 2i \pi z^3}\right) - \frac{(1 - (-1)^M)}{2}\,
\ln \!\left(1- e^{\pm 2i \pi z^3} \right) \;\;,
\label{eightyeight}\end{eqnarray}
and 
\begin{eqnarray}
SD_M^{-}\!\left(z^3 \right) = (-1)^M \Bigl( \lfloor M/2 \rfloor + \frac{1- (-1)^M}{2} \Bigr) 2\pi |z^3| - \frac{1}{2} \,
\ln \Bigl( 1- e^{-2\pi |z^3|} \Bigr) \;\;.
\label{eightynine}\end{eqnarray}
The upper- and lower-signed versions of Eq.\ (\ref{eightyeight}) are valid for $(M-1/2) \pi/3 < \theta < (M+1/2) \pi/3$ and
$-(M+1/2) \pi/3 < \theta< -(M-1/2) \pi/3$, respectively, while Eq.\ (\ref{eightynine}) is valid for $\theta= \pm(M+1/2) \pi/3$. 
Therefore, for $-\pi< \theta \leq \pi$ or the principal branch for $z$, the Stokes lines occur at $\pm \pi/6$, $\pm \pi/2$, and
$\pm 5 \pi/6$. We shall investigate these forms after we have considered the results for the Stokes sectors first.

The Borel-summed asymptotic forms that are valid for the Stokes sectors can be expressed as
\begin{eqnarray}
\ln \Gamma \left(z^3\right)= {\small \begin{cases} 
F \left( z^3 \right) + TS_N \! \left( z^3 \right) + R_N^{+}\left( z^3 \right)  
-\ln \Bigl(-\,e^{2i \pi z^3}\Bigr) 
-\ln \Bigl(1-\,e^{2i \pi z^3}\Bigr), &    5\pi/6 < \theta \leq \pi, \\
F \left( z^3 \right) + TS_N \! \left( z^3 \right) + R_N^{+} \!\left( z^3 \right) -\ln \Bigl(-\,e^{2i \pi z^3}\Bigr),
&  \!\!\!\!\! \! \pi/2 < \theta< 5\pi/6, \\
F \left( z^3 \right) + TS_N \! \left( z^3 \right) + R_N^{+} \!\left( z^3 \right) 
-\ln \Bigl(1-\,e^{2i \pi z^3}\Bigr), &  \!\!\! \pi/6 < \theta \leq \pi/2, \\
F \left( z^3 \right) + TS_N \! \left( z^3 \right) + R_N^{+} \!\left( z^3 \right), & \!\!\!\!\!\!\!\!\!\!\!\!\!\!\!\!\!\!\!\!\!\! 
-\pi/6< \theta< \pi/6, \\  
F \left( z^3 \right) + TS_N \! \left( z^3 \right) + R_N^{+} \!\left( z^3 \right)  
-\ln \Bigl(1-\,e^{-2i \pi z^3}\Bigr), &   \!\!\!\!\!\!\!\!\!\!\!\!\!\!\!  -\pi/2 < \theta < -\pi/6, \\
F \left( z^3 \right) + TS_N \! \left( z^3 \right) + R_N^{+} \!\left( z^3 \right) -\ln \Bigl(-\,e^{-2i \pi z^3}\Bigr), &
\!\!\!\!\!\!\!\!\!\!\!\!\!\!\! -5\pi/6 < \theta< -\pi/2,\\ 
F \left( z^3 \right) + TS_N \! \left( z^3 \right) + R_N^{+} \!\left( z^3 \right) - \! \ln \Bigl(-\,e^{-2i \pi z^3}\Bigr) 
\! - \! \ln \Bigl(1-\,e^{-2i \pi z^3} \! \Bigr), &   \!\!\!\!\!\!\!  -\pi < \theta < -5\pi/6. \\
\end{cases}}
\label{ninety}\end{eqnarray}
If we compare the above results with the corresponding MB-regularized asymptotic forms given by Eqs.\ (\ref{seventysix})-(\ref{eightythree}), 
then we see that while they involve the same logarithmic terms, the main difference is that these terms emerge in different sectors of the
principal branch. For example, the Borel-summed asymptotic form for $\theta$ lying between $-\pi/6$ and $\pi/6$  does not possess a 
logarithmic term. Similarly, the first MB-regularized asymptotic form given by Eq.\ (\ref{seventysix}) does not possess a logarithmic 
term, but on this occasion $\theta$ lies between $-\pi/3$ and $\pi/3$. Moreover, according to the third result in Eq.\ (\ref{ninety}), 
there is a logarithmic term in the Borel-summed asymptotic form with $\theta$ lying in $(\pi/6,\pi/2)$. This term is identical to the 
logarithmic term appearing in the MB-regularized asymptotic form given by Eq.\ (\ref{seventyseven}), although the latter result is valid 
over a different sector, viz.\ $(0,2\pi/3)$.  

We have seen from Sec.\ 3 that in order to obtain very accurate results from the Borel-summed forms for $\ln \Gamma(z)$, the convergent 
sum in the remainder $R_N^{+}$ given by Eq.\ (\ref{twentytwo}) had to be truncated at a very large value despite the fact that $|z|$ 
was chosen to be relatively large. Typically, we replaced the upper limit of infinity in the sum by $10^5$ and chose $|z|$ to be equal
to 3. Selecting a much smaller value of $|z|$ increases the computation time significantly. Therefore, we shall select a larger value
of $|z|$ than in Table\ \ref{tab8}, but we cannot choose an integer because of the singularities that occur in the logarithmic
term of  $\ln \left( 1-\exp(\pm 2i \pi z^3 \right))$, when $|z|$ is an integer and $\theta = k\pi/3$. We were not confronted with this problem 
in Sec.\ 3 because $z$ was not replaced by $z^3$ in the Borel-summed asymptotic forms for the regularized value of $\ln \Gamma(z)$. Consequently,
we shall choose $|z|$ to be equal to 5/2. 

In setting an upper limit of $10^5$ in the convergent sum for the remainder $R_N^{+}$, we saw that the time taken to perform a
calculation increased markedly to over 4 hours on a Sony VAIO laptop with 2 GB RAM. In order to accelerate the process of obtaining
the results from many calculations, a script file of the program was created so that it could be run on a SunFire server mentioned 
previously. Even though each calculation took longer than a calculation on the Sony VAIO laptop, the saving in time occurred because many 
calculations were able to be performed simultaneously. Hence in order to carry out a proper comparison between Eqs.\ 
(\ref{seventysix})-(\ref{eightythree}) and the Borel-summed asymptotic forms given by Eq.\ (\ref{ninety}), we shall create another script 
file, which will be a composite of the third and fourth programs in the appendix. This program appears as Program 5 in the appendix. 

\begin{table}
\centering
{\small
\begin{tabular}{|c|c|c|c|} \hline
 $\theta$ & $N$ & Quantity  &  Value \\ \hline
& &$ F(z^3)$ &                    -14.88486001926988218316890967 - 30.6490797188237317042533096605$\, i$ \\
& & $TS_4(z^3)$ &                  0.001187233093636153159546970 + 0.00520018521369322415111941337$\, i$ \\
$-\pi/7$ & 4 & MB Int (M1=0) &     2.63163543021301667181503$ \times 10^{-12}$ - 6.6918000701164118190594$ \times 10^{-15}\, i$ \\
& & Log. Term (M1)   &  0 \\ 
& & Total via M1 &                -14.88367278617361439457914968 - 30.6438795645051371993532040127$\,i$ \\
& & MB Int (M2=-1) &               2.63163543021301667181503$ \times 10^{-12}$ - 6.6918000701164118190594$ \times 10^{-15}\, i$ \\
& & Log. Term (M2)  &             -2.67692877187577721034710$ \times 10^{-42}$ - 3.9146393659521104430365$ \times 10^{-43}\, i$ \\
& & Total via M2    &             -14.88367278617361439457914968 - 30.6438795645051371993532040127$\,i$ \\
& & Borel Rem       &              2.63163543021301667181503$ \times 10^{-12}$ - 6.6918000701164118190594$ \times 10^{-15}\, i$ \\
& & Borel Log Term  & 0 \\
& & Borel Total     &             -14.88367278617361439457914968 - 30.6438795645051371993532040127$\,i $ \\ 
& & LogGamma[zcube] &             -14.88367278617361439457914968 - 30.6438795645051371993532040127$\,i $ \\ \hline
& &$ F(z^3)$ &                     22.44243671730881854487411054 - 15.9297607498710220596760086394$\,i $ \\
& & $TS_2(z^3)$ &                  0.004927357506726862699350310 + 0.00204097830594714544921845324$\,i $ \\
$5\pi/8$ & 2 & MB Int (M1=1) &    -2.78985280328019887884107$ \times 10^{-7}$ - 6.7196224438310148304858$ \times 10^{-7}\,i$ \\
& & Log. Term (M1=1)   &          -37.56985811777164196399874512 + 0.40452594974678470711825665596$\,i $ \\ 
& & Total via M1    &             -15.12249432194137688444517216 - 15.5231944937805345902100165788$\,i $ \\   
& & MB Int (M2=2)   &             -2.78985280372387906470607$ \times 10^{-7}$ - 6.71962244402097138440407$ \times 10^{-7}\,i$ \\
& & Log. Term (M2)  &             -37.56985811777164191963072654 + 0.40452594974678472611391204778$\,i $ \\ 
& & Total via M2    &             -15.12249432194137688444517216 - 15.5231944937805345902100165788$\,i $ \\   
& & Borel Rem       &             -2.78985280372387820649852$ \times 10^{-7}$ - 6.71962244402096931250776$ \times 10^{-7}\,i$ \\
& & Borel Log Term  &             -37.56985811777164191963072654 + 0.40452594974678472611391204778$\,i $ \\ 
& & Borel Total &                 -15.12249432194137688444508634 - 15.5231944937805345902098093891$\,i $ \\ \hline  
& &$ F(z^3)$ &                     26.8706304919944588244828769703 + 0 $\,i\;\;\;\;\;\;\;\;\;\;\;\;\;\;\;\;\;\;\;\;\;\;\;\;\;\;\;\;\;\;\;\;\;\;\;\;\;\;\;\;\;\;\;\;\;\;\;\;\;$ \\
& & $TS_2(z^3)$ &                  0.0053333333333333333333333333 + 0 $\,i\;\;\;\;\;\;\;\;\;\;\;\;\;\;\;\;\;\;\;\;\;\;\;\;\;\;\;\;\;\;\;\;\;\;\;\;\;\;\;\;\;\;\;\;\;\;\;\;\;$ \\
$2\pi/3$ & 2 & MB Int (M1=2) &    -7.27328204587763887038193$ \times 10^{-7}$ + 0 $\,i\;\;\;\;\;\;\;\;\;\;\;\;\;\;\;\;\;\;\;\;\;\;\;\;\;\;\;\;\;\;\;\;\;\;\;\;\;\;\;\;\;\;\;\;\;\;\;\;\;\;$ \\ 
& & Log. Term (M1=2)   &    $\;\;\;\;\;\;\;\;\;\;\;\;\;\;\;\;\;\;\;\;\;\;\;\;\;\;\;\;\;\;\;\;\;\;\;\;\;\;\;\;\;\;\;\;\;\;\;\;\;\;\;\;\;$ - 0.785398163397448309615660845819$\,i $ \\ 
& & Total via M1 &                 26.875963097999587570052547525 - 0.785398163397448309615660845819$\,i $ \\
& & Borel Rem &                   -7.27328204587763662777752$ \times 10^{-7}$ + 0 $\,i\;\;\;\;\;\;\;\;\;\;\;\;\;\;\;\;\;\;\;\;\;\;\;\;\;\;\;\;\;\;\;\;\;\;\;\;\;\;\;\;\;\;\;\;\;     \;\;\;\;\;$ \\
& & Borel Log Term &       $\;\;\;\;\;\;\;\;\;\;\;\;\;\;\;\;\;\;\;\;\;\;\;\;\;\;\;\;\;\;\;\;\;\;\;\;\;\;\;\;\;\;\;\;\;\;\;\;\;\;\;\;\;$ - 0.785398163397448309615660845819$\, i$ \\
& & Borel Total &                  26.875963097999587570052323265 - 0.785398163397448309615660845819$\,i $ \\ \hline
& &$ F(z^3)$ &                    -45.81050523277279074549656662 - 7.88812400847329603213924107797$\, i$ \\
& & $TS_4(z^3)$ &                 -0.003771750463193365200048475 - 0.00377072066430421316496938343$\, i$ \\
$11\pi/12$ & 4 & MB Int (M1=2) &   1.8403031527099105137119$ \times 10^{-12}$ - 1.86174503255187256069532$ \times 10^{-12}\,i$ \\ 
& & Log Term (M1=2) &              69.42004590872447260962314046 - 2.83658512384077187501765734573$\, i$ \\
& & Total via M1 &                 23.60576892549032880207923528 - 10.7284798529802338653544196797$\, i$ \\
& & MB Int (M2=3) &                1.8403031527099105130346$ \times 10^{-12}$ - 1.86174503255187256048211$ \times 10^{-12}\, i$ \\
& & Log Term (M2=3) &              69.42004590872447260962314046 - 2.83658512384077187501765734573$\, i$ \\
& & Total via M2 &                 23.60576892549032880207923528 - 10.7284798529802338653544196797$\, i$ \\
& & Borel Rem &                    1.8403031527099105130346$ \times 10^{-12}$ - 1.86174503255187256048211$ \times 10^{-12}\, i$ \\
& & Borel Log Term &               69.42004590872447260962314046 - 2.83658512384077187501765734573$\, i$ \\
& & Borel Total  &                 23.60576892549032880207923528 - 10.7284798529802338653544196797$\, i$ \\ \hline
\end{tabular}
}
\normalsize
\caption{Determination of $\ln\Gamma(z^3)$ for $|z| \!=\! 5/2$ and various values of $\theta$ and $N$}
\label{tab9}
\end{table}

Table\ \ref{tab9} presents a small sample of the results obtained by running Program 5 in the appendix on a SunFire server. 
The program prints out the results for the MB-regularized asymptotic forms first, then the results from Borel-summed asymptotic
forms given by Eq.\ (\ref{ninety}) and finally, the results from Mathematica's LogGamma routine, whenever possible. Most of the 
calculations took between 4 and 6 hrs to execute, although some did take much longer such as the $\theta \!=\! 4\pi/7$ and 
$N \!=\! 3$ calculation, which took 7.5 hours. 

The first calculation in Table\ \ref{tab9} presents the results obtained for $\theta \!=\! -\pi/7$ and $N \!=\! 4$. Because the optimal
point of truncation occurs approximately at $N_{OP} \!=\! 8$ according to the discussion immediately below ``Eq." (\ref{fiftyseven}),
we expect that the Stirling approximation or $F(z^3)$ to yield a good approximation to the actual value of $\ln \Gamma((5/2)^3 
\exp(-3i\pi/7))$. This turns out to be case when we compare the first value with the LogGamma[zcube] at the bottom of the calculation.
Then we see that the value of $F(z^3)$ is accurate to the third decimal place or the first four figures. As a consequence, the
truncated sum is small and only makes a contribution at the third decimal place. For this calculation there are two valid MB-regularized
asymptotic forms, viz.\ Eqs.\ (\ref{seventysix}) and (\ref{eightyone}), which are denoted by M1=0 and M2=-1, respectively. MB Int
denotes the value obtained from the MB integrals in these asymptotic forms. Both are of the order of $10^{-12}$. In the case of M1=0, 
there is no logarithmic term accompanying the MB integral and thus, it has been set equal to zero, while for M2=-1, there is a contribution,
but it is almost negligible, since it is the order of $10^{-42}$. This means that according to the accuracy and precision goals set 
in the NIntegrate routines, the M1=0 and M2=-1 calculations are virtually identical to one another, which is reflected by the 
fact that the MB integrals display identical values in the table. Consequently, the totals representing the sum of $F(z^3)$ and their
respective MB integrals and logarithmic terms, are identical to one another. Therefore, we expect these values to be more accurate
than the results from the Borel-summed asymptotic forms since the latter have been truncated. In fact, we see that the results from both
MB-regularized asymptotic forms are identical to the value obtained via Mathematica's LogGamma function. 

The result labelled Borel Rem represents the Borel-summed remainder or $R_N^{+}(z^3)$ in the fourth asymptotic form of Eq.\ (\ref{ninety}) 
when it is truncated at $10^5$. Despite the truncation it is identical to the values obtained from the MB integrals above it. In actual fact, 
the Borel Rem value was identical to the first 34 decimal figures of the MB integrals, which is well outside the accuracy and precision goals 
in the script file. Bearing in mind, that the remainder is very small, this means that only the first 13 or so decimal figures of each remainder
calculation will contribute to the totals. To observe the effect of the remaining figures and the logarithmic term for the M1=-1 calculation,
we need to extend WorkingPrecision, AccuracyGoal and PrecisionGoal to much higher values in the various calls to the NIntegrate routine,
which will come at the expense of the time of execution. Nevertheless, the results for $\theta \!=\! -\pi/7$ and $N \!=\! 4$ confirm the
validity of the fourth form in Eq.\ (\ref{ninety}), although it is not sufficient.

The second calculation in Table\ \ref{tab9} presents the results obtained for $\theta \!=\! 5\pi/8$ and $N \!=\!2$. The first observation to 
be made about this calculation is that there is no value from the Mathematica's LogGamma routine, which will apply to all subsequent calculations 
in the table. Hence the only means of obtaining the value of $\ln\Gamma(z^3)$ is by using either the MB-regularized or the Borel-summed asymptotic
forms for the function. Next we observe that unlike the previous calculation, the value of $F(z^3)$ is nowhere near the actual of value of the 
function when it is compared with any of the three totals, i.e. Total via M1, Total via M2 and and Borel Total, even though $|z|$ is
relatively large. In standard Poincar$\acute{\rm e}$ asymptotics this is explained by the fact that we have crossed into another Stokes 
sector. Consequently, the reason for the huge discrepancy between the final totals and the value of $F(z^3)$ is due to the Stokes discontinuous 
term or the logarithmic term in the second asymptotic form of Eq.\ (\ref{ninety}), whose value appears in the row labelled Borel Log Term. Appearing
immediately below the value of $F\left((5/2)^3 \exp(15\pi i/8)\right)$ is the value of the truncated sum $TS_2\left((5/2)^3 \exp(15 \pi i/8) \right)$.
As in the previous calculation this value is small in comparison with the value of the Stirling approximation, only beginning to make
a contribution at the third decimal place when the actual value is four orders of magnitude higher. In this calculation we once again have 
two MB-regularized asymptotic forms for the remainder, which are given by Eqs.\ (\ref{seventyseven}) where M1=1, and (\ref{seventyeight}) 
where M2=2. Both these asymptotic forms for the remainder yield small values, but because $N$ is further from the optimal point of truncation,
they are larger in magnitude (of the order of $10^{-7}$) than in the previous calculation. Therefore, they agree with each other to a 
lesser number of decimal figures, 10 rather than 24. Below the MB Int (M1=1) value appears the value of logarithmic term in Eq.\ (\ref{seventyseven}),
which dominates the calculation of the real part of $\ln \Gamma \left( (5/2)^3 \exp(15 \pi i/8) \right)$, but not the imaginary part. The
latter is dominated by the imaginary part of the Stirling approximation. Appearing immediately below the MB Int (M2=2) value is the value
of the logarithmic term in Eq.\ (\ref{seventyeight}). Although this term is slightly different to the logarithmic term in Eq.\ (\ref{seventyseven}),
both values still agree with each other to 15 decimal figures. Once again, we find that the logarithmic term dominates the real part of 
the final value of $\ln \Gamma\left( (5/2)^3 \exp(15 \pi i/8) \right)$, while the imaginary part is dominated by the imaginary part of
the Stirling approximation. Below the Total via M2 value are the various terms appearing in the second asymptotic form of Eq.\ (\ref{ninety}).
The first value is the value of the remainder $R_2^{+}\left( (5/2)^3 \exp(15 \pi i/8)\right)$, which is not only of the same magnitude
as the MB Int values, viz.\ $10^{-7}$, but also agrees with them to a large number of decimal figures as in the previous calculation.
This is not so surprising because the Borel Log Term value is identical to the logarithmic term in Eq.\ (\ref{seventyeight}).  Finally, when
the Borel Rem and Borel Log Term values are added to the Stirling approximation and the value of the truncated sum, we obtain the Borel Total,
which is identical to Total via M1 and Total via M2 values.

The third calculation in Table\ \ref{tab9} presents the results obtained for $\theta \!=\! 2\pi/3$ and $N \!=\! 2$. For this value of $\theta$ 
there is only one valid MB-regularized asymptotic form, which is given by Eq.\ (\ref{seventyeight}). Moreover, the values of the Stirling 
approximation, the truncated sum and the MB integral are real, while the Log Term yields the imaginary contribution, which equals $-\pi/4$. 
That is, $ \Im \ln \Gamma \left( |z|^3 e^{2i\pi} \right) = -\pi/4 $. As expected, the value of the MB integral is very small of the order 
of $10^{-7}$. Therefore, unlike the previous calculation, the Stirling approximation now provides a very accurate value for the real part of 
$\ln \Gamma\left((5/2)^3 \exp(2i\pi) \right)$. Below the Total via M1 value appears the remainder of Borel-summed version for 
$\ln \Gamma \left( |z|^{3} \exp(2i \pi) \right)$, which has been obtained by evaluating $R_N^{+}(z^3)$ in the second result in Eq.\ (\ref{ninety}).
This result has the same logarithmic term as in Eq.\ (\ref{seventyeight}). Thus, we expect that the calculation via the Borel-summed 
form to be identical to the value obtained via the MB-regularized form. If we look closely at the two totals, we the see that they agree
to 23 decimal figures, not the expected 30 specified by the accuracy and precision goals. Since the Stirling approximation, the truncated sum  
and the logarithmic term all agree with each other, the discrepancy must occur in the evaluation of the MB integral and Borel-summed
remainder. Since the former represents only one intergal, it is expected to be more than the result obtained from $R_2^{+}(z^3)$, which has been
truncated to $10^5$. Hence this is an example where the remainder in the Borel-summed version needs to be truncated at a much higher value
in order to achieve the desired accuracy. In addition, when it comes to reliability, it is often much better to use the results 
obtained via MB-regularized asymptotic forms than those via Borel summation as was first observed in Ref.\ \cite{kow95}.

The final calculation in Table\ \ref{tab9} are the results obtained for $\theta \!=\! 11\pi/12$ and $N \!=\! 4$. In this calculation there are 
two MB-regularized asymptotic forms for $\ln \Gamma \left( |z|^3 \exp(11i \pi/4) \right)$, which are given by Eqs.\ (\ref{seventyeight}) where 
M1=2 and (\ref{seventynine}) where M2=3. The Borel-summed asymptotic form for this calculation is given by the first form in Eq.\ (\ref{ninety}). 
Once again, we find that the real part is nowehere near the real part of the actual value for the function, although the imaginary part is not too 
far away. However, if the logarithmic term for all three forms is added to the Stirling approximation, then one obtains a good approximation for 
$\ln \Gamma \left( (5/2)^3 \exp(11 i \pi/4) \right)$. In this instance the logarithmic term for the Borel-summed form is identical to that in Eq.\ 
(\ref{seventynine}), but by comparing it with the value obtained via the M1=2 form, the extra term in the former is very small indeed, only affecting 
the 20 decimal place of the term, which is itself of the order of $10^{-12}$. Consequently, all three totals agree with each other as in the other 
examples displayed in the table.

There is yet another study that we need to perform before we can say with certainty that all the Borel-summed asymptotic forms in Thm.\ 2.1 agree 
with their respective MB-regularized asymptotic forms given in Thm.\ 4.1. Specifically, we need to investigate those asymptotic forms in Thm.\ 2.1 
that apply to the Stokes lines to see if they yield the same results as their MB-regularized counterparts. To obtain the Borel-summed 
asymptotic forms, we introduce $R_N^{-}(z^3)$ and $SD_M^{-}(z^3)$ as given by Eqs.\ (\ref{eightyseven}) and (\ref{eightynine}), respectively, into 
Eq.\ (\ref{eightyfive}). By putting $M$ equal to 0, 1 and 2, we obtain the specific forms for the Stokes lines where $\theta$ equals $\pm \pi/6$, 
$\pm \pi/2$, and $\pm 5\pi/6$. Hence after a little manipulation, we arrive at 
\begin{eqnarray}
\ln \Gamma \left(z^3\right)= {\small \begin{cases} 
F \left( z^3 \right) + TS_N \! \left( z^3 \right) + R_N^{-}\left( z^3 \right)  
- \frac{1}{2} \; \ln \Bigl(1-\,e^{-2 \pi |z^3|}\Bigr), &  \theta = \pm \pi/6 , \\
F \left( z^3 \right) + TS_N \! \left( z^3 \right) + R_N^{-}\left( z^3 \right)  
- \frac{1}{2} \; \ln \Bigl(1-\,e^{-2 \pi |z^3|}\Bigr)  & \\
- 2\pi |z^3|, &    \theta = \pm \pi/2 , \\
F \left( z^3 \right) + TS_N \! \left( z^3 \right) + R_N^{-}\left( z^3 \right)  
- \frac{1}{2} \; \ln \Bigl(1 -\,e^{-2 \pi |z^3|}\Bigr) & \\
+ 2 \pi |z^3| , &  \!\!\!  \theta = \pm 5\pi/6 . \\
\end{cases}}
\label{ninetyone}\end{eqnarray}
Note the similarity of the Stokes discontinuity terms with the corresponding terms or $S_{MB}(z^3)$ in the MB-regularized asymptotic forms 
given in Eqs.\ (\ref{seventyseven}), (\ref{seventysevenb}), (\ref{seventyninea})-(\ref{eightyone}), (\ref{eightytwob})-(\ref{eightytwod}) and 
(\ref{eightythreea}). In fact, the major difference occurs with the logarithmic term, which is represented by either a zero or full residue 
contribution in the MB-regularized asymptotic forms, while it is always represented by a semi-residue contribution in the Borel-summed 
asymptotic forms. Hence the above asymptotic forms possess a factor of 1/2 outside their logarithmic terms.

Table\ \ref{tab10} presents a small sample of the results obtained by running the final program in the appendix. Because the program utilizes 
the Borel-summed asymptotic forms in the calculations, which we have already seen are very time-consuming in an accurate evaluation, it appears 
as a Mathematica script file, thereby enabling all the results for each Stokes line to be calculated simultaneously on a SunFire alpha 
server. Each calculation performed on the alpha server took almost 30 hours, primarily because of the heavy computation required in evaluating
the remainder of the Borel-summed asymptotic forms, which were truncated at $10^5$ terms as done previously. The table displays the results obtained 
for $|z|=9/10$, although $|z|=1/10$ was also considered in this study. Interestingly, the Borel-summed asymptotic forms yielded more accurate 
results than the corresponding MB-regularized asymptotic forms for $|z|=1/10$. Specifically, for $\theta=|\pi/6|$, the Borel-summed 
asymptotic forms yielded values that agreed with the LogGamma routine in Mathematica to 33 decimal places, whereas the results obtained
via the MB-regularized asymptotic forms only agreed to 13 decimal places. This situation arose because the truncated sum and MB integrals in 
the MB-regularized asymptotic forms were of the order of $10^{23}$ when $|z|=1/10$. Hence a cancellation of a large number of decimal 
figures occurred before the final value for $\ln \Gamma(z^3)$ was obtained. To circumvent this problem, we need to increase the working precision 
and the precision and accuracy goals in the the MB-regularized asymptotic forms substantially when dealing with very low values of $|z|$, although
this will come at a cost in the time required to do the calculations. It should also be borne in mind that the variable in this study is,
in reality, $z^3$. Hence the calculations are actually being performed for a value of $10^{-3}$ when $|z|=1/10$. On the other hand, the MB 
integrals yielded values of the order of $10^{-3}$ for $|z|=9/10$, which means that there was no significant cancellation of decimal figures 
occurring with the results displayed in Table\ \ref{tab10}.

\begin{table}
\small
\centering
\begin{tabular}{|c|c|c|} \hline
$\theta$ & Method & Value \\ \hline
 &   Eq. (\ref{seventysix}) &        -0.0629795852996006019126614 - 1.86781980997058048039434088$\,i$ \\
$\pi/6$ & Eq. (\ref{seventyseven}) &    -0.0629795852996006019126614 - 1.86781980997058048039434088$\,i$ \\
& Top, Eq. (\ref{ninetyone})&             -0.0629795852996006019126614 - 1.86781980997058048039434088$\,i$ \\       
 &  Mathematica &     -0.0629795852996006019126614 - 1.86781980997058048039434088$\,i$ \\ \hline    
 & Eq. (\ref{eightytwob}) &           -4.6434216742335191435911954 - 1.86781980997058048039434088$\, i$ \\
$-\pi/2$ & Eq. (\ref{eightytwoc}) &   -4.6434216742335191435911954 - 1.86781980997058048039434088$\, i$ \\
& Middle, Eq. (\ref{ninetyone}) &           -4.6434216742335191435911954 - 1.86781980997058048039434088$\, i$  \\ \hline
& Eq. (\ref{eighty}) &             4.5174625036343179397658872 - 1.86781980997058048039434088$\, i $ \\
$5\pi/6$ & Eq. (\ref{eightyone}) &    4.5174625036343179397658872 - 1.86781980997058048039434088$\, i$ \\
& Bottom, Eq. (\ref{ninetyone}) &            4.5174625036343179397658872 - 1.86781980997058048039434088$\, i$  \\ \hline
\end{tabular}
\normalsize
\caption{Values of $\ln\Gamma \! \left( z^3\right)$ for $|z| =9/10$ at the Stokes lines where $\theta$ equals $\pi/6$, 
$-\pi/2$ and $5 \pi/6$}
\label{tab10}
\end{table}

The first column of Table\ \ref{tab10} displays the value of $\theta$ or the particular Stokes lines being considered. These are
the Stokes lines at: (1) $\theta \!=\! \pi/6$,  (2) $\theta \!=\! -\pi/2$ and (3) $\theta \!=\! 5\pi/6$. As indicated previously,
$\ln \Gamma(z^3)$ cannot be evaluated for the last two lines by Mathematica. Consequently, there is an extra result for the calculations
at $\theta \!=\! \pi/6$. The second column of Table\ \ref{tab10} displays the equation that was used to calculate the value of
$\ln \Gamma(z^3)$, while the third column displays the actual values to 27 decimal places. We see that not only do the two
different MB-regularized asymptotic forms agree with one another at each Stokes line, but they also agree with the results
obtained from the the Borel-summed asymptotic forms in Eq.\ (\ref{ninetyone}). In addition, we see that for the Stokes line
at $\theta \!=\! \pi/6$, the Borel-summed and MB-regularized asymptotic forms agree with the value obtained from the LogGamma routine
in Mathematica.

\section{Conclusion}
In this work Stirling's approximation or formula for the special function $\ln \Gamma(z)$ has been exactified. Exactification
is defined as the process of obtaining the values of a function from an asymptotic expansion. Two steps are required for this
process to be successful. First, one must ensure that the asymptotic expansion is complete, which means the inclusion of frequently 
neglected subdominant or decaying exponential terms, over all arguments or phases of the power series variable in the expansion. The 
combination of a complete asymptotic expansion and the sector or line over which it is valid is referred to as an asymptotic 
form. The second step in exactification is to regularize all asymptotic series appearing in the asymptotic forms, where
regularization is defined as the removal of the infinity in the remainder of an asymptotic series so as to make the series summable. 
By itself, regularization represents an abstract mathematical technique when applied to divergent series, but it is necessary in
asymptotics because it represents the means of correcting the asymptotic method used in the derivation of the complete asymptotic
expansion. By carrying out the two steps, it has been possible to obtain exact values of $\ln \Gamma(z)$ to incredible accuracy, 
often more than 30 decimal figures, which was the restriction placed on the computing system used in this work. The computing
system was composed of powerful computers and mathematical software in the form of Mathematica \cite{wol92}. In those rare cases where 
30 figure accuracy was not obtained, it was generally attributed to having chosen too high or too low a value for the truncation parameter 
$N$ in the asymptotic forms, which resulted in very large values for the remainder. In first situation a large number of decimal figures 
in the remainder was cancelled by the large values of the truncated asymptotic series in the asymptotic forms. As a result of this 
cancellation, the capacity of the computing system to provide values within the accuracy and precision goals set in the programs was 
reduced. However, by extending the working precision, in particular, one was able to obtain the desired accuracy again, but at the 
expense of a substantial increase in computing time. The second case of very low values of the truncation parameter is discussed further
below.

In this work two techniques has been used to regularize the asymptotic series in a complete asymptotic expansion: (1) Borel summation 
and (2) Mellin-Barnes (MB) regularization. Though different, they both yield the same regularized values of an asymptotic series. In
the case of Borel summation the asymptotic forms for a function are different across adjacent Stokes sectors, whose boundaries
are separated by lines of discontinuity known as Stokes lines or rays. The asymptotic forms arising from the MB regularization of
an asymptotic series are dependent upon domains of convergence for specific MB integrals, which represent the remainder of the 
asymptotic series. Unlike Borel-summed asymptotic forms, the domains of convergence overlap, which means in turn that often
two MB-regularized forms can be used to obtain the regularized value of an asymptotic series in the overlapping sectors. Therefore, 
they serve as a check on whether the regularization process has been conducted properly.

Because this work has been concerned with Stirling's approximation for $\ln \Gamma(z)$ rather than $\Gamma(z)$, it has resulted
in the derivation of asymptotic forms quite unlike those derived previously in Refs.\ \cite{kow95}, \cite{kow002}, and
\cite{kow001}-\cite{kow11c}. In this instance the asymptotic forms have been derived as a result of summing an infinite number
of a particular generalized Type I terminant denoted by $S^{I}_{2,-1}(N,(1/2n\pi z)^2)$ over $n$. The Borel-summed regularized 
value of this generalized Type I terminant is composed of a Cauchy integral and a finite number of exponential terms depending
on how many Stokes sectors have been crossed in accordance with the Stokes phenomenon. As discussed in Ref.\ \cite{kow09}, the 
exponential terms arise from residues of the singularities in the Cauchy integral. Normally, they would not pose a problem, but 
because of the infinite sum over $n$ in the asymptotic forms, they result in logarithmically divergent series, which need to be 
regularized according to Lemma\ 2.2. 

A similar situation develops when evaluating the MB-regularized value, which involves the same 
generalized Type I terminant. Although the sum over $n$ in the MB integrals reduces to a more compact form where the infinite sum is 
replaced by the Riemann zeta function, we still obtain logarithmically divergent series similar to those in the Borel-summed asymptotic 
forms. In this case the exponential terms arise from the differences in the MB integrals when they are valid over the common 
or overlapping sector of their domains of convergence. Although these series can be MB-regularized, Lemma\ 2.2 was used again to 
regularize them as it produces compact logarithmic terms, which do not require numerical integration. That is, although the results in 
Thm.\ 4.1 are not purely MB-regularized asymptotic forms, they are more compact since they avoid the computation of extra MB integrals. 
Furthermore, the logarithmic terms, which are denoted by $SD_M^{SL}(z)$ or $SD_M^{SS}(z)$ and $SD_{MB}(M,z)$, respectively, in the 
Borel-summed and MB-regularized asymptotic forms, are not only necessary for exhibiting the multivaluedness of $\ln \Gamma(z)$, but also
enable values of the special function to be evaluated beyond the principal branch of $-\pi < {\rm arg}\,z \leq \pi$, which the LogGamma 
routine in Mathematica \cite{wol92} is unable to do. 
    
We have also seen in this work that the MB-regularized asymptotic forms are far more expedient for obtaining values of $\ln \Gamma(z)$
than their Borel-summed counterparts because the latter possess an infinite convergent sum of exponential integrals, which is not as
amenable to numerical computation as the MB integrals for the remainder of the asymptotic series $S(z)$ given by Eq.\ (\ref{sixteen}). As a 
consequence, the Borel-summed asymptotic forms must be truncated at a very large number, $10^5$ here, in order to match the accuracy of 
the MB-regularized remainders, which, as mentioned above, benefit from having their infinite sum replaced by the Riemann zeta function. 
It should also be noted that the selection of a very low value for the truncation parameter could affect the accuracy of the Borel-summed 
remainder, but for these cases the problem could be overcome by increasing the limit well beyond $10^5$, again at the expense of 
computer time. Thus, the Borel-summed asymptotic forms for the regularized remainder take significantly longer 
to compute than their MB-regularized counterparts. Another disadvantage with the Borel-summed asymptotic forms is that they are vastly 
different along Stokes lines compared with Stokes sectors and thus, separate programs need to be written, whereas only one program is 
required for obtaining the values of $\ln \Gamma(z)$ via its MB-regularized asymptotic forms in Thm.\ 4.1.

In conclusion, this work has profound implications for both (hyper)asymptotics and mathematics in general. First, we have seen that 
it is permissible to differentiate the asymptotic expansion for a differentiable function provided it is complete, which contradicts 
widely accepted dogma in standard Poincar$\acute{{\rm e}}$ asymptotics \cite{whi73}. Next, as a result of regularization, one can
obtain exact values of a multivalued function from its complete asymptotic expansion for all values of the argument of the main
variable. Consequently, asymptotics has been transformed into a proper mathematical discipline capable of eliciting exact values rather 
than one suffering from the drawbacks of vagueness, limited accuracy and non-specific ranges of validity/applicability. Finally, this 
work has demonstrated that mathematicians should no longer fear divergence by trying to avoid it at all costs, but should confront
it by taking on the challenge to tame it since this will produce new and interesting mathematics.     

\section{Acknowledgement}
The author wishes to thank Prof. M.L. Glasser, Clarkson University, for his interest and encouragement of this work. Interesting
discussions with Dr R.B. Paris, University of Abertay, have also been used in producing this work.

The corresponding author states that there is no conflict of interest.

\newpage
\section{Appendix}
In this appendix we present the Mathematica modules \cite{wol92} that have been used in the various numerical studies
appearing in this work. They are displayed here so that the adept reader can verify the results in the numerous tables. 
Although condensed to minimize space, they are nevertheless readable since Mathematica constructs have been set in boldface.

\subsection*{Program 1}
The Mathematica program/module used to evaluate values of $\ln \Gamma(z)$ over the Stokes sectors in the principal branch
of the complex plane for $z$ by using the Eq.\ (\ref{sixtyeight}) for the remainder $R_N^{+}(z)$ appears below: 
\vspace{0.5 cm}
\newline
c1[k$_{-}$] := -2 {\bfseries Zeta}[2 k]/{\bfseries Pi}$^{\wedge}$(2 k) \newline
Smand[z$_{-}$] := (-1)$^{\wedge}$k {\bfseries Gamma}[2 k - 1] c1[k]/(2 z)$^{\wedge}$(2 k) \newline
F[modz$_{-}$, theta$_{-}$] := (modz {\bfseries Exp}[{\bfseries I} theta] - 1/2) ({\bfseries Log}[modz] + {\bfseries I} theta) 
- modz {\bfseries Exp}[{\bfseries I} theta] + Log[2 {\bfseries Pi]}/2\newline
Strlng[TP$_{-}$, modz$_{-}$, theta$_{-}$, limit$_{-}$] := {\bfseries Module}[\{\}, z = modz {\bfseries Exp}[{\bfseries I} theta]; 
\newline e0 = F[modz, theta]; \newline {\bfseries Print}["The value of the leading terms or Stirling's approximation \
is ", {\bfseries N}[e0, 50] // {\bfseries FullForm}]; e1 = z {\bfseries Sum}[Smand[z], {k, 1, TP - 1}];tot = 0; 
\newline {\bfseries Print}["The value of the truncated sum when N equals ", TP, " is ", {\bfseries N}[e1, 50] //{\bfseries FullForm}]; 
{\bfseries Do}[e2 = {\bfseries Exp}[-2 k {\bfseries Pi} z {\bfseries I}]*{\bfseries Gamma}[2 - 2 TP, -2 k {\bfseries Pi} z {\bfseries I}]/k; 
 e3 = -{\bfseries Exp}[2 k {\bfseries Pi} z {\bfseries I}]*{\bfseries Gamma}[2 - 2 TP, 2 k {\bfseries Pi} z {\bfseries I}]/k; 
 tot = tot + e2 + e3 , \{k, 1, limit\}]; \newline rem = {\bfseries Gamma}[2 TP - 1]*tot/(2 {\bfseries Pi I}); 
{\bfseries Print}["For z=", z, " and N=", TP,  " the value of the remainder is found to equal ", 
{\bfseries  N}[rem, 50] // {\bfseries FullForm}]; \newline {\bfseries Which}[-{\bfseries Pi} $<$ theta $<$ -{\bfseries Pi}/2, 
 e4 = -{\bfseries Log}[1 - {\bfseries Exp}[-2 {\bfseries Pi I} z]], -{\bfseries Pi}/2 $<$ theta $<$ {\bfseries Pi/2}, e4 = 0, 
{\bfseries Pi}/2 $<$ theta $<$= {\bfseries Pi}, e4 = -{\bfseries Log}[1 - {\bfseries Exp}[2 {\bfseries Pi I} z]]]; \newline
{\bfseries Print}[" For z=", z, " the value of the Stokes discontinuity term is equal to ", {\bfseries N[}e4, 50] // {\bfseries FullForm}];
e5 = e0 + e1 + rem + e4; \newline {\bfseries  Print}[" Combining all contributions for z=", z, " yields ", 
{\bfseries N}[e5, 50] // {\bfseries FullForm}]; \newline e5 = {\bfseries LogGamma}[z];\newline 
{\bfseries Print}["For z=", z,  " the value of Log($\Gamma$(z)) is equal to ",{\bfseries N}[e5, 50] //{\bfseries FullForm}]]
\vspace{0.5 cm}
\newline
The first three lines must be entered separately prior to the main module strlng, which requires four 
input values in order to execute. The first is the truncation parameter, denoted here by TP, while the second 
and third are respectively modz or $|z|$ and $\theta$ or ${\rm arg}\,z$. The final value called limit represents 
the value at which the infinite series in Eq.\ (\ref{sixtyeight}) is truncated. This variable appears only as the 
upper limit in the Do loop. The larger it is, the more accurate the remainder becomes, but at the expense of the 
time taken to execute the program. Note that e2 and e3 represent the terms in the remainder involving the 
incomplete gamma function. The Stokes discontinuity terms are evaluated by means of the Which statement, which 
evaluates the variable e4 according to whichever Stokes sector $\theta$ lies in. LogGamma[z] represents 
the value obtained via Mathematica's own routine for $\ln \Gamma(z)$. 

\subsection*{Program 2}
The previous program relies on the fact that an intrinsic numerical routine exists that can evaluate the incomplete 
gamma function quickly and extremely accurately. The following module evaluates $\ln \Gamma(z)$ by using the numerical 
integration routine NIntegrate in Mathematica. 
\vspace{0.5 cm}
\newline
c1[k$_{-}$] := -2 {\bfseries Zeta}[2 k]/{\bfseries Pi}$^{\wedge}$(2 k) \newline
Smand[z$_{-}$] := (-1)$^{\wedge}$k {\bfseries Gamma}[2 k - 1] c1[k]/(2 z)$^{\wedge}$(2 k) \newline
Intgrd[N$_{-}$, z$_{-}$, k$_{-}$] := y$^{\wedge}$(2 N - 2) {\bfseries Exp}[-y]/(y$^{\wedge}$2 + (2 {\bfseries Pi} k z)$^{\wedge}$2) 
\newline Strlng2[TP$_{-}$, modz$_{-}$, theta$_{-}$, limit$_{-}$] := 
{\bfseries Module}[\{\}, z = modz {\bfseries Exp}[{\bfseries I} theta]; \newline e0 = F[modz, theta]; \newline 
{\bfseries Print}["The value of the leading terms or Stirling's approximation 
is ", {\bfseries N}[e0, 25] // {\bfseries FullForm}]; \newline e1 = z {\bfseries Sum}[Smand[z], \{k, 1, TP-1\}]; \newline 
{\bfseries Print}["The value of the truncated sum when N equals ", TP, " is ", 
{\bfseries N}[e1, 25] // FullForm]; \newline Rem = 0; \newline
{\bfseries Do}[e2 = {\bfseries NIntegrate}[Intgrd[TP, z, k], \{y, 0, 1/2, 1, 2, 5, 10, 100, 1000, {\bfseries Infinity}\}, 
{\bfseries MinRecursion} -$>$ 3, {\bfseries MaxRecursion} -$>$ 10, {\bfseries WorkingPrecision} -$>$ 60, 
{\bfseries AccuracyGoal} -$>$ 30, {\bfseries PrecisionGoal} -$>$ 30]; e2 = e2/k$^{\wedge}$(2 TP - 2); Rem = Rem + e2, \{k, 1, limit\}]; 
\newline Rem = 2 (-1)$^{\wedge}$(TP + 1) z Rem/(2 {\bfseries Pi} z)$^{\wedge}$(2 TP - 2) ; \newline 
{\bfseries Print}["The value for the remainder when N equals ", TP, " is ", 
{\bfseries N[}Rem, 25] // {\bfseries FullForm}];\newline 
{\bfseries Which}[ -{\bfseries Pi} $<$ theta $<$ -{\bfseries Pi}/2, e3 = -{\bfseries Log}[1 - {\bfseries Exp}[-2 {\bfseries Pi I} z]], 
-{\bfseries Pi}/2 $<$ theta $<$ {\bfseries Pi}/2, e3 = 0, {\bfseries Pi}/2 $<$ theta $<$= {\bfseries Pi}, 
e3 = -{\bfseries Log}[1 - {\bfseries Exp}[2 {\bfseries Pi I} z]]];\newline 
{\bfseries Print}[" For z=", z," the value of the Stokes discontinuity term is equal to ", 
{\bfseries N[}e3, 25] // {\bfseries FullForm}]; \newline e4 = e0 + e1 + Rem + e3; \newline 
{\bfseries Print}["Combining all contributions for z=", z, " yields ", 
{\bfseries N}[e4, 25] // {\bfseries FullForm}]; \newline e5 = {\bfseries LogGamma}[z]; \newline 
{\bfseries Print}["For z=", z," the value of Log($\Gamma$(z)) is equal to ", 
{\bfseries N}[e5, 25] // {\bfseries FullForm}]]

\subsection*{Program 3}
The following Mathematica module uses the same statements for F(z) and Smand[z] as in the first program, which must be 
introduced prior to the module called strlngline. The main difference between this and the previous program is that the 
Do loop in the second program evaluates the remainder by using the NIntegrate routine according to Eq.\ (\ref{sixty}), whereas 
one needs to specify where the singularity occurs on a Stokes line. \newline
To calculate the values on the Stokes lines, we require a third program, which is:\vspace{0.5 cm} \newline
c1[k$_{-}$] := -2 {\bfseries Zeta}[2 k]/{\bfseries Pi}$^{\wedge}$(2 k) \newline
Intgrd[N$_{-}$, z$_{-}$, snglrty$_{-}$] := y$^{\wedge}$(2 N - 2) {\bfseries Exp}[-y]/(y$^{\wedge}$2 - snglrty$^{\wedge}$2)\newline
Smand[z$_{-}$] := (-1)$^{\wedge}$k {\bfseries Gamma}[2 k - 1] c1[k]/(2 z)$^{\wedge}$(2 k) \newline
F[modz$_{-}$, theta$_{-}$] := (modz {\bfseries Exp}[I theta] - 1/2) ({\bfseries Log}[modz] + {\bfseries I} theta) - \newline
   modz {\bfseries Exp}[{\bfseries I} theta] + {\bfseries Log}[2 {\bfseries Pi}]/2\newline
strlngline[TP$_{-}$, modz$_{-}$, theta$_{-}$, limit$_{-}$] := {\bfseries Module}[\{e0, e1, e2, e3, e4, e5\}, 
z = modz {\bfseries Exp}[{\bfseries I }theta]; e0 = F[modz, theta]; 
{\bfseries Print}["The value of the leading terms or Stirling's approximation 
is ", {\bfseries N}[e0, 25] // {\bfseries FullForm]}; e1 = z {\bfseries Sum}[Smand[z], {k, 1, TP-1}]; \newline 
{\bfseries Print}["The value of the truncated sum when N equals ", TP, " is ", 
{\bfseries N}[e1, 25] // {\bfseries FullForm]}; Rem = 0; \newline
{\bfseries Do}[sglrty = 2 k {\bfseries Pi} modz; {\bfseries Which}[0 $<$ sglrty $<$ 1, 
e2 = {\bfseries NIntegrate}[Intgrd[TP, z, sglrty], \{y, 0, sglrty, 1, 2, 5, 10, 100, 1000, {\bfseries Infinity} \}, 
{\bfseries Method}-$>$``{\bfseries PrincipalValue}", {\bfseries MinRecursion} -$>$ 3, {\bfseries MaxRecursion} -$>$ 10, 
{\bfseries WorkingPrecision} -$>$ 80, {\bfseries AccuracyGoal} -$>$ 30, {\bfseries PrecisionGoal} -$>$ 30], 1 $<$ sglrty $<$ 2, 
e2 = {\bfseries NIntegrate}[Intgrd[TP, z, sglrty], \{y, 0, 1/2, 1, sglrty, 2, 5, 10, 100, 1000, {\bfseries Infinity} \}, 
{\bfseries Method}-$>$``{\bfseries PrincipalValue}", {\bfseries MinRecursion} -$>$ 3, {\bfseries MaxRecursion} -$>$ 10,  
{\bfseries WorkingPrecision} -$>$ 80, {\bfseries AccuracyGoal} -$>$ 30, {\bfseries PrecisionGoal} -$>$ 30], 2 $<$ sglrty $<$ 5, 
e2 = {\bfseries NIntegrate}[Intgrd[TP, z, sglrty], \{y, 0, 1/2, 1, 2, sglrty, 5, 10, 100, 1000, {\bfseries Infinity} \}, 
{\bfseries Method}-$>$``{\bfseries PrincipalValue}", {\bfseries MinRecursion} -$>$ 3, {\bfseries MaxRecursion} -$>$ 10, 
{\bfseries WorkingPrecision} -$>$ 80, {\bfseries AccuracyGoal} -$>$ 30, {\bfseries PrecisionGoal} -$>$ 30], 5 $<$ sglrty $<$ 10, 
e2 = {\bfseries NIntegrate}[Intgrd[TP, z, sglrty], \{y, 0, 1/2, 1, 2, 5, sglrty, 10, 100, 1000, {\bfseries Infinity} \}, 
{\bfseries Method}-$>$``{\bfseries PrincipalValue}", {\bfseries MinRecursion} -$>$ 3, {\bfseries MaxRecursion} -$>$ 10, 
{\bfseries WorkingPrecision} -$>$ 80, {\bfseries AccuracyGoal} -$>$ 30, {\bfseries PrecisionGoal} -$>$ 30], 10 $<$ sglrty $<$ 100, 
e2 = {\bfseries NIntegrate}[Intgrd[TP, z, sglrty], \{y, 0, 1/2, 1, 2, 5, 10, sglrty, 100, 1000, {\bfseries Infinity}\}, 
{\bfseries Method}-$>$``{\bfseries PrincipalValue}",{\bfseries MinRecursion} -$>$ 3, {\bfseries MaxRecursion} -$>$ 10, 
{\bfseries WorkingPrecision} -$>$ 80, {\bfseries AccuracyGoal} -$>$ 30, {\bfseries PrecisionGoal} -$>$ 30], 100 $<$ sglrty $<$ 1000, 
e2 = {\bfseries NIntegrate}[Intgrd[TP, z, sglrty], \{y, 0, 1/2, 1, 2, 5, 10, 100, sglrty, 1000, {\bfseries Infinity} \}, 
{\bfseries Method}-$>$``{\bfseries PrincipalValue}",{\bfseries MinRecursion} -$>$ 3, {\bfseries MaxRecursion} -$>$ 10, 
{\bfseries WorkingPrecision} -$>$ 80, {\bfseries AccuracyGoal} -$>$ 30, {\bfseries PrecisionGoal} -$>$ 30], sglrty $>$ 1000, 
e2 = {\bfseries NIntegrate}[Intgrd[TP, z, sglrty], \{y, 0, 1/2, 1, 2, 5, 10, 100, 1000, sglrty, {\bfseries Infinity}\}, 
{\bfseries Method}-$>$ ``{\bfseries PrincipalValue}",{\bfseries MinRecursion} -$>$ 3, {\bfseries MaxRecursion} -$>$ 10, 
{\bfseries WorkingPrecision} -$>$ 80, {\bfseries AccuracyGoal} -$>$ 30, {\bfseries PrecisionGoal} -$>$ 30]];\newline 
e2 = e2/k$^{\wedge}$(2 TP - 2);  Rem = Rem + e2, \{k, 1, limit\}]; \newline
Rem = 2  z Rem/(2 {\bfseries Pi} modz)$^{\wedge}$(2 TP - 2) ;\newline 
{\bfseries Print}[" The value of the remainder when N equals ", TP, " is ",  {\bfseries N}[Rem, 25] // {\bfseries FullForm}]; \newline 
{\bfseries Which}[ theta == -{\bfseries Pi}/2, e3 = -{\bfseries Log}[1 - {\bfseries Exp}[-2 {\bfseries Pi} modz]]/2, theta == {\bfseries Pi}/2, 
e3 = -{\bfseries Log}[1 - {\bfseries Exp}[-2 {\bfseries Pi} modz]]/2];\newline 
{\bfseries Print}[" For z=", z, " the value of the Stokes discontinuity term is equal to ", {\bfseries N}[e3, 25] // {\bfseries FullForm}]; 
e4 = e0 + e1 + e3 + Rem; \newline
{\bfseries Print}["Combining all contributions for z=", z, " yields ",  {\bfseries N}[e4, 25] // {\bfseries FullForm}]; 
e5 = {\bfseries LogGamma}[z]; 
{\bfseries Print}["For z=", z, " the value of Log($\Gamma(z)$) is equal to ", {\bfseries N}[e5, 25] // {\bfseries FullForm}]] \vspace{0.5 cm} \newline
The above program uses the same statements for the cosecant numbers, the truncated sum and Stirling's approximation or $F(z)$ 
as in the previous programs. However, the integrand Intgrd is different in accordance with Eq.\ (\ref{twentythree}). In addition,
the Do loop is very different in that there is now a Which statement, which is equivalent to nesting if-then-else
statements in C/C++. It has been introduced here since the singularity in the remainder called sglrty is 
different for each value of k in the Do loop. Because the integral for the remainder has been divided into smaller intervals, 
the Which statement determines the interval where the singularity is situated so that the interval can be split into two 
intervals with the singularity acting as a limit in both of them. It should also be noted that the option Method-$>$``PrincipalValue" 
has been invoked in each call to NIntegrate to avoid convergence problems.   

\subsection*{Program 4}
In Sec.\ 4 we derived via MB regularization general forms of the regularized value of the main asymptotic series $S(z)$ in $\ln \Gamma(z)$.
These forms varied according to the domains of convergence of the MB integrals in these results. In order to verify several of these
general forms numerically over the principal branch of the complex plane, the variable $z$ was replaced by $z^3$, thereby yielding
Eqs.\ (\ref{seventysix})-Eq.\ (\ref{eightythree}). The Mathematica program that evaluates the regularized value of $\ln \Gamma(z^3)$ via
these results appears below: \vspace{0.5 cm} \newline
Intgrd[modz$_{-}$, theta$_{-}$, s$_{-}$, M$_{-}$] := (2 {\bfseries Pi} modz$^{\wedge}$3)$^{\wedge}$(-2 s) {\bfseries Zeta}[2 s] 
{\bfseries Gamma}[2 s - 1] ({\bfseries Exp}[2 (M {\bfseries Pi} - 3 theta) {\bfseries I} s]/(Exp[-{\bfseries I Pi} s] - 
{\bfseries Exp}[{\bfseries I Pi} s])) \newline
c1[k$_{-}$] := -2 {\bfseries Zeta}[2 k]/{\bfseries Pi}$^{\wedge}$(2 k) \newline
F[modz$_{-}$, theta$_{-}$] := (modz {\bfseries Exp}[{\bfseries I} theta] - 1/2) ({\bfseries Log}[modz] + {\bfseries I} theta) - modz 
{\bfseries Exp}[{\bfseries I} theta] + {\bfseries Log}[2 {\bfseries Pi}]/2 \newline
Smand[z$_{-}$] := (-1)$^{\wedge}$k  {\bfseries Gamma}[2 k - 1] c1[k]/(2 z)$^{\wedge}$(2 k) \newline
MBloggam[modz$_{-}$, theta$_{-}$, TP$_{-}$, c$_{-}$] :=  {\bfseries Module}[\{M1, M2\}, s = c + {\bfseries I} r; s1 = c - {\bfseries I} r; 
M1 = 0; M2 = 0; z = modz {\bfseries Exp}[{\bfseries I} theta]; zcube = z$^{\wedge}$3; e0 = F[{\bfseries Abs}[zcube], {\bfseries Arg}[zcube]]; 
{\bfseries Print}["The value of Stirling's approximation is ", {\bfseries N}[e0, 30] // {\bfseries FullForm}];\newline 
e1 = zcube {\bfseries Sum}[Smand[zcube], \{k, 1, TP - 1\}]; {\bfseries Print}["The value of the truncated sum when N equals ", TP, " is ", 
{\bfseries N}[e1, 30] // {\bfseries FullForm}];\newline 
{\bfseries Which}[-{\bfseries Pi} $<$ theta $<$ -2 {\bfseries Pi}/3, M1 = -2; M2 = -3, theta == -2 {\bfseries Pi}/3, M1 = -2; 
{\bfseries Print}["Only one value of M applies."], -2 {\bfseries Pi}/3 $<$ theta $<$ -{\bfseries Pi}/3, M1 = -1; M2 = -2, theta == -{\bfseries Pi}/3, 
M1 = -1; {\bfseries Print}["Only one of value of M applies."], -{\bfseries Pi}/3 $<$ theta $<$ 0, M1 = 0; M2 = -1, theta == 0, M1 = 0; 
{\bfseries Print}["Only one value of M applies."], 0 $<$ theta $<$ {\bfseries Pi}/3, M1 = 0; M2 = 1, theta == {\bfseries Pi}/3, M1 = 1; 
{\bfseries Print}["Only one value of M applies."], {\bfseries Pi}/3 $<$ theta $<$ 2 {\bfseries Pi}/3, M1 = 1; M2 = 2, 
theta == 2 {\bfseries Pi}/3, M1 = 2; {\bfseries Print}["Only one value of M applies."], 
2 {\bfseries Pi}/3 $<$ theta $<$ {\bfseries Pi}, M1 = 2; M2 = 3, theta == {\bfseries Pi}, M1 = 3; {\bfseries Print}["One value of M applies."]]; 
{\bfseries Print}["The values of M1 and M2 are respectively ", M1, " and ", M2];\newline
e2 = {\bfseries NIntegrate}[{\bfseries I} Intgrd[modz, theta, s, M1], \{r, 0, 1/2, 1, 2, 5, 10, 100, 1000, {\bfseries Infinity} \}, 
{\bfseries MinRecursion} -$>$ 3, {\bfseries MaxRecursion} -$>$ 10, {\bfseries WorkingPrecision} -$>$ 60, {\bfseries AccuracyGoal} -$>$ 30, 
{\bfseries PrecisionGoal} -$>$ 30];\newline 
{\bfseries Print}["For M1=", M1, " the value of the upper half of the line contour for the MB integral is ", {\bfseries N}[e2, 30] // 
{\bfseries FullForm}];\newline 
e3 = {\bfseries NIntegrate}[{\bfseries I} Intgrd[modz, theta, s1, M1], \{r, 0, 1/2, 1, 2, 5, 10, 100, 1000, {\bfseries Infinity} \}, 
{\bfseries MinRecursion} -$>$ 3, {\bfseries MaxRecursion} -$>$ 10, {\bfseries WorkingPrecision} -$>$ 60, {\bfseries AccuracyGoal} -$>$ 30, 
{\bfseries PrecisionGoal} -$>$ 30];\newline 
{\bfseries Print}["For M1=", M1, " the value of the lower half of the line contour for the MB integral is ", {\bfseries N}[e3, 30] // 
{\bfseries FullForm}];\newline 
e4 = -2 zcube (e2 + e3);\newline
{\bfseries Print}["For M1=", M1, " the total MB integral including the pre-factor is ",  {\bfseries N}[e4, 30] // {\bfseries FullForm}];\newline 
{\bfseries Which}[-{\bfseries Pi} $<$ theta $<$ -2 {\bfseries Pi}/3 \&\& M1 == -3, logterm = {\bfseries Log}[1 - {\bfseries Exp}[-2 {\bfseries Pi I} zcube]] 
- {\bfseries Log}[- Exp[-2 {\bfseries Pi I} zcube]], -{\bfseries Pi} $<$ theta $<$ -{\bfseries Pi}/3 \&\& M1 == -2, 
logterm = -{\bfseries Log}[- {\bfseries Exp}[-2 {\bfseries Pi I} zcube]], -2 {\bfseries Pi}/3 $<$ theta $<$ 0 \&\& M1 == -1, 
logterm = -{\bfseries Log}[1 - {\bfseries Exp}[-2 {\bfseries Pi I} zcube]], -{\bfseries Pi}/3 $<$ theta $<$ {\bfseries Pi}/3 \&\& M1 == 0, logterm = 0, 
0 $<$ theta $<$ 2 {\bfseries Pi}/3 \&\& M1 == 1, logterm = -{\bfseries Log}[1 - {\bfseries Exp}[2 {\bfseries Pi I} zcube]], 
{\bfseries Pi}/3 $<$ theta $<$ {\bfseries Pi} \&\& M1 == 2, logterm = {\bfseries Log}[-{\bfseries Exp}[-2 {\bfseries Pi I} zcube]], 
2 {\bfseries Pi}/3 $<$ theta $<$= {\bfseries Pi} \&\& M1 == 3, logterm = {\bfseries Log}[-{\bfseries Exp}[-2 {\bfseries Pi I} zcube]] - 
{\bfseries Log}[1 - {\bfseries Exp}[2 {\bfseries Pi I} zcube]]]; \newline 
{\bfseries Print}["For M1 =", M1, " and theta =", theta," the value of the logarithmic term in the regularized value is ", 
{\bfseries N}[logterm, 30] // {\bfseries FullForm}]; tot1 = e0 + e1 + e4 + logterm; \newline 
{\bfseries Print}["For M1=", M1, " and z=", z, " the MB-regularized value of log($\Gamma$(z$^{\wedge}$3)) is ", {\bfseries N}[tot1, 30] 
// {\bfseries FullForm}];\newline 
{\bfseries If}[M2 != 0, e5 = {\bfseries NIntegrate}[{\bfseries I} Intgrd[modz, theta, s, M2], \{r, 0, 1/2, 1, 2, 5, 10, 100, 1000, {\bfseries Infinity} \}, 
{\bfseries MinRecursion} -$>$ 3, {\bfseries MaxRecursion} -$>$ 10, {\bfseries WorkingPrecision} -$>$ 60, {\bfseries AccuracyGoal} -$>$ 30, 
{\bfseries PrecisionGoal} -$>$ 30];
{\bfseries Print}["For M2=", M2," the value of the upper half of the line contour for the MB integral is ", {\bfseries N}[e5, 30] // {\bfseries FullForm}];
\newline e6 = {\bfseries NIntegrate}[{\bfseries I} Intgrd[modz, theta, s1, M2], \{r, 0, 1/2, 1, 2, 5, 10, 100, 1000, {\bfseries Infinity} \}, 
{\bfseries MinRecursion} -$>$ 3, {\bfseries MaxRecursion} -$>$ 10, {\bfseries WorkingPrecision} -$>$ 60, {\bfseries AccuracyGoal} -$>$ 30, 
{\bfseries PrecisionGoal} -$>$ 30];
{\bfseries Print}["For M2=", M2, " the value of the lower half of the line contour for the MB integral is ", {\bfseries N}[e6, 30] // {\bfseries FullForm}]; 
\newline e7 = -2 z$^{\wedge}$3 (e5 + e6); {\bfseries Print}["For M2=", M2, " the total MB integral including the pre-factor is ",  
{\bfseries N}[e7, 30] // {\bfseries FullForm}]; \newline
{\bfseries Which}[-{\bfseries Pi} $<$ theta $<$ -2 {\bfseries Pi}/3 \&\& M2 == -3, logterm2 = -{\bfseries Log}[1 - {\bfseries Exp}[-2 {\bfseries Pi I} zcube]] 
-{\bfseries Log}[- {\bfseries Exp}[-2 {\bfseries Pi I} zcube]], -{\bfseries Pi} $<$ theta $<$ -{\bfseries Pi}/3 \&\& M2 == -2, logterm2 = 
-{\bfseries Log}[- {\bfseries Exp}[-2 {\bfseries Pi I} zcube]], -2 {\bfseries Pi}/3 $<$ theta $<$ 0 \&\& M2 == -1, logterm2 = 
-{\bfseries Log}[1 - {\bfseries Exp}[-2 {\bfseries Pi I} zcube]] , -{\bfseries Pi}/3 $<$ theta $<$ {\bfseries Pi}/3 \&\& M2 == 0, logterm2 = 0, 
0 $<$ theta $<$ 2 {\bfseries Pi}/3 \&\& M2 == 1, logterm2 = -{\bfseries Log}[1 - {\bfseries Exp}[2{\bfseries  Pi I} zcube]], 
{\bfseries Pi}/3 $<$ theta $<$ {\bfseries Pi} \&\& M2 == 2, logterm2 = {\bfseries Log}[-{\bfseries Exp}[-2 {\bfseries Pi I} zcube]],
2 {\bfseries Pi}/3 $<$ theta $<=$ {\bfseries Pi} \&\& M2 == 3, logterm2 = {\bfseries Log}[-{\bfseries Exp}[-2 {\bfseries Pi I} zcube]] 
- {\bfseries Log}[1 - {\bfseries Exp}[2 {\bfseries Pi I} zcube]]];\newline
{\bfseries Print}["For M2 =", M2, " and theta =", theta, " the value of the logarithmic term in the regularized value is 
", {\bfseries N}[logterm2, 30] // {\bfseries FullForm}]; tot2 = e0 + e1 + e7 + logterm2; 
{\bfseries Print}["For M2=", M2, " and z=", z, " the MB-regularized value of log($\Gamma$(z$^{\wedge}$3)) is ", {\bfseries N}[tot2, 30] 
// {\bfseries FullForm}]]; \newline
{\bfseries If}[-{\bfseries Pi}/3 $<$ theta $<=$ {\bfseries Pi}/3, e8 = {\bfseries LogGamma}[zcube]; 
{\bfseries Print}[ "For $|z|$=", modz, " and $\theta$=", theta, " the value of log($\Gamma$(z$^{\wedge}$3)) is ", {\bfseries N}[e8, 30] 
// {\bfseries FullForm}]]]

\subsection*{Program 5}
In Sec.\ 5 it was stated that the higher or lower Borel-summed regularized values of $\ln \Gamma(z)$ appearing in Thm.\ 2.1, viz.\ those
Stokes sectors and lines where $|M|>1$, had not been verified in Sec.\ 3, because they could not be checked against a different result 
other than the routine in Mathematica for LogGamma[z]. Consequently, the numerical analysis in Sec.\ 3 was limited to the principal branch 
of the complex plane for $z$. As a result of MB regularization we now have different forms for the regularized value of $\ln \Gamma(z^3)$,
which can be checked against the Borel-summed regularized values with $z$ replaced by $z^3$. Since the Borel-summed regularized values possess
an infinite convergent sum, we found that the sum had to be truncated at a large number, e.g. $10^5$, in order to achieve accurate results.
Typically, each calculation takes of the order of 6 hours on a Sony VAIO laptop with 2 GB RAM. In order to perform a vast number of these 
calculations simultaneously, it is best to carry them out on a server with as many processors as possible. Appearing below is the script 
program that was run on a SunFire alpha server to produce the results in Table\ \ref{tab9}. This program is a composite of Programs 3 and 4 
except that the Which statement in the former program has been expanded to include more Stokes sectors as a result of altering $z$ to $z^3$. 
The variable logterm3 is responsible for obtaining the appropriate logarithmic terms in the Borel-summed values given in Eq.\ (\ref{ninety}).
\vspace{0.5cm} \newline
{\bfseries Print}[\$ScriptCommandLine] \newline
params={\bfseries Rest@ToExpression}[\$ScriptCommandLine] \newline
{\bfseries Print}[params]\vspace{0.5 cm}\newline
Intgrd[modz$_{-}$, theta$_{-}$, s$_{-}$,M$_{-}$] := (2 {\bfseries Pi} modz$^{\wedge}$3)$^{\wedge}$(-2 s) {\bfseries Zeta}[2 s] 
{\bfseries Gamma}[2 s - 1] ({\bfseries Exp}[2 (M {\bfseries Pi} - 3 theta) {\bfseries I} s]/({\bfseries Exp}[-{\bfseries I Pi} s] - 
{\bfseries Exp}[{\bfseries I Pi} s])) \newline
Intgrd2[N$_{-}$, z$_{-}$, k$_{-}$] := y$^{\wedge}$(2 N - 2) {\bfseries Exp}[-y]/(y$^{\wedge}$2 + (2 {\bfseries Pi} k z)$^{\wedge}$2) \newline
c1[k$_{-}$] := -2 {\bfseries Zeta}[2 k]/{\bfseries Pi}$^{\wedge}$(2 k)\newline
F[modz$_{-}$, theta$_{-}$] := (modz {\bfseries Exp}[{\bfseries I} theta] - 1/2) ({\bfseries Log}[modz] + {\bfseries I} theta) - modz 
{\bfseries Exp}[{\bfseries I} theta] + {\bfseries Log}[2 {\bfseries Pi}]/2
Smand[z$_{-}$] := (-1)$^{\wedge}$k  {\bfseries Gamma}[2 k - 1] c1[k]/(2 z)$^{\wedge}$(2 k) \newline
MBBrlStr[modz$_{-}$, theta$_{-}$, TP$_{-}$, limit$_{-}$] := {\bfseries Module}[\{M1, M2\}, c = TP - 1/4; s = c + {\bfseries I} r; 
s1 = c - {\bfseries I} r; M1 = 0; M2 = 0; z = modz {\bfseries Exp}[{\bfseries I} theta]; zcube = z$^{\wedge}$3; 
e0 = F[{\bfseries Abs}[zcube], {\bfseries Arg}[zcube]]; \newline
{\bfseries Print}["The value of Stirling's approximation is ", {\bfseries N}[e0, 30] // {\bfseries FullForm}]; 
e1 = zcube {\bfseries Sum}[Smand[zcube], \{k, 1, TP - 1\}]; \newline
{\bfseries Print}["The value of the truncated sum when N equals ", TP, " is ", {\bfseries N}[e1, 30] // {\bfseries FullForm}];\newline 
{\bfseries Which}[-{\bfseries Pi} $<$ theta $<$ -2 {\bfseries Pi}/3, M1 = -2; M2 = -3, theta == -2 {\bfseries Pi}/3, M1 = -2; 
{\bfseries Print}["Only one value of M applies."], -2 {\bfseries Pi}/3 $<$ theta $<$ -{\bfseries Pi}/3, 
 M1 = -1; M2 = -2, theta == -{\bfseries Pi}/3, M1 = -1; {\bfseries Print}["Only one of value of M applies."], -{\bfseries Pi}/3 $<$ theta $<$ 0, 
M1 = 0; M2 = -1, theta == 0, M1 = 0; {\bfseries Print}["Only one value of M applies."], 0 $<$ theta $<$ {\bfseries Pi}/3, M1 = 0; 
 M2 = 1, theta == {\bfseries Pi}/3, M1 = 1; {\bfseries Print}["Only one value of M applies."], {\bfseries Pi}/3 $<$ theta $<$ 2 {\bfseries Pi}/3, 
M1 = 1; M2 = 2, theta == 2 {\bfseries Pi}/3, M1 = 2; {\bfseries  Print}["Only one value of M applies."], 2 {\bfseries Pi}/3 $<$ theta $<$ {\bfseries Pi}, 
M1 = 2;M2 = 3, theta == {\bfseries Pi}, M1 = 3; {\bfseries Print}["One value of M applies."]];\newline 
{\bfseries Print}["The values of M1 and M2 are respectively ", M1, " and ", M2];\newline
 e2 = {\bfseries NIntegrate}[{\bfseries I} Intgrd[modz, theta, s, M1], \{r, 0, 1/2, 1, 2, 5, 10, 100, 1000, {\bfseries Infinity} \}, 
{\bfseries MinRecursion} -$>$ 3, {\bfseries MaxRecursion} -$>$ 10, {\bfseries WorkingPrecision} -$>$ 60, {\bfseries AccuracyGoal} -$>$ 30, 
{\bfseries PrecisionGoal} -$>$ 30]; 
{\bfseries  Print}["For M1=", M1," the value of the upper half of the line contour for the MB integral is ", {\bfseries N}[e2, 30] // {\bfseries FullForm}]; 
e3 = {\bfseries NIntegrate}[{\bfseries I} Intgrd[modz, theta, s1, M1], \{r, 0, 1/2, 1, 2, 5, 10, 100, 1000, {\bfseries Infinity} \}, 
{\bfseries MinRecursion} -$>$ 3, {\bfseries MaxRecursion} -$>$ 10, {\bfseries WorkingPrecision} -$>$ 60, {\bfseries AccuracyGoal} -$>$ 30, 
{\bfseries PrecisionGoal} -$>$ 30]; \newline
{\bfseries  Print}["For M1=", M1," the value of the lower half of the line contour for the MB integral is ", {\bfseries N}[e3, 30] // {\bfseries FullForm}]; 
 e4 = -2 zcube (e2 + e3); \newline 
{\bfseries Print}["For M1=", M1, " the total MB integral including the pre-factor is ", {\bfseries N}[e4, 30] // {\bfseries FullForm}];\newline 
{\bfseries  Which}[-{\bfseries Pi} $<$ theta $<$ -2 {\bfseries Pi}/3 \&\& M1 == -3, logterm = {\bfseries Log}[1 - {\bfseries Exp}[-2 {\bfseries Pi I} zcube]] 
-{\bfseries Log}[- {\bfseries Exp}[-2 {\bfseries Pi I} zcube]], -{\bfseries Pi} $<$ theta $<$ -{\bfseries Pi}/3 \&\& M1 == -2,logterm = 
-{\bfseries Log}[- {\bfseries Exp}[-2 {\bfseries Pi I} zcube]], -2 {\bfseries Pi}/3 $<$ theta $<$ 0 \&\& M1 == -1, logterm = 
-{\bfseries Log}[1 - {\bfseries Exp}[-2 {\bfseries Pi I} zcube]], -{\bfseries Pi}/3 $<$ theta $<$ {\bfseries Pi}/3 \&\& M1 == 0, logterm = 0, 
0 $<$ theta $<$ 2 {\bfseries Pi}/3 \&\& M1 == 1,logterm = -{\bfseries Log}[1 - {\bfseries Exp}[2 {\bfseries Pi I} zcube]], 
{\bfseries Pi}/3 $<$ theta $<$ {\bfseries Pi} \&\& M1 == 2, logterm = {\bfseries Log}[-{\bfseries Exp}[-2 {\bfseries Pi I} zcube]], 
2 {\bfseries Pi}/3 $<$ theta $<$= {\bfseries Pi} \&\& M1 == 3, logterm = {\bfseries Log}[-Exp[-2 {\bfseries Pi I} zcube]] - 
{\bfseries  Log}[1 - {\bfseries Exp}[2 {\bfseries Pi I} zcube]]]; \newline 
{\bfseries  Print}["For M1 =", M1, " and theta =", theta," the value of the logarithmic term in the regularized value is ", 
{\bfseries N}[logterm, 30] // {\bfseries FullForm}]; tot1 = e0 + e1 + e4 + logterm; \newline
{\bfseries Print}["For M1=", M1, " and z=", z," the MB-regularized value of log($\Gamma$](z$^{\wedge}$3)) is ", {\bfseries N}[tot1, 30] 
// {\bfseries FullForm}]; \newline
{\bfseries  If}[M2 != 0, e5 = {\bfseries NIntegrate}[{\bfseries I} Intgrd[modz, theta, s, M2], \{r, 0, 1/2, 1, 2, 5, 10, 100, 1000, {\bfseries Infinity} \}, 
{\bfseries MinRecursion} -$>$ 3, {\bfseries MaxRecursion} -$>$ 10, {\bfseries WorkingPrecision} -$>$ 60, {\bfseries AccuracyGoal} -$>$ 30, 
{\bfseries PrecisionGoal} -$>$ 30];
{\bfseries Print}["For M2=", M2," the value of the upper half of the line contour for the MB integral is ", {\bfseries N}[e5, 30] // {\bfseries FullForm}]; 
 e6 = NIntegrate[I Intgrd[modz, theta, s1, M2], \{r, 0, 1/2, 1, 2, 5, 10, 100, 1000, {\bfseries Infinity} \}, 
{\bfseries MinRecursion} -$>$ 3, {\bfseries MaxRecursion} -$>$ 10, {\bfseries WorkingPrecision} -$>$ 60, {\bfseries AccuracyGoal} -$>$ 30, 
{\bfseries PrecisionGoal} -$>$ 30]; \newline
{\bfseries  Print}["For M2=", M2, " the value of the lower half of the line contour for the MB integral is ", {\bfseries N}[e6, 30] // {\bfseries FullForm}]; 
 e7 = -2 z$^{\wedge}$3 (e5 + e6); \newline
 {\bfseries Print}["For M2=", M2," the total MB integral including the pre-factor is ", {\bfseries N}[e7, 30] // {\bfseries FullForm}];\newline 
 {\bfseries Which}[-{\bfseries Pi} $<$ theta $<$ -2 {\bfseries Pi}/3 \&\& M2 == -3, logterm2 = -{\bfseries Log}[1 - {\bfseries Exp}[-2 {\bfseries Pi I} zcube]] 
-{\bfseries Log}[- {\bfseries Exp}[-2 {\bfseries Pi I} zcube]], -{\bfseries Pi} $<$ theta $<$ -{\bfseries Pi}/3 \&\& M2 == -2, logterm2 = 
-{\bfseries Log}[- {\bfseries Exp}[-2 {\bfseries Pi I} zcube]], -2 {\bfseries Pi}/3 $<$ theta $<$ 0 \&\& M2 == -1, logterm2 = 
-{\bfseries Log}[1 - {\bfseries Exp}[-2 {\bfseries Pi I} zcube]], -{\bfseries Pi}/3 $<$ theta $<$ {\bfseries Pi}/3 \&\&  M2 == 0, logterm2 = 0, 
0 $<$ theta $<$ 2 {\bfseries Pi}/3 \&\& M2 == 1,logterm2 = -{\bfseries Log}[1 - {\bfseries Exp}[2 {\bfseries Pi I} zcube]], 
{\bfseries  Pi}/3 $<$ theta $<$ {\bfseries Pi} \&\& M2 == 2, logterm2 = {\bfseries Log}[-{\bfseries Exp}[-2 {\bfseries Pi I} zcube]],
2 {\bfseries Pi}/3 $<$ theta $<$= {\bfseries Pi} \&\& M2 == 3, logterm2 = {\bfseries Log}[-{\bfseries Exp}[-2 {\bfseries Pi I} zcube]] 
- {\bfseries Log}[1 - {\bfseries Exp}[2 {\bfseries Pi I} zcube]]];\newline
{\bfseries Print}["For M2 =", M2, " and theta =", theta," the value of the logarithmic term in the regularized value is ", {\bfseries N}[logterm2, 30] 
// {\bfseries FullForm}]; tot2 = e0 + e1 + e7 + logterm2; \newline
{\bfseries  Print}["For M2=", M2, " and z=", z," the MB-regularized value of log($\Gamma$(z$^{\wedge}$3)) is ", {\bfseries N}[tot2, 30] // 
{\bfseries FullForm}]]; Rem = 0; \newline 
{\bfseries Do}[e8 = {\bfseries NIntegrate}[Intgrd2[TP, zcube, k], \{y, 0, 1/2, 1, 2, 5, 10, 100, 1000, {\bfseries Infinity} \}, 
{\bfseries MinRecursion} -$>$ 3, {\bfseries MaxRecursion} -$>$ 10, {\bfseries WorkingPrecision} -$>$ 60, {\bfseries AccuracyGoal} -$>$ 30, 
{\bfseries PrecisionGoal} -$>$ 30]; e8 = e8/k$^{\wedge}$(2 TP - 2); Rem = Rem + e8, \{k, 1, limit\}]; 
 Rem = 2 (-1)$^{\wedge}$(TP + 1) zcube Rem/(2 Pi zcube)$^{\wedge}$(2 TP - 2); \newline
{\bfseries  Print}["The value of the Borel-summed remainder when N equals ", TP," is ", {\bfseries N}[Rem, 25] // {\bfseries FullForm}]; \newline 
{\bfseries } Which[-{\bfseries Pi} $<$ theta $<$ -5 {\bfseries Pi}/6,logterm3 = -{\bfseries Log}[- {\bfseries Exp}[-2 {\bfseries Pi I} zcube]] 
- {\bfseries Log}[1 - {\bfseries Exp}[-2 {\bfseries Pi I} zcube]], -5 {\bfseries Pi}/6 $<$ theta $<$ -{\bfseries Pi}/2, 
 logterm3 = -{\bfseries Log}[- {\bfseries Exp}[-2 {\bfseries Pi I} zcube]], -{\bfseries Pi}/2 $<$ theta $<$ -{\bfseries Pi}/6, logterm3 = 
-{\bfseries Log}[1 - {\bfseries Exp}[-2 {\bfseries Pi I} zcube]] , -{\bfseries Pi}/6 $<$ theta $<$ {\bfseries Pi}/6,logterm3 = 0, 
{\bfseries Pi}/6 $<$ theta $<$ {\bfseries Pi}/2, logterm3 = -{\bfseries Log}[1 - {\bfseries Exp}[2 {\bfseries Pi I} zcube]], {\bfseries Pi}/2 
$<$ theta $<$ 5 {\bfseries Pi}/6, logterm3 = -{\bfseries Log}[-{\bfseries Exp}[2 {\bfseries Pi I} zcube]], 5 {\bfseries Pi}/6 $<$ theta 
$<$= {\bfseries Pi},logterm3 = -{\bfseries Log}[-{\bfseries Exp}[2 {\bfseries Pi I} zcube]] - {\bfseries Log}[1 - 
{\bfseries Exp}[2 {\bfseries Pi I} zcube]]]; \newline
{\bfseries Print}["For $|{\rm z}|$ =", modz, " and Theta=", theta," the value of the Borel log term is ", {\bfseries N}[logterm3, 30] 
// {\bfseries FullForm}]; 
e9 = e0 + e1 + logterm3 + Rem; \newline 
{\bfseries Print}["For $|{\rm z}|$ =", modz, " and Theta=", theta," the Borel-summed value of log($\Gamma$(z$^{\wedge}$3)) is ", {\bfseries N}[e9, 30] 
// {\bfseries FullForm}];
\newline 
{\bfseries If}[-{\bfseries Pi}/3 $<$ theta $<$= {\bfseries Pi}/3, e10 = {\bfseries LogGamma}[zcube]; 
{\bfseries Print}[ "For $|{\rm z}|$=", modz, " and Theta=", theta, " the value of log($\Gamma$(z$^{\wedge}$3)) is ",  {\bfseries N}[e10, 30] 
// {\bfseries FullForm}]]] \vspace{0.5cm}\newline
{\bfseries Timing}[MBBrlStr@@params]

\subsection*{Program 6}
The final program that has been used in this work is a Mathematica batch program that is designed to evaluate 
the MB-regularized asymptotic forms of $\ln \Gamma(z^3)$ for $\theta$ lying on the Stokes lines, viz.\ Eqs.\ 
(\ref{seventyseven}), (\ref{seventysevenb}), (\ref{seventyninea}), (\ref{eighty}), (\ref{eightyone}), 
(\ref{eightytwob})-(\ref{eightytwod}) and (\ref{eightythree}) and then to compare them with the corresponding 
Borel-summed asymptotic forms in Eq.\ (\ref{ninetyone}). For the special case, where $| \theta|=\pi/6$, the code 
evaluates $\ln \Gamma(z^3)$ via the LogGamma routine in Mathematica. The code is able to consider the various 
Stokes lines as a result of {\bfseries Which} statements that allow the Stokes discontinuity terms in the 
Borel-summed asymptotic forms and $S_{MB}(M,z^3)$ in the MB-regularized asymptotic forms to be calculated for 
any value of $\theta$ and $M$.
\vspace{0.5 cm} \newline
{\bfseries Print}[\$ScriptCommandLine] \newline
params={\bfseries Rest@ToExpression}[\$ScriptCommandLine] \newline
{\bfseries Print}[params]\vspace{0.5cm} \newline
Intgrd[modz$_{-}$, theta$_{-}$, s$_{-}$,M$_{-}$] := (2 {\bfseries Pi} modz$^{\wedge}$3)$^{\wedge}$(-2 s) {\bfseries Zeta}[2 s] {\bfseries Gamma}
[2 s - 1] \newline{\bfseries Exp}[2 (M {\bfseries Pi} - 3 theta) {\bfseries I} s]/({\bfseries Exp}[-{\bfseries I Pi} s] - 
{\bfseries Exp}[{\bfseries I Pi} s])) \newline
Intgrd3[N$_{-}$,snglrty$_{-}$]:=y$^{\wedge}$(2 N-2) {\bfseries Exp}[-y]/(y$^{\wedge}$2-snglrty$^{\wedge}$2)\newline
c1[k$_{-}$]:=-2 {\bfseries Zeta}[2 k]/{\bfseries Pi}$^{\wedge}$(2 k)\newline
F[modz$_{-}$,theta$_{-}$]:=(modz {\bfseries Exp}[{\bfseries I} theta]-1/2) ({\bfseries Log}[modz]+ {\bfseries I} theta) -modz 
{\bfseries Exp}[{\bfseries I} theta]+ {\bfseries Log}[2 {\bfseries Pi}]/2 \newline
Smand[z$_{-}$]:=(-1)$^{\wedge}$k  {\bfseries Gamma}[2 k-1]c1[k]/(2 z)$^{\wedge}$(2k)\newline
StrStksLn[modz$_{-}$,theta$_{-}$,TP$_{-}$,limit$_{-}$]:={\bfseries Module}[\{M1,M2,e8,Rem,sglrty\},c=TP-1/4;M1=0;\newline M2=0;
s=c+{\bfseries I} r;s1=c-{\bfseries I} r;M1=0;M2=0;z=modz {\bfseries Exp[}I theta];zcube=z$^{\wedge}$3;\newline mdzcube=modz$^{\wedge}$3;
e0=F[{\bfseries Abs}[zcube],{\bfseries Arg}[zcube]]; \newline {\bfseries Print}["The value of Stirling's approximation is ",{\bfseries N}
[e0,30]//{\bfseries FullForm}];\newline e1=zcube {\bfseries Sum}[Smand[zcube],\{k,1,TP-1\}];
{\bfseries Print}["The value of the truncated sum when N equals ",TP," is ",{\bfseries N}[e1,30]//{\bfseries FullForm}];\newline
{\bfseries Which}[theta==-5 {\bfseries Pi}/6,M1=-2;M2=-3,theta==-{\bfseries Pi}/2,M1=-1;M2=-2,
theta==-{\bfseries Pi}/6,M1=0;M2=-1,theta=={\bfseries Pi}/6,M1=0;M2=1,theta=={\bfseries Pi}/2,M1=1;M2=2,theta==5 {\bfseries Pi}/6,M1=2;M2=3];\newline 
{\bfseries Print}["The values of M1 and M2 are respectively ",M1," and ",M2]; \newline
e2={\bfseries NIntegrate}[I Intgrd[modz,theta,s,M1],\{r,0,1/2,1,2,5,10,100,1000,{\bfseries Infinity}\},{\bfseries MinRecursion}-$>$3,
{\bfseries MaxRecursion}-$>$10,{\bfseries WorkingPrecision}-$>$60, {\bfseries AccuracyGoal}-$>$30,{\bfseries PrecisionGoal}-$>$30];
{\bfseries Print}["For M1=",M1," the value of the upper half of the line contour for the MB integral is ",{\bfseries N[}e2,30]//{\bfseries FullForm}];
\newline  e3={\bfseries NIntegrate}[I Intgrd[modz,theta,s1,M1],\{r,0,1/2,1,2,5,10,100,1000,{\bfseries Infinity}\},{\bfseries MinRecursion}-$>$3,
{\bfseries MaxRecursion}-$>$10,{\bfseries WorkingPrecision}-$>$60, {\bfseries AccuracyGoal}-$>$30,{\bfseries PrecisionGoal}-$>$30];
{\bfseries Print}["For M1=",M1," the value of the lower half of the line contour for the MB integral is ",{\bfseries N}[e3,30]//{\bfseries FullForm}]; 
\newline e4=-2 zcube(e2+e3); {\bfseries Print}["For M1=",M1," the total MB integral including the pre-factor is ", {\bfseries N}[e4,30]//{\bfseries FullForm}];
\newline {\bfseries Which}[{\bfseries Abs}[M1]==3 \&\& {\bfseries Abs}[theta]== 5{\bfseries Pi}/6 ,exterm=-{\bfseries Log}[1- {\bfseries Exp}[-2 {\bfseries Pi} mdzcube]]
+2 {\bfseries Pi} mdzcube,{\bfseries Abs}[M1]==2 \&\& {\bfseries Abs}[theta] ==5 {\bfseries Pi}/6,exterm=2 {\bfseries Pi} mdzcube,\newline {\bfseries Abs}[M1]==2 \&\& 
{\bfseries Abs}[theta] == {\bfseries Pi}/2,exterm=-2 {\bfseries Pi} mdzcube, {\bfseries Abs}[M1]==1 \&\& {\bfseries Abs}[theta] == {\bfseries Pi}/2,
exterm=-{\bfseries Log}[1- {\bfseries Exp}[-2 {\bfseries Pi} mdzcube]] - 2 {\bfseries Pi} mdzcube,{\bfseries Abs}[M1]==1 \&\& {\bfseries Abs}[theta] == 
{\bfseries Pi}/6,exterm= -{\bfseries Log}[1- {\bfseries Exp}[-2 {\bfseries Pi} mdzcube]],M1==0,exterm=0]; \newline 
{\bfseries Print}["For M1 =",M1," and theta =",theta," the value of the logarithmic term in the regularized value is ", {\bfseries N}
[exterm,30]//{\bfseries FullForm}];\newline tot1=e0+e1+e4+exterm; \newline
{\bfseries Print}["For M1=",M1," and z=",z," the MB-regularized value of log(Gamma(z$^{\wedge}$3)) is ",
{\bfseries N}[tot1,30]//{\bfseries FullForm}];\newline
{\bfseries If}[M2!=0,e5={\bfseries NIntegrate}[{\bfseries I} Intgrd[modz,theta,s,M2],\{r,0,1/2,1,2,5,10,100,1000,Infinity\},\newline 
{\bfseries MinRecursion}-$>$3,{\bfseries MaxRecursion}-$>$10,{\bfseries WorkingPrecision}-$>$60, {\bfseries AccuracyGoal}-$>$30,\newline
{\bfseries PrecisionGoal}-$>$30];
{\bfseries Print}["For M2=",M2," the value of the upper half of the line contour for the MB integral is ",{\bfseries N}[e5,30]//{\bfseries FullForm}];
\newline e6={\bfseries NIntegrate}[I Intgrd[modz,theta,s1,M2],{r,0,1/2,1,2,5,10,100,1000,{\bfseries Infinity}},{\bfseries MinRecursion}-$>$3,
{\bfseries MaxRecursion}-$>$10,{\bfseries WorkingPrecision}-$>$60, {\bfseries AccuracyGoal}-$>$30,{\bfseries PrecisionGoal}-$>$30];
{\bfseries Print}["For M2=",M2," the value of the lower half of the line contour for the MB integral is ",{\bfseries N}[e6,30]//{\bfseries FullForm}];
e7=-2 z$^{\wedge}$3(e5+e6); \newline {\bfseries Print}["For M2=",M2," the total MB integral including the pre-factor is ", {\bfseries N}[e7,30]//{\bfseries FullForm}];
\newline {\bfseries Which}[{\bfseries Abs}[M2]==3 \&\& {\bfseries Abs}[theta]== 5{\bfseries Pi}/6 ,exterm2=-{\bfseries Log}[1- {\bfseries Exp}[-2 {\bfseries Pi} mdzcube]]
+2 {\bfseries Pi} mdzcube,{\bfseries Abs}[M2]==2 \&\& {\bfseries Abs}[theta] ==5 {\bfseries Pi}/6,exterm2=2 {\bfseries Pi} mdzcube,\newline {\bfseries Abs}[M2]==2 \&\& 
{\bfseries Abs}[theta] == {\bfseries Pi}/2,exterm2=-2 {\bfseries Pi} mdzcube, {\bfseries Abs}[M2]==1 \&\& {\bfseries Abs}[theta] == {\bfseries Pi}/2,
exterm2=-{\bfseries Log}[1- {\bfseries Exp}[-2 {\bfseries Pi} mdzcube]] - 2 {\bfseries Pi} mdzcube,{\bfseries Abs}[M2]==1 \&\& {\bfseries Abs}[theta] == 
{\bfseries Pi}/6,exterm2= -{\bfseries Log}[1- {\bfseries Exp}[-2 {\bfseries Pi} mdzcube]],M2==0,exterm2=0]; \newline 
{\bfseries Print}["For M2 =",M2," and theta =",theta," the value of the logarithmic term in the regularized value is ", 
{\bfseries N}[exterm2,30]//{\bfseries FullForm}];tot2=e0+e1+e7+exterm2;\newline {\bfseries Print}["For M2=",M2," and z=",z," the MB-regularized value of 
log(Gamma(z$^{\wedge}$3)) is ",{\bfseries N}[tot2,30]//{\bfseries FullForm}]];
Rem=0;\newline {\bfseries Do}[sglrty=2 k {\bfseries Pi} mdzcube;{\bfseries Which}[0$<$sglrty$<$1,\newline e8={\bfseries NIntegrate}[Intgrd3[TP,sglrty],
\{y,0,sglrty,1,2,5,10,100,1000,{\bfseries Infinity}\},{\bfseries Method}-$>$\newline"{\bfseries PrincipalValue}",{\bfseries MinRecursion}-$>$3,
{\bfseries MaxRecursion}-$>$10,{\bfseries WorkingPrecision}-$>$60,\newline {\bfseries AccuracyGoal}-$>$30,{\bfseries PrecisionGoal}-$>$30],1$<$ sglrty$<$2,
\newline e8={\bfseries NIntegrate}[Intgrd3[TP,sglrty],\{y,0,1/2,1,sglrty,2,5,10,100,1000,{\bfseries Infinity}\},{\bfseries Method}-$>$\newline "{\bfseries PrincipalValue}",
{\bfseries MinRecursion}-$>$3,{\bfseries MaxRecursion}-$>$10,{\bfseries WorkingPrecision}-$>$60,\newline {\bfseries AccuracyGoal}-$>$30,{\bfseries PrecisionGoal}-$>$30],
2$<$ sglrty$<$5, \newline e8={\bfseries NIntegrate}[Intgrd3[TP,sglrty],\{y,0,1/2,1,2,sglrty,5,10,100,1000,{\bfseries Infinity}\},{\bfseries Method}-$>$
\newline "{\bfseries PrincipalValue}",{\bfseries MinRecursion}-$>$3,{\bfseries MaxRecursion}-$>$10,{\bfseries WorkingPrecision}-$>$60,\newline 
{\bfseries AccuracyGoal}-$>$30,{\bfseries PrecisionGoal}-$>$30],5$<$ sglrty$<$10,
\newline e8=NIntegrate[Intgrd3[TP,sglrty],\{y,0,1/2,1,2,5,sglrty,10,100,1000,{\bfseries Infinity}\},{\bfseries Method}-$>$\newline"{\bfseries PrincipalValue}",
{\bfseries MinRecursion}-$>$3,{\bfseries MaxRecursion}-$>$10,{\bfseries WorkingPrecision}-$>$60,\newline{\bfseries AccuracyGoal}-$>$30,{\bfseries PrecisionGoal}-$>$30],
10$<$ sglrty$<$100,
\newline e8={\bfseries NIntegrate}[Intgrd3[TP,sglrty],\{y,0,1/2,1,2,5,10,sglrty,100,1000,{\bfseries Infinity}\},{\bfseries Method}-$>$\newline "{\bfseries PrincipalValue}",
{\bfseries MinRecursion}-$>$3,{\bfseries MaxRecursion}-$>$10,{\bfseries WorkingPrecision}-$>$60,\newline{\bfseries AccuracyGoal}-$>$30,{\bfseries PrecisionGoal}-$>$30],
100$<$ sglrty$<$1000,\newline e8={\bfseries NIntegrate}[Intgrd3[TP,sglrty],\{y,0,1/2,1,2,5,10,100,sglrty,1000,{\bfseries Infinity}\},{\bfseries Method}-$>$
\newline "PrincipalValue",{\bfseries MinRecursion}-$>$3,{\bfseries MaxRecursion}-$>$10,{\bfseries WorkingPrecision}-$>$60,
\newline {\bfseries AccuracyGoal}-$>$30,{\bfseries PrecisionGoal}-$>$30],sglrty$>$1000,\newline e8={\bfseries NIntegrate}[Intgrd3[TP,sglrty],\{y,0,1/2,1,2,5,
10,100,1000,sglrty,{\bfseries Infinity}\},{\bfseries Method}-$>$"{\bfseries PrincipalValue}",
{\bfseries MinRecursion}-$>$3,{\bfseries MaxRecursion}-$>$10,{\bfseries WorkingPrecision}-$>$60,\newline {\bfseries AccuracyGoal}-$>$30,{\bfseries PrecisionGoal}-$>$30]];
\newline e8=e8/k$^{\wedge}$(2 TP-2);Rem=Rem+e8,\{k,1,limit\}];Rem= 2 zcube Rem/(2 {\bfseries Pi} mdzcube)$^{\wedge}$(2 TP-2);
\newline {\bfseries Print}[" The value for the remainder when N equals ",TP," is ",
{\bfseries N}[Rem,25]//{\bfseries FullForm}];
{\bfseries Which}[{\bfseries Abs}[theta]=={\bfseries Pi}/6,e9=-(1/2) {\bfseries Log}[1- {\bfseries Exp}[-2 {\bfseries Pi} mdzcube]],
{\bfseries Abs}[theta]=={\bfseries Pi}/2,\newline e9=-(1/2){\bfseries Log}[1- {\bfseries Exp}[-2 {\bfseries Pi} mdzcube]]- 2 {\bfseries Pi} mdzcube,
{\bfseries Abs}[theta]==5 {\bfseries Pi}/6,e9=(-1/2)\newline {\bfseries Log}[1- {\bfseries Exp}[-2 {\bfseries Pi} mdzcube]]+2 {\bfseries Pi} mdzcube];
\newline {\bfseries Print}["For z=",z," the value of the Stokes discontinuity term is equal to ",{\bfseries N}[e9,25]//{\bfseries FullForm}];
\newline e10=e0+e1+e9+Rem;\newline 
{\bfseries Print}["Combining all contributions for z=",z," yields ",{\bfseries N}[e10,25]//{\bfseries FullForm}]; 
\newline {\bfseries If}[{\bfseries Abs}[theta] $<$= {\bfseries Pi}/3,e11={\bfseries LogGamma}[zcube]; \newline
{\bfseries Print}["For $|$z$|$=",modz," and theta =",theta," the value of log(Gamma(z$^{\wedge}$3)) via LogGamma[z] is ",
{\bfseries N[}e11,30]//{\bfseries FullForm}]]]\vspace{0.5cm} \newline
{\bfseries Timing}[StrStksLn@@params]

\end{document}